

\documentstyle{amsppt}
\magnification1200 
\NoBlackBoxes
\def\phi{\varphi}
\def\RE{\text{\rm Re }}
\def\IM{\text{\rm Im }}
\def\GT{\Gamma_{\Theta}(S)}

\def\hbar{\overline{h}}

\def\Lam{\Lambda}
\def\hchi{\hat \chi}
\def\hsigma{\hat \sigma}
\def\hbar{\overline{h}}
\pageheight{9 true in}
\pagewidth{6.5 true in}

\topmatter
\title
The Spectrum of Multiplicative Functions
\endtitle
\author
Andrew Granville and K. Soundararajan
\endauthor
\dedicatory{Dedicated to Richard Guy on his 80th birthday, for all
the inspiring problems that he has posed}
\enddedicatory
\rightheadtext{The Spectrum of Multiplicative Functions}

\address
Department of Mathematics, University of Georgia, Athens, GA , USA
\endaddress
\email 
andrew\@sophie.math.uga.edu
\endemail
\address
School of Mathematics,  Institute for Advanced Study, 
Princeton, NJ 08540, USA
\endaddress 
\email 
ksound\@math.ias.edu
\endemail
\thanks{The first author is a Presidential Faculty Fellow.  He is also
supported, in part, by the National Science Foundation.  The second 
author is 
supported by the American Institute of Mathematics (AIM), 
and, in part, by the National Science Foundation (DMS 97-29992).}\endthanks

\toc

\head 1  Introduction:  Definitions and properties of the spectrum \endhead

\head 2  The natural and logarithmic densities of $m$-th power residues 
\endhead

\subhead 2a  The proportion of $m$-th power residues up to $x$ 
\endsubhead 

\subhead 2b  Logarithmic proportions of $m$-th power residues \endsubhead
 
\head 3  Basic properties of integral equations \endhead

\subhead 3a  Existence and uniqueness of solutions and first estimates
\endsubhead

\subhead 3b  Inclusion-Exclusion inequalities
\endsubhead

\head 4  Proof of the Structure Theorem \endhead

\subhead 4a  Variation of averages of multiplicative functions \endsubhead

\subhead 4b  Some useful identities \endsubhead

\subhead 4c  Removing the impact of small primes \endsubhead

\subhead 4d  Completing the proof of the Structure Theorem  \endsubhead
 
\head 5  The spectrum of $[-1,1]$ \endhead

\subhead 5a Preliminaries \endsubhead

\subhead 5b Bounding $\int_0^{u_0} |\sigma(u)| du $ \endsubhead

\subhead 5c Proof of Theorem 5.1 for large $\tau (>29/100)$ \endsubhead

\subhead 5d The range $u_0 \le u \le 2u_0/\sqrt{e}$ \endsubhead

\subhead 5e The range $2u_0/\sqrt{e} \le u\le (1+1/\sqrt{e}) u_0$ \endsubhead

\subhead 5f Completing the proof of Theorem 5.1 \endsubhead

\head 6  The Euler product spectrum \endhead

\subhead 6a Proof of Theorem 4 \endsubhead

\subhead 6b Proof of Corollary 3 \endsubhead

\subhead 6c Proof of Theorem 3$'$ \endsubhead

\head 7  Angles and projections of the spectrum \endhead

\subhead 7a Proof that Ang$(\Gamma(S)) \ll \text{Ang}(S)$:\ Theorem 6(i) \endsubhead

\subhead 7b The maximal projection of $S= \{\pm 1, \pm i\}$:\ Theorem 7(i)
\endsubhead

\subhead 7c Towards the proofs of Theorems 5, 6(ii), and 7(ii) \endsubhead

\subhead 7d Proof of Theorem 5 \endsubhead 

\subhead 7e Proof of Theorem 6(ii) \endsubhead 

\subhead 7f Proof of Theorem 7(ii) \endsubhead

\head 8  Generalized notions of the spectrum: the logarithmic spectrum 
\endhead
        
\subhead 8a Generalities on logarithmic means \endsubhead

\subhead 8b  Bounding $\Gamma_0(S)$:  Proof of Theorem 8  \endsubhead

\subhead 8c  Proof of Corollary 4 \endsubhead

\head 9  Quadratic Residues and Nonresidues revisited:  Theorem 9
\endhead

\endtoc

\endtopmatter

\head 1. Introduction: Definitions and Properties of the Spectrum \endhead 

\noindent Let $S$ be a subset of the unit disc ${\Bbb {\Bbb U}}$, and let 
${\Bbb T}$ be the unit circle.  Let ${\Cal F}(S)$ denote the class of 
completely (totally)
multiplicative functions\footnote{That is, $f(mn)=f(m)f(n)$ for all positive integers $m,n$} $f$ such that $f(p)\in S$ for all primes $p$.  
Our main concern is: 

\centerline {\sl What numbers arise as mean-values of functions in ${\Cal F}(S)$?}

Precisely, we define 
$$
\Gamma_N(S) = \biggl\{ \frac 1N \sum\Sb n \le N\endSb f(n) : f\in {\Cal F}(S)
\biggr\}
\ \ \text{and} \ \ \Gamma(S) =\lim_{N \to \infty} \Gamma_N(S).
$$
Here and henceforth, if we have a sequence of subsets $J_N$ of the unit 
disc ${\Bbb {\Bbb U}}:= \{ |z| \le 1\}$, then by writing $\lim_{N \to \infty} J_N =J$
we mean that $z\in J$ if and only if 
there is a sequence of points $z_N\in J_N$ with $z_N\to z$ as $N \to \infty$.
We call $\Gamma(S)$ the {\sl spectrum} of the set $S$ and the object of this 
paper is to understand the spectrum.  
Although we can determine the spectrum explicitly
only in one interesting case (where $S=[-1,1]$), we are able, in general, to qualitatively describe it and obtain some of its
geometric structure.  For example, qualitatively, the 
spectrum may be described in terms of Euler products and solutions to 
certain integral equations. Geometrically, we can always 
determine the boundary points of the spectrum (that is, the elements of 
$\Gamma(S)\cap {\Bbb T}$) and show that the spectrum is connected. 
Moreover we can bound the spectrum, and make conjectures about some of
its properties, though we have no precise idea of what it usually
looks like.

We begin with a few immediate consequences of our definition:\

$\bullet$ \ \ $\Gamma(\{1\}) =\{1\}$.  

$\bullet$ \ \  If $S_1 \subset S_2$ then $\Gamma(S_1) \subset \Gamma(S_2)$.

$\bullet$ \ \  $\Gamma(S)$ is a closed subset of the unit disc ${\Bbb {\Bbb U}}$.

$\bullet$ \ \  $\Gamma(S) = \Gamma(\overline{S})$ 
(where $\overline{S}$ denotes the closure of $S$).  

\noindent {\sl Henceforth, we shall assume that $S$ is always closed.}

One of our main results, which formed the original motivation to
study the questions discussed herein,  is a precise description of the 
spectrum of $[-1,1]$.

\proclaim{Theorem 1} The spectrum of the interval $[-1,1]$ is the 
interval $\Gamma([-1,1])=[\delta_1,1]$ where 
$$
\delta_1 = 1-2\log(1+\sqrt{e}) + 4\int_1^{\sqrt{e}} \frac{\log t}{t+1} dt 
= - 0.656999 \ldots.
$$
\endproclaim

Theorem 1 tells us that for any real-valued completely multiplicative function
$f$ with $|f(n)| \le 1$,
$$
\sum_{n\le x} f(n) \ge (\delta_1 +o(1))x. \tag{1.1}
$$
In 1994, Roger Heath-Brown conjectured that there is some constant $c > -1$ 
such that $\sum_{n\le x} f(n) \ge (c+o(1)) x$.  Richard Hall [6] proved this
conjecture, and, in turn, conjectured (as did Hugh L.~Montgomery independently)
the stronger estimate (1.1).  Both Hall and Montgomery noticed that the
estimate (1.1) is best possible by taking 
$$
f(q) = \cases
1 &\text{for primes } q\le x^{1/(1+\sqrt{e})}\\
-1 &\text{for primes } x^{1/(1+\sqrt{e})} \le q \le x.\\
\endcases
\tag{1.2}
$$
In this example, the reader can verify (or see [6]) that equality holds 
in (1.1).
Our proof shows that this is essentially the only case when equality holds in
(1.1):

\proclaim{Corollary 1}  Let $x$ be sufficiently large, and let $f$ be 
any real-valued completely multiplicative function with $-1\le f(n)\le 1$.  
Then 
$$
\sum_{n\le x} f(n) \ge (\delta_1+o(1))x.
$$
Equality holds above if and only if 
$$
\sum_{p\le x^{1/(1+\sqrt{e})}} \frac{1-f(p)}{p} + \sum_{x^{1/(1+\sqrt{e})}
\le p\le x} \frac{1+f(p)}{p} = o(1).
$$
\endproclaim

By applying this Corollary to the completely multiplicative function $f(n) =
\fracwithdelims() {n}{p}$, for some prime $p$, we deduce that the 
number of integers below $x$ that are quadratic residues $\pmod p$ is
$$
\frac{1}{2} \sum_{n\le x} \Big( 1+ \fracwithdelims(){n}{p}\Big)
\ge \frac{1+\delta_1}{2} x +o(x) = (\delta_0  + o(1)) x,
$$
say.  In fact, the constant $\delta_0 =0.171500\ldots$\footnote{From 
the definition of $\delta_1$ in Theorem 1, we can derive the 
following curious expression for $\delta_0$: 
 $$
\delta_0= 1 - \frac{\pi^2}{6} - \log(1 + \sqrt{e}) \log\frac{e}{1 +
\sqrt{e}} + 2 \sum_{n= 1}^{\infty} \frac{1}{n^2} \frac{1}{(1 +
\sqrt{e})^n}.
$$}
More colloquially we have:
\smallskip
\centerline{\sl If $x$ is sufficiently large then, for all primes $p$, more than}

\centerline{\sl  $17.15\% $ of the integers up to $x$ are quadratic residues $\pmod p$.}
\smallskip 
\noindent  The constant $\delta_0$ here is best possible.  To see this, we choose $p$ 
such that $\fracwithdelims() qp$ is given as in the Hall-Montgomery example 
(1.2); 
that infinitely many such primes exist follows from quadratic reciprocity
and Dirichlet's theorem on primes in arithmetic progressions. 

Naturally, one wonders if similar results hold for $m$-th power residues.  
We partially answer this question by demonstrating that for any prime $\ell$,
the set of integers below $x$ that are $m$-th power residues $\pmod \ell$
has positive density, and that its logarithmic density exceeds $1/2^{m-1}$.

\proclaim{Theorem 2}  For integers $m\ge 2$, define
$$
\gamma_m = \liminf_{x\to \infty} \inf_{\ell} 
\frac{1}{x} \sum\Sb n\le x\\ n\equiv a^m\pmod \ell\endSb 1,
\qquad \text{and   }
\gamma_m^{\prime} = \liminf_{x\to \infty}\inf_{\ell} 
\frac{1}{\log x} \sum\Sb n\le x\\ n\equiv a^m \pmod \ell \endSb 
\frac 1n.
$$
Then $\gamma_2=\delta_0$, $\gamma_2^{\prime}=1/2$, and for $m\ge 3$, 
$$
0<\gamma_m\le \rho(m) \left( =\frac{1}{m^{m+o(m)}}\right) <  
\frac{1}{2^{m-1}} \le \gamma_m^{\prime} \le \min_{\beta\ge 0} 
\frac{1}{e^{\beta}} \sum_{k=0}^{\infty} \frac{\beta^{km}}{(km)!}
\left( \sim  \frac{1}{e^{m/e}}\right) .
$$
Here $\rho(u)$ is the Dickman-de Bruijn function, defined by $\rho(u) =1$
for $0 \le u \le 1$, and $u\rho^{\prime}(u) = -\rho(u-1)$ for all
$u \ge 1$.  
\endproclaim

We do not know the exact values of $\gamma_m$ and $\gamma_m^{\prime}$ for 
any $m\ge 3$.  By calculating 
numerically the minimum over $\beta$ in Theorem 2, we found that 
$\gamma_3^{\prime} \le 0.3245$, $\gamma_4^{\prime} \le 0.2187$,
$\gamma_5^{\prime} \le 0.14792$, and $\gamma_6^{\prime} \le 0.1003$. Theorem 2
implies

\smallskip
\centerline{\sl For given integer $m\geq 2$, there exists a constant $\pi_m>0$ such that}

\centerline{\sl if $x$ is sufficiently large then, for all primes $p$, more than }

\centerline{\sl  $\pi_m\% $ of the integers up to $x$ are $m$th power residues $\pmod p$.}
\smallskip

We now proceed to a more systematic treatment of the spectrum.
For a given $f \in{\Cal F}(S)$, the mean-value of $f$ (that 
is, $\lim_{x\to \infty} x^{-1} \sum_{n\le x} f(n)$),  if it exists, 
is obviously an element of the spectrum $\Gamma(S)$.  We begin by trying to 
understand the subset of the spectrum consisting of such mean-values.

Let $f$ be any multiplicative function with $|f(n)|\le 1$ for all $n$.  
Throughout this paper we define
$$
\Theta(f,x):= \prod_{p\le x} \biggl( 1+ \frac{f(p)}{p} +\frac{f(p^2)}{p^2}
+ \ldots \biggl) \biggl(1-\frac{1}{p} \biggr).
$$
In [13], A. Wintner showed, by a simple convolution argument, that 
if $\sum_p |1-f(p)|/p$ converges then
$$
\lim_{x\to \infty} \frac{1}{x}\sum_{n\le x} f(n) =
 \Theta(f,\infty); \tag{1.3}
$$
and so, if $f\in {\Cal F}(S)$ then $\Theta(f,\infty)\in {\Gamma}(S)$.
As an application of Wintner's result, we can take $f(p)=1$ for all $p>x$,
as long as $1\in S$,
so that $\Theta(f,x)\in {\Gamma}(S)$ for all $x$. 
Thus, if $1\in S$, we define the {\sl Euler product spectrum} of $S$ as
$$
\GT =  \lim_{x\to \infty} \{ \Theta(f,x): f\in{\Cal F}(S)\},
$$
which is a closed subset of $\Gamma(S)$.  If $1\notin S$ then 
define $\GT =\{0\}$. 
  
Proving an old conjecture of Erd{\H o}s and Wintner, 
Wirsing [14] showed that every real multiplicative function with 
$|f(n)|\le 1$ has a mean-value.  In fact, he proved that (1.3) always 
holds for such functions.  Thus, when $S\subset [-1,1]$ Wirsing's Theorem 
gives that $\GT$ is precisely the set of mean-values of elements in 
${\Cal F}(S)$.  In view of Wintner's result, 
the critical point in Wirsing's Theorem is to show that 
if $f$ is real valued and $\sum_p (1-f(p))/p$ diverges, then 
$x^{-1}\sum_{n\le x} f(n)\to 0$.  

The situation is more delicate for complex valued multiplicative 
functions.  For example, the function $f(n)=n^{i\alpha}$ ($\alpha$ a non-zero 
real) does not have a mean-value; indeed $\sum_{n\le x} f(n) \sim x^{1+i\alpha}
/(1+i\alpha)$.  Note that here $\sum_p (1-\text{Re } p^{i\alpha})/p$ 
diverges but $x^{-1}\sum_{n\le x} n^{i\alpha}$ does not tend to $0$.
Hal{\' a}sz [2] excluded this example by requiring that 
the set $\{f(p)\}$ be everywhere dense on ${\Bbb T}$.   
In fact, he proved that
if $\sum_{p} (1-\text{Re } f(p)p^{-i\xi})/p$ diverges (which  
obviously does not hold for the troublesome example $n^{i\alpha}$) 
for all real $\xi$ then $x^{-1}\sum_{n\le x} f(n) \to 0$; and he quantified 
how fast this tends to $0$.  

\proclaim{Lemma 1 (Hal{\' a}sz)}  Let $f$ be a multiplicative function with 
$|f(n)|\le 1$ for all $n$, and set 
$$
M(f,x) = \min_{|t|\le \log x} \sum_{p\le x} \frac{1-\RE f(p)p^{-it}}{p}.
$$
Then
$$
\sum_{n\le x} f(n) \ll x e^{-M(f,x)/16}.
$$
\endproclaim

Hal{\' a}sz comments that the factor $1/16$ in the exponent can be replaced
by the optimal constant $1$.  
Over the years Hal{\' a}sz' Theorem 
has been considerably refined ([5,7]), and recently Hall [5] found the 
following useful formulation.

\proclaim{Lemma 1$^{\prime}$ (Hall)} 
Let $D$ be a convex subset of ${\Bbb {\Bbb U}}$ containing 
$0$.  If $f\in {\Cal F}(D)$ then 
$$
\sum_{n\le x} f(n) \ll x \exp\biggl(-\eta(D) \sum_{p\le x} \frac{1-\text{Re }
f(p)}{p}\biggr), 
$$
where $\eta(D)$ is a constant determined by the geometry of $D$ (see [5]).  
In particular, if $\lambda(D)$ denotes the perimeter length of 
$D$ then $\eta(D) \ge (1-\lambda(D)/2\pi)/2$.
\endproclaim

As a byproduct of our investigations here, we 
have been able to obtain explicit quantitative versions of Lemma 
1 (with the strong exponent $1$ there) and Lemma 1$'$.  These 
will appear elsewhere.   

Lemmas 1 and 1$^\prime$ are 
important tools in all our subsequent work here.
We now note two immediate consequences of these results:
\  If $1\notin S$ then
(recalling that $S$ is closed)
$$
\sum_{p\le x} \frac{1-\RE f(p)}{p} \gg \sum_{p\le x} \frac{1}{p} 
= \log \log x +O(1).
$$
Further, we see easily that $S$ can be contained in a convex region with 
perimeter length $<2\pi$.  By Lemma 1$^{\prime}$, it follows that
$\Gamma(S) = \{0\}$.  Thus we state
\smallskip
$\bullet$ \ \  If $1\not\in S$ then $\Gamma(S) = \{0\}$.
  
\noindent {\sl Henceforth, we shall assume $1\in S$.}

Our second consequence characterizes the subsets $S$ of ${\Bbb {\Bbb U}}$ with the 
property that (1.3) holds for all $f\in{\Cal F}(S)$.  Wirsing's result 
states that subsets of $[-1,1]$ have this property.  To formulate 
our characterization fluidly, 
and for subsequent results, we introduce the notion of the {\sl angle} 
of a set.  For any $V\subseteq {\Bbb {\Bbb U}}$, define 
$$
\text{Ang}(V) := \sup\Sb v\in V\\ v \neq 1\endSb |\arg(1-v)|.  \tag{1.4}
$$
Note that each such $1-v$ has positive real part, so  
$0 \le \text{Ang}(V) \le \pi /2$.  We adopt the convention that Ang$(\{1\})
= \text{Ang} (\emptyset) = 0$.  Sometimes we will speak of the
angle of a point $z\in {\Bbb {\Bbb U}}$ ($z\neq 1$); by this we mean $\text{Ang}(z)
=|\arg (1-z)|$.  

\proclaim{Corollary 2}  Suppose $S \subset {\Bbb {\Bbb U}}$ and Ang$(S)<\pi/2$.
Then (1.3) holds for every $f\in {\Cal F}(S)$; that is, every $f\in{\Cal F}(S)$
has a mean-value.   Thus,
$$
\GT = \biggl\{ \lim_{N\to \infty} \frac{1}{N}\sum_{n\le N} f(n): \,\,\, 
f\in{\Cal F}(S)\biggr\} = 
\biggl\{ \Theta(f,\infty): \,\,\, f\in{\Cal F}(S)\biggr\}.
$$
\endproclaim

If $S\subset [-1,1]$ then Ang$(S)=0$, and so Corollary 2 generalizes Wirsing's
result.  If $\alpha \neq 0$ is real then Ang$(\{p^{i\alpha}\})=\pi/2$, and 
thus Corollary 2 avoids the example $f(n)=n^{i\alpha}$.  Corollary 2 follows 
from Wintner's result in the case that $\sum_p |1-f(p)|/p$ converges.  
If Ang$(S) <\pi/2$, and $f(p)\in S$ then $|1-f(p)| \asymp 1-\RE f(p)$.  
So the divergence of $\sum_p |1-f(p)|/p$ is equivalent to the 
divergence of $\sum_p (1-\RE f(p))/p$, and so by Lemma 1$^{\prime}$
the mean-value of $f$ is $0$.  This proves Corollary 2.

In general, the 
spectrum contains more elements than simply the Euler products.
For example, the spectrum of Euler products for $S =[-1,1]$ is simply 
the interval $[0,1]$.  However, as Theorem 1 shows, the spectrum of $S$ 
is more exotic.  We now describe a family of integral equations whose 
solutions belong to the spectrum.  In Theorem 3, we shall show that 
all points of the spectrum may be obtained by suitably combining an
Euler product and a solution to one of these integral equations.

Recall that we assume $S$ is closed and $1\in S$.  
We define $\Lambda(S)$ to be the set of 
values $\sigma(u)$ obtained as follows.  For a subset $S$ of the unit disc 
we denote by $S^*$ the convex hull of $S$.   Let $K(S)$ denote 
the class of measurable functions $\chi: [0,\infty) \to S^*$ 
with $\chi(t)= 1$ for $0\le t\le 1$.  
We prove in Theorem 3.3 (below) that associated to each $\chi$ 
there is a unique $\sigma: [0,\infty) \to {\Bbb U}$ satisfying 
the following integral equation 
$$
\align
u\sigma(u)= \sigma * \chi (u)=  \int_0^u \sigma(u-t)\chi(t) dt
 \ \ &\text{\rm for} \ \ u > 1 ,\tag{1.5} \\
\text{\rm with the initial condition} \  \ \sigma(u)=1 \ \ \ 
&\text{\rm for} \ \ 0\le u \le 1. \\
\endalign
$$
Here, and throughout, $f*g$ denotes the convolution of the two functions 
$f$ and $g$: that is, $f*g(x)=\int_0^x f(t)g(x-t)dt$.

That the integral equation (1.5) is relevant to the study of 
multiplicative functions was already observed by Wirsing [14].  This 
connection may be seen from the following Proposition.

\proclaim{Proposition 1} Let $f$ be a multiplicative function with $|f(n)|
\le 1$ for all $n$ and $f(n)=1$ for 
$n\le y$.  Let $\vartheta(x) =\sum_{p\le x} \log p$ and define 
$$
\chi(u) = \chi_f(u) = \frac{1}{\vartheta(y^u) } \sum_{p\le y^u} f(p)\log p. 
$$
Then $\chi(t)$ is a measurable function taking values in the unit disc
and with $\chi(t)=1$ for $t\le 1$. Let 
$\sigma(u)$ be the corresponding unique solution to (1.5).
Then 
$$
\frac{1}{y^u} \sum_{n\le y^u} f(n) =\sigma(u) +O\biggl(\frac{u}{\log y}\biggr).
$$
\endproclaim

The converse to Proposition 1 is also true:

\proclaim{Proposition 1 (Converse)} Let $S\subset {\Bbb {\Bbb U}}$ and 
$\chi \in K(S)$ be given.  Given $\epsilon >0$ and $u\ge 1$ there 
exist arbitrarily large $y$ and $f\in {\Cal F}(S)$ with $f(n)=1$ for 
$n\le y$ and 
$$
\biggl| \chi(t) - \frac{1}{\vartheta(y^t)} \sum_{p\le y^t} f(p)\log p
\biggr| \le \epsilon \ \ \text{for almost all } 0 \le t\le u.
$$
Consequently, if $\sigma(u)$ is the solution to (1.5) for this $\chi$ then 
$$
\sigma(t) = \frac{1}{y^t} \sum_{n\le y^t} f(n) + O(u^{\epsilon}- 1) 
+O\biggl(\frac{u}{\log y}\biggr) \ \ \text{for all } t\le u.
$$
\endproclaim

If $J$ and $K$ are two subsets of the unit disc, we define $J\times K$ to
be the set of elements $z=jk$ where $j\in J$ and $k\in K$. 

\proclaim{Theorem 3 (The Structure Theorem)}  
For any closed subset $S$ of ${\Bbb {\Bbb U}}$ with $1\in S$, 
$\Gamma(S) =\Gamma_{\Theta}(S) \times \Lambda(S)$. 
\endproclaim

Researchers in the field have previously used results like Proposition
1 and Theorem 3 in special, usually extreme, cases
(see [8, 10, 14], for instance), but this appears to be the first attempt to provide such a result in this generality.  The idea of the proof of Theorem 3 
is to decompose $f \in {\Cal F}(S)$ into two parts: $f_{s}(p)=f(p)$ 
for $p\le y$ and $f_s(p)=1$ for $p>y$, and $f_l(p)=1$ for $p\le y$ 
and $f_l(p)=f(p)$ for $p>y$.  For appropriately chosen $y$, the average 
of $f$ until $x$ is approximated by the product of the averages of $f_s$ and
$f_l$.  If $y$ is small enough compared with $x$,
then the average of $f_s$ is approximated by $\Theta(f_s,\infty)\in \GT$.  
Proposition 1 shows that if $y$ is not too small, the average of $f_l$ 
is approximated by the solution to an integral equation.  Combining these,
one gets that $\Gamma(S) \subset \GT \times \Lambda(S)$.  The proof that
$\GT \times \Lambda(S) \subset \Gamma(S)$ is similar, invoking the converse 
of Proposition 1.

As the case $S=[-1,1]$ illustrates, $\GT$ represents the easy part of 
the spectrum while $\Lambda(S)$ is more mysterious.  Here Theorems 1 and 3 
tell us that $\Lambda(S) \subset [\delta_1,1]$. That is, given any 
$\chi \in K([-1,1])$ we have $\sigma(u) \ge \delta_1$ for all $u$ 
(where $\sigma$ is the 
corresponding solution to (1.5)).  An important example is the function 
$\chi(t)=1$ for $t\le 1$ and $\chi(t)=-1$ for $t>1$.  Denote by $\rho_{-}(u)$
the corresponding solution to (1.5).  Then $\rho_-(u)$ satisfies a 
differential-difference equation very similar to that satisfied by the
Dickman-de Bruijn function.  Namely, $\rho_-(u)=1$ for $u\le 1$ and 
for $u>1$,
$$
u\rho_-^{\prime}(u)=-2\rho_-(u-1).
$$
It is not hard to verify that $\rho_-(u)$ decreases for $u$ in 
$[1,1+\sqrt{e}]$
and increases for $u>1+\sqrt{e}$.  The absolute minimum $\rho_-(1+\sqrt{e})$
is guaranteed by Theorem 1 to be $\ge \delta_1$ and in fact $\rho_-(1+\sqrt{e})
=\delta_1$.   By continuity, $\rho_-(u)$ takes on all values in the 
interval $[\delta_1,1]$ showing that $\Lambda(S) \supset [\delta_1,1]$.

We now describe properties of $\GT$, which are also inherited by 
$\Gamma(S)$:\ In many cases, we get an explicit description of $\GT$. 
To state our results we introduce the set ${\Cal E}(S)$ defined as follows:
\ If $1\in S\subset {\Bbb {\Bbb U}}$  define
$$
{\Cal E}(S) =\{ e^{-k(1-\alpha)}:\ \ k \ge 0, \ \ \alpha \ \ \text{is in the 
convex hull of } S\},
$$
so that  ${\Cal E}(S)$ consists of various ``spirals'' connecting $1$ to $0$. 

\proclaim{Theorem 4} {\rm (i)}\ For all closed subsets $S$ of ${\Bbb U}$ 
with $1\in S$,
$$
{\Cal E}(S) \times [0,1] \supset \GT = \GT \times {\Cal E}(S) \supset 
{\Cal E}(S).
$$
If $z\in \GT$ then $|z| \le \exp(-|\arg(z)| \cot( \text{Ang}(S)))$.  

\noindent{\rm (ii)}\  If the 
convex hull of $S$ contains a real point other than 
$1$, then 
$$
\GT = {\Cal E}(S) = {\Cal E}(S) \times [0,1].
$$   
In particular,  $\GT$ is {\rm starlike}; that is, $\GT$ 
contains each line joining $0$ to a point $z\in \GT$.
\endproclaim 

We may describe the set ${\Cal E}(S)$ explicitly as follows:  
If $S$ does not contain any element with positive imaginary part 
then define ${\Cal I}^+ = \emptyset$.  If $S$ does contain 
elements with positive imaginary part, then let $z^+$ be an element of 
$S$ with Im$(z^+) >0$ and such that Ang$(z^+)$ 
is the largest among all $z \in S$ with positive imaginary part.  
Define now ${\Cal I}^+$ to be the interior of the closed curve 
$\{e^{-k(1-z^+)}: \ \ 0 \le k\le 2\pi/|\IM z^+|\} \cup 
[e^{-2\pi (1-\RE z^+)/|\IM z^+|},1]$.  Similarly, define 
${\Cal I}^-$ by focussing on elements of $S$ with negative 
imaginary part.  Then ${\Cal E}(S)$ is contained in ${\Cal I}^+ \cup 
{\Cal I}^-$; and if the convex hull of $S$ contains 
a real point other than $1$, then ${\Cal E}(S)= {\Cal I}^+\cup {\Cal I}^-$.

It is easy to see that 
\smallskip
$\bullet$ \ \  If Ang$(S)=\pi/2$, then ${\Bbb {\Bbb U}}
=\overline{{\Cal E}(S)}$, 
so that $\Gamma(S)={\Bbb U}$.

\noindent {\sl Thus the spectrum is of interest only when Ang$(S)<\pi/2$.}

Combining Theorems 3 and 4 
enables us to deduce some basic properties of the spectrum (see \S 6b for the 
proof of this Corollary).

\proclaim{Corollary 3}  Let $S$ be 
a closed subset of ${\Bbb U}$ with $1\in S$. 
\item{\rm (i)} Then $\Gamma(S) =\Gamma(S) \times {\Cal E}(S)$.  
Consequently, the spectrum of $S$ is connected.  If 
the convex hull of $S$ contains a real point other than $1$, then
the spectrum is starlike, and contains the shape 
$\{ z: \ |z| \le \exp(-|\arg(z)| \cot( \text{Ang}(S))) \}$.
\item{\rm (ii)}  If $\alpha \in S$ then $1-
(1-\alpha)\log u \in \Lambda(S)$ for 
all $1\le u\le 2$.   If $\pi/2 > \text{Ang}(S) > 0$ then $\Gamma(S)$ 
contains elements not in $\GT$.   
\item{\rm (iii)}  If $1$, $e^{i\alpha}$ and $e^{i\beta}$ 
are distinct elements 
of $S$ then ${\Cal E}(S)$, and so $\Gamma(S)$, contains the disc centered at 
the origin with radius $\exp(-2\pi/(|\cot (\alpha/2)-\cot (\beta/2)|))$. 
\endproclaim

We have seen that the sets 
$\GT$ and $\Gamma(S)$ have the property that multiplying by ${\Cal E}
(S)$ leaves them unchanged.  It turns out that $\Lambda (S)$ also 
has this property, leading to  the following variant of Theorem 3, 
which reveals that $\Lambda(S)$ typically
contains all the information about the spectrum (see \S 6c for the 
proof of this Theorem).

\proclaim{Theorem 3$^{\prime}$} If 
$S$ is a closed subset of ${\Bbb U}$ with $1\in S$
then 
$$
\Lam(S) = \Lam(S) \times {\Cal E}(S),
$$
and 
$$
\Lambda(S) \subset \Gamma(S) \subset \Lambda(S) \times [0,1].
$$
If the convex hull of $S$ contains a real point different from $1$ then
$\Gamma(S) =\Lambda(S)$.
\endproclaim

Next we bound the spectrum and determine $\Gamma(S) \cap {\Bbb T}$.

\proclaim{Theorem 5} Suppose $S$ is a closed subset of ${\Bbb U}$ with 
$1\in S$.
The spectrum of $S$ is ${\Bbb U}$ if and only if {\rm Ang}$(S) =\pi/2$.  
If {\rm Ang}$(S) =
\theta < \pi/2$, then there exists a positive constant $A(\theta)$,
depending only on $\theta$, such that $\Gamma(S)$ is contained in a disc 
centered at $A(\theta)$ with radius $1-A(\theta)$.   In fact, 
$A(\theta) = (28/411) \cos^2 \theta$ is permissible.  Thus 
$$
\Gamma(S) \cap {\Bbb T} = \cases
\{1\} &\text {if Ang}(S)<\pi/2\\
{\Bbb T} &\text {if Ang}(S)=\pi/2.\\
\endcases
$$
\endproclaim

Applied to the set $S=[-1,1]$, Theorem 5 shows that there exists $c>-1$ such
that $\Gamma(S) \subset [c,1]$.  Thus Theorem 5 generalises Hall's 
result on Heath-Brown's conjecture.  

If $z\in S$ is such that Ang$(z)=$Ang$(S)=\theta$ then,
taking $k=\pi/|\IM z|$, we have 
$-\exp( -\pi \cot \theta) = e^{-k(1-z)}\in {\Cal E}(S) \subset 
\Gamma(S)$. Therefore 
$A(\theta) \leq (1-\exp( -\pi \cot \theta))/2 \leq \frac{\pi}{2}\cos\theta$.

By a simple calculation, we can show that ${\Cal E}(S)$, $\GT$ and $S$ all
have the same angle.  From Theorems 3 and 3$^{\prime}$ 
we see that Ang$(\Gamma(S))=$ Ang$(\Lambda(S)) \ge $ Ang$(S)$. We believe that
these angles are all equal:

\proclaim{Conjecture 1} The angle of the set equals the angle of the 
spectrum.  Thus
$$
\text{Ang}(\Gamma(S)) = \text{Ang} (\Lambda(S))= 
\text{Ang}(\GT) = \text{Ang}({\Cal E}(S)) = \text{Ang}(S).
$$
\endproclaim

Given $0\le \theta \le \pi/2$ define $H^{\theta}$ to be 
the subset of ${\Bbb U}$ inside the lines $\text{arg}(1-z) =\pm \theta$: 
thus, $H^{\theta}$ is the set of all points $z$ with Ang$(z)\le \theta$.  
If Conjecture 1 holds then taking $S=H^{\theta}$ there 
we deduce that $\Gamma(H^{\theta}) \subset H^{\theta}$.  Conversely, 
if $\Gamma(H^{\theta}) \subset H^{\theta}$ then for any $S\subset {\Bbb U}$ 
with Ang$(S)=\theta$ we must have $S\subset H^{\theta}$ 
and so $\Gamma(S) \subset \Gamma(H^{\theta}) \subset H^{\theta}$.  
It follows at once that Ang$(\Gamma(S)) = $ Ang$(S)$.  Thus 
Conjecture 1 is equivalent to the following:

\proclaim{Conjecture 1$'$}  With $H^{\theta}$ as defined above 
$\Gamma (H^{\theta}) \subset H^{\theta}$.
\endproclaim

We support Conjecture 1 by showing that Ang$(S)$ and 
 Ang$(\Gamma(S))$ 
are comparable in the situations Ang$(S)\to 0$ and Ang$(S) \to \pi/2$.  

\proclaim{Theorem 6}  Suppose $S \subset {\Bbb U}$ and {\rm Ang}$(S)
=\theta =\pi/2-\delta$.
\item{\rm (i)} Then, {\rm Ang}$(\Gamma(S)) \ll $ {\rm Ang}$(S)$.
\item{\rm (ii)} Further,
$$
\frac {\pi}2 -\delta =\text{Ang}(S) \le \text{Ang}
(\Gamma(S)) \le \frac{\pi}{2}
-\frac{\sin \delta }{2}.
$$
\endproclaim

The first part of the Theorem says that Ang$(S)$ and Ang$(\Gamma(S))$ are
comparable when Ang$(S)$ is small.  The second part of the Theorem is 
mainly interesting 
in the complementary case when Ang$(S)$ is close to $\pi/2$.
In fact, when $\delta$ is small we see that we are away from the truth 
only by a factor of $2$ (as $\sin \delta \sim \delta$).  

\remark {Example} 
Let $k\ge 3$ and $S_k$ denote the set of $k$-th roots of unity.
If $f\in {\Cal F}(S_k)$ then $f(n)\in S_k$ for all $n$.  Hence $\Gamma(S_k)$ 
is contained in the convex hull of $S_k$: that is, in the regular $k$-gon with
vertices the $k$-th roots of unity.  Notice that this implies Ang$(\Gamma(S_k)
) \le $ Ang$(S_k)$, so that Ang$(S_k) = $ Ang$(\Gamma(S_k))$ by 
Theorem 6(ii), 
supporting Conjecture 1.
Applying Corollary 3(iii) with the two points $e^{\pm 2\pi i/k}$,
we conclude that $\Gamma(S_k)$ is starlike and contains the disc centered at 
$0$ with radius $\exp(-\pi \tan(\pi/k))$.  Even in this simple case
we have not been able to determine the spectrum $\Gamma(S_k)$, though
we do know that $2\pi^3/3k^2 + o(1/k^2)\le \pi - $Area$(\Gamma(S_k))\le
2\pi^3/k+o(1/k)$.
\endremark

We define {\sl the projection of} (a complex number) $z$ {\sl in the 
direction $e^{i\alpha}$} to be $\RE (e^{-i\alpha} z)$.  Theorem 1 may be 
re-interpreted as stating that if $z \in \Gamma(\{\pm 1\})$ then the 
projection of $z$ in the direction $-1$ is $\le -\delta_1$. 
 Evidently if $1\in S$
then $1\in \Gamma(S)$ so there is always a $z\in \Gamma(S)$ whose projection
in the direction 1, is 1, and thus uninteresting to us.  This motivates
us to define the {\sl maximal projection of the spectrum $\Gamma(S)$
of a set} $S \subset {\Bbb T}$ as
$$
\max_{1\neq \zeta \in S} \max_{z \in \Gamma(S)} \RE (\zeta^{-1} z).
$$

\proclaim{Conjecture 2}  Let $S$ be a closed subset of ${\Bbb T}$ 
with $1\in S$.  If $\text{Ang}(S) = \theta$ 
then the maximal projection of $\Gamma(S)$ is 
$$
\max_{1\neq \zeta \in S} \max_{z\in \Gamma (S) } \RE (\zeta^{-1} z) = 1- 
(1+\delta_1) \cos^2 \theta.
$$
\endproclaim

One half of this conjecture is easy to establish: namely, the 
maximal projection is $\ge 1 - (1+\delta_1) \cos^2 \theta$.  
To see this, let
$z=x^{-1} \sum_{n\leq x} f(n)$ where $f$ is the 
completely multiplicative function defined by $f(p)=1$ for all 
$p\leq x^{1/(1+\sqrt{e})}$, and $f(p) =\zeta$ for 
$x^{1/(1+\sqrt{e})} \le p \le x$, where $\zeta \in {\Bbb T}$ and 
Ang$(\zeta)=\theta$.  Then, a simple calculation (analogous to the 
calculation in the Hall-Montgomery example (1.2)) gives that 
the projection of $z$ along $\zeta$ is $1-(1+\delta_1)\cos^2\theta +o(1)$.

\proclaim{Theorem 7} 
\item{(i)} Conjecture 2 is true for the sets 
$S=\{ 1,-1\}$ and $S=\{ 1,-1, i, -i\}$.
\item{(ii)} For any closed subset $S$ of ${\Bbb T}$ with $1\in S$, 
the maximal projection of $\Gamma(S)$ is  
$\leq 1- (56/411) \cos^2 \theta$, where $\theta=$ {\rm Ang}$(S)$.
\endproclaim

To facilitate comparison between Theorem 7 and Conjecture 2, we observe that
$1+\delta_1 =0.3430\ldots$ whereas $56/411 = 0.1362\ldots$. 
Thus Theorem 7 is not too far away from the (conjectured) truth.

Let $S\subset {\Bbb T}$ and $\theta=\text{Ang}(S)$ and define $\alpha:=e^{i(\pi-2\theta)}$.
Let $S_\theta:=\{ 1, \alpha, \overline{\alpha} \}$ and 
$S^{\theta}$ be $\{ 1\}$ together with the arc of 
${\Bbb T}$ anticlockwise from $\alpha$ to $\overline{\alpha}$.
Then $S\subset S^\theta$, and if $S$ is symmetric about the real axis then
$S_\theta\subset S$. Conjecture 2 is equivalent to the conjecture that $\Gamma(S)\subset \Gamma(S^\theta)$ is contained inside the
arc of the circle of radius $r_\theta:=1 - (1+\delta_1) \cos^2 \theta$,
centered at the origin, going
anticlockwise from $r_\theta \alpha$ to $r_\theta \overline{\alpha}$, and inside the tangent lines to the circle from these 
two points going to the right.  We suspect that one should be able to  restrict $\Gamma(S^\theta)$ more than as in Conjectures 1 and 2, particularly on the left side ($\RE(z)<0$) of the plane.

If $S_\theta\subset S$ then  $\Gamma(S_\theta)\subset \Gamma(S)$. Collecting
several results above, we have seen that $\Gamma(S_\theta)$
contains the interior of the shape given by the line joining $1$ to 
$1-(1-\alpha)\log 2$, the contour $c(u)=1-(1-\alpha)\log u + ((1-\alpha)^2/2)
\int_{t=1}^{u-1} (\log (u-t)/t) dt$ for $2\leq u\leq 1+\sqrt{e}$, and the 
spiral $c(1+\sqrt{e})e^{-t(1-\alpha)},\ t\geq 0$ until it hits  the real
axis, along with their complex conjugates.

In \S 8 we investigate other notions of spectrum. For fixed $\sigma>0$,
the spectrum of
$$
 \lim_{x\to \infty} \Big\{ 
 \sum_{n\le x} \frac{f(n)}{n^\sigma} \bigg/  \sum_{n\le x} \frac{1}{n^\sigma} : \qquad f\in {\Cal F}(S) \Big\}
$$
is evidently determined by the Euler products if $\sigma>1$, and turns
out to be the same as $\Gamma(S)$ for $0<\sigma<1$, as we show at the beginning
of section 8. Thus the only new and interesting case is where $\sigma=1$,
which gives the  {\sl logarithmic spectrum}, $\Gamma_0(S)$.
As might be expected, the logarithmic spectrum is easier to study than $\Gamma(S)$.  In fact $\Gamma_0(S)$ lies inside the convex hull of
$\Gamma(S)$.
Our next result allows us to bound $\Gamma_0(S)$ independently of $\Gamma(S)$.

\proclaim{Theorem 8} Suppose $S$ is a closed subset of ${\Bbb U}$ with $1\in S$,
and let ${\Cal R}$ denote the 
closure of the 
convex hull of the points $\prod_{i=1}^{n} \frac{1+s_i}{2}$, for all $n\ge 1$,
and all choices of points $s_1$, $\ldots$, $s_n$ lying in the convex hull of
$S$.  Then $\Gamma_0(S)$ is contained in ${\Cal R}$.
\endproclaim

As a consequence of Theorem 8, we have $\Gamma_0([-1,1])=[0,1]$, 
and also, lending credence to Conjecture 1, that Ang$(S)= \text{Ang}
(\Gamma_0(S))$.

\proclaim{Corollary 4} Let $S$ be a closed subset of ${\Bbb U}$ with $1\in S$.  
\item{\rm (i)} $\Gamma_0([-1,1]) = [0,1]$.
\item{\rm (ii)} If $\alpha \in S$ then $1-(1-\alpha)(\log u-1+\frac 1u) 
\in \Gamma_0(S)$ 
for all $1\le u\le 2$.  If $0<\text{Ang}(S) <\pi/2$ then $\Gamma_0(S)$ 
contains elements not in $\GT$.
\item{\rm (iii)} $\text{Ang}(S) =\text{Ang}(\Gamma_0(S))$.
\item{\rm (iv)} Suppose $\text{Ang}(S)=\frac {\pi}{2} -\delta$ with $\delta >0$. 
If $z\in \Gamma_0(S)$ then $|z| \le (\cos \delta)^{|\arg (z)|/\delta}$ 
where we choose $|\arg(z)| \in [0,\pi]$.
\endproclaim 

Most of the ideas above generalize to the spectrum of {\sl all} multiplicative
functions in ${\Bbb U}$; that is where $f(mn)=f(m)f(n)$ for all pairs of
coprime integers $m,n$.  Thus the mean-value of $f$ depends now on the 
(independent values of) $f(p^k)$ with $k\geq 2$ as well as the $f(p)$.
A priori it is not obvious what range we should allow for the $f(p^k)$; it seems
that the most useful choices are $f(p^k)\in S=S_m$, when $S$ is the set of
$m$th roots of unity and, otherwise,  $f(p^k)\in {\Bbb U}$ for
all $k\geq 2$. We call this new spectrum $\hat{\Gamma}(S)$, and note that 
$\Gamma(S)\subset \hat{\Gamma}(S)$. Moreover we define  
$\hat{\Gamma}_{\Theta}(S)$ to be the set of values $\Theta(f,x)$ as before.
Now Theorem 1, Corollary 1, Lemmas 1, Corollary 2, and Propositions 1 all
hold, Theorem 3 with $\hat{\Gamma}(S) =\hat{\Gamma}_{\Theta}(S) \times \Lambda(S)$.  The most significant change is that the analogue to Theorem 4
is not true since $\hat{\Gamma}_{\Theta}(S)$ is not necessarily a subset of
${\Cal E}(S)\times [0,1]$. For example, if $S=\{ 1,-1,i,-i\}$ take
$f$ for which $f(p^k)=i$ for each $k\geq 1$, and $f(q^k)=0$ if $q\ne p$, so
that $z=\Theta(f,\infty)=1-1/p+i/p$ does not satisfy  $|z|\leq e^{-|\arg(z)|}$.
Changes thus need to be made in subsequent results, which are easy but messy,
and the theory necessarily loses
some of its elegance since, now, $\hat{\Gamma}(S)$ rarely equals $\Lambda(S)$.
Note also that Conjecture 1 is untrue for $\hat{\Gamma}(S)$ since if
$f(p)=\beta^2$ and $f(p^k)=\beta$ for each $k\geq 2$ where $\beta=e^{i(\pi/2-\text{Ang}(S))}$, and $f(q^k)=0$ if $q\ne p$, then
Ang$( \Theta(f,\infty) )>$Ang$(S)$.

Define for $B>0$ 
$$
\alpha(B) = \limsup_{|D| \to \infty} \frac{1}{(\log |D|)^B} 
\sum_{n\le (\log |D|)^B}\fracwithdelims() {D}{n} 
$$
and 
$$
\beta(B) = \liminf_{|D| \to \infty} \frac{1}{(\log |D|)^B} 
\sum_{n\le (\log |D|)^B}
\fracwithdelims() {D}{n},$$
where, $D$ represents a fundamental discriminant. 
Plainly $\alpha(B)=1$ for $B\le 1$, and in [1] we showed that 
$\alpha(B)\ge \rho(B)$ where $\rho$ is the Dickman-de Bruijn function.  
Further, we showed there that if the Generalized Riemann Hypothesis holds 
then $\alpha(B) \le \rho(B/2)$.  The exact value of $\alpha(B)$ is 
not known for any $B>1$, though we do conjecture that 
$\alpha (B)= \rho(B)$ for all $B>0$.  Regarding $\beta$, we see from Theorem 1 
that $\beta(B)\ge \delta_1$ for all $B$ and, in view of the Hall-Montgomery 
example,  $\beta(B) =\delta_1$ for $B\le 1$.  
Hybridizing this consequence of Theorem 1 and our result on $\alpha(B)$, 
Mark Watkins asked us whether $\beta(B)<0$ for all $B$.  
We see below that this is indeed so.

\proclaim{Theorem 9}  Given $u\ge 1$, let ${\Cal C}(u)$ 
denote the set of all measurable functions $\chi$ such that $\chi(t)=1$ 
for $t\le 1$, $\chi(t) \in [-1,1]$ for $1\le t\le u$, and 
$\chi(t) =0$ for $t>u$.  Define 
$$
\gamma(B) = \min_{u \ge 1} \min\Sb \chi \in {\Cal C}(u)\endSb \sigma(Bu),
$$
for all $B>0$, where $\sigma$ refers to the solution to (1.5).  Then $\beta(B)\le \gamma(B)$ for all $B>0$,  where 
$-\rho(B) \le \gamma(B) <0$.  
\endproclaim

Assuming the GRH, we can show that $\beta(B)\ge \gamma(B/2)$.  In [1], we 
gave our reasons for believing that $\alpha(B)=\rho(B)$; these 
also lead us to believe that $\beta(B)=\gamma(B)$ for all $B$.

To help orient the reader we supply a brief overview of the 
following sections, and describe the logical dependencies 
among them.  The reader interested in a proof of the structure theorem 
can skip \S 2 and proceed to \S 3a and \S 4.  After this a perusal 
of \S 3b and \S 5 would lead to a proof of Theorem 1.  The bulk of our general
results on the spectrum are covered in \S 6 and \S 7; both these sections 
build upon the work of \S 3 and \S 4.  Next \S 8 deals with other notions of 
spectrum, chiefly the logarithmic spectrum.  Again the material of \S 3 and 
\S 4 is assumed here.  Finally \S 2 and \S 9 may be read independently of the 
rest of the paper.

\head 2. The natural and logarithmic densities 
of $m$th power residues up to $x$
\endhead

\subhead 2a. The proportion of $m$th power residues up to $x$ \endsubhead

\noindent As noted in the introduction, it is clear that $\gamma_2=\delta_0$.
Given a set of $m$-th roots of unity $\alpha_p$ for each prime $p\le x$,
we see (by the Chebotarev density theorem) that there are infinitely 
many primes $\ell \equiv 1\pmod m$ such that there is a character $\chi 
\pmod \ell$ of order $m$ for which $\chi(p)=\alpha_p$ for all $p\le x$.  
Choose $\alpha_p=1$ for $p\le x^{1/m}$, and $\alpha_p = e^{2\pi i/m}$ 
for $x^{1/m}< p \le x$.  Then an integer $n\le x$ is an $m$-th power residue 
$\pmod \ell$ if and only if all its prime divisors are $\le x^{1/m}$.  
It is well-known that the number of such integers is $(\rho(m)+o(1))x
=m^{-m+o(m)}x$.  This gives the upper bound $\gamma_m \le \rho(m) 
= m^{-m+o(m)}$. 

We now show that $\gamma_m >0$ for $m\ge 3$.  To this end, we require the
following result of Hildebrand [9].

\proclaim{Lemma 2.1 (Hildebrand)}  Fix $\theta>0$. In the two limits below 
the $\sup$ and $\inf$ are taken over all completely multiplicative
functions $f$ with $0\leq f(n) \leq 1$, such that 
$\Theta(f,x) = e^{-\theta} + o(1)$. We have
$$
\lim_{x \to \infty} \inf \frac 1x \sum_{n\le x} f(n) =  \rho(e^{\theta})
\ \ \text{and} \ \ 
\lim_{x \to \infty} \sup \frac 1x \sum_{n\le x} f(n) \leq  e^{-\theta}
\int_0^{e^\theta} \rho(t) dt.
$$
The lower bound is attained when $f(p)=1$ for all 
$p\leq x^{e^{-\theta}}$, and  $f(p)=0$ for all larger primes $p$.
\endproclaim

The exact value of the lim sup above is still not known, though it must be
at least the average value, $\ge e^{-\theta}$.  Also 
$\int_0^\infty \rho(t) dt =e^{\gamma}$, and so the upper bound 
given above is not too far from the truth.  In our application, 
it is the lim inf result that is useful.

\proclaim{Proposition 2.2}  Suppose  $m$  is a given positive integer
and  $c \ge 0$  is a given constant.  For any sufficiently large
integer  $n$, and prime  $\ell > n$, with  $\ell \equiv 1 \pmod{m}$
suppose that for some divisor  $M$  of  $m$  one has
$$
\sum_{p\in P} \frac{1}{p} \leq c
$$
where  $P$  is the set of primes  $\le x$  that are not  $M$th
power residues  $\mod{\ell}$.  Then

\noindent {\rm Either} more than  $\frac{M}{2m}\sum_{n\le x,\ (n,P) = 1} 1$  
integers up to  $x$, that are coprime to $P$, are
$m$th power of residues  $\mod{\ell}$;

\noindent {\rm Or} there exists a divisor  $d > 1$  of  $m/M$  such that
$$
\sum_{q\in Q} \frac{1}{q} \le \kappa (c,m)
$$
where  $Q$  is the set of primes  $\le n$  that are not  $Md$th
power residues  $\mod{\ell}$.  Here  $\kappa (c,m)$  is a constant that
depends only on  $c$  and  $m$.
\endproclaim

\demo{Proof}  Let  $G$  be a set of coset representatives for the
characters  $\pmod{\ell}$  of order dividing  $m$  modulo the
characters  $\pmod{\ell}$  of order dividing  $M$.  Note that if
$n$  is an  $M$th power  $\pmod{\ell}$  then
$$
\sum_{x\in G} \chi(n) = \cases
|G| = m/M &\text{if  $n$  is an $m$th power  $\pmod{\ell}$}\cr
0 &\text{otherwise}.\cr
\endcases 
$$
So suppose for each  $\chi \in G$, except the identity  $\chi_0$  one
has
$$
\biggl| \sum\Sb n\le x \\ (n,P) = 1\endSb \chi(n)\biggr| \le
\frac{1}{2|G|} \sum\Sb n\le x \\ (m,P) = 1 \endSb 1 \tag{2.1}
$$
Then the number of  $m$th powers  $\mod{\ell}$  up to  $x$  is
$$
\align
&\ge \frac{1}{|G|} \sum_{\chi \in G} \sum\Sb n\le x \\ (m,P) =
1\endSb \chi(n) \geq 
 \frac{1}{|G|} \biggl( \sum\Sb n\le x \\ (n,P) = 1 \endSb 1 -
(|G| - 1) \cdot \frac{1}{2|G|} \sum\Sb n \le x \\ (n,P) = 1 \endSb
1\biggr) \\
&\ge \frac{1}{2|G|} \sum\Sb n\le x \\ (n,P) = 1 \endSb 1 =
\frac{M}{2m} \sum\Sb n\le x \\ (n,P) =1 \endSb 1
\endalign
$$
On the other hand, if (2.1) does not hold for some  $\chi \in G$,
$\chi \ne \chi_0$, then suppose  $\chi$  has order  $d > 1$  in
$G$.  (Thus  $\chi(p) \ne 1$  if and only if  $p \in Q$.)

Then, by Lemma 1$^{'}$ (with  $D$  the convex hull formed by the
$m$th roots of unity),
$$
\frac{M}{2m} \sum \Sb n\le x \\ (n,P) = 1 \endSb 1 \ll_m x
\exp\biggl( -c_m \sum_{p\in Q} \frac{1}{p}\biggr) .
$$
By the first part of Lemma 2.1,  $\operatornamewithlimits{\sum}\limits\Sb
n\le x \\ (n,P) = 1 \endSb 1 \gtrsim x\rho(c)$.  Thus
$\operatornamewithlimits{\sum}\limits_{p\in Q} \frac{1}{p} \ll_{c,m}
1$  which gives the result.
\enddemo

\demo{Proof of Theorem 2}  
First we may assume that $x$ is sufficiently large for the
argument below to work. Second we may assume that $m$ divides $\ell-1$
else we  replace $m$ by gcd$(\ell-1,m)$. Third we may assume that
$\ell>x$ else the proportion of such integers is certainly $\gg 1/m$ from 
elementary considerations.

First take the Proposition with  $c_1 = 0$, $P_1 =
\emptyset$, $M_1 = 1$.  Either the result follows immediately with
$\gamma_m \geq \frac{1}{2m}$  or there exists an integer  $d$, as
described.  Let  $P_2 = Q$, $M_2 = d$, and  $c_2 = \kappa (0,m)$.  If so,
apply the Proposition again; either we get  $\geq \frac{M_2}{2M}
\operatornamewithlimits{\sum}\limits\Sb n\le x \\ (n,P_2) = 1\endSb
1$  such integers as desired, and this is  $\gtrsim \frac{M_2}{2m}
\rho(c_2)x$  by Hildebrand's Lemma; or we get another integer
$d_2$  as described.  If so apply the Proposition again and again
with
$$
P_{k+1} = Q_k, \quad M_{k+1} = d_k \quad\text{and}\quad c_{k+1} =
k
\kappa (c_k,m) .
$$
The process eventually terminates (since  $M_1\mid M_2 \mid M_3 \mid 
\dots \mid m$  and each  $M_{k+1} > M_k$); when it does we get
$$
\gtrsim \frac{M_k}{2m} \rho(c_k)x \quad\text{integers up to}\quad x
$$
which are  $m$th power residues  $\mod{\ell}$.  Thus  $\gamma_m$
exceeds the minimum of the  $\frac{M_k}{2m}\rho(c_k)$  over all possible
such sequences  $M_1\mid M_2 \mid \dots \mid m$  (of which there
are evidently only finitely many).
\enddemo

\subhead 2b. Logarithmic proportions of $m$th power residues \endsubhead

\noindent It is plain that $\gamma_m^{\prime} \le 1/m$, and so in particular
$\gamma_2^{\prime}\le 1/2$.  Let $\beta >1 $ be a parameter to be chosen 
shortly, and put $\alpha_p = e^{2\pi i/m}$ if $(\log \log x)^{1/\beta}
\le p\le \log \log x$, and $\alpha_p=1$ for all other primes $p\le x$.  
Choose $\ell \equiv 1\pmod m$ such that there is a character $\chi\pmod \ell$
of order $m$ with $\chi(p)=\alpha_p$ for all $p\le x$.  Let $P$ denote the
product of those primes $p\le x$ with $\alpha_p\neq 1$.  We may 
write every $n\le x$ uniquely as $NR$ where $p|N\implies p|P$, and
$p|R \implies p\nmid P$.  Note that $\chi(n)=\chi(N) =1$ if and only if the
number of primes dividing $N$ (counted with multiplicity) is a 
multiple of $m$.  Thus 
$$
\align
\gamma_m^{\prime} &\le \frac{1}{\log x} \sum\Sb n\le x\\ \chi(n)=1
\endSb \frac 1n \le \frac{1}{\log x} \sum\Sb R\le x\\(R,P)=1\endSb \frac{1}{R}
\sum\Sb N\le x, \chi(N)=1\\ p|N\implies p|P \endSb \frac{1}{N}\\
&\le \biggl(\frac{\phi(P)}{P} +o(1)\biggr) \sum_{k=0}^{\infty} 
\frac{1}{(km)!} \biggl(\sum_{p|P} \frac{1}{p-1}\biggr)^{km}.\\
\endalign
$$
Letting $x\to \infty$, and using the prime number theorem we obtain
$$
\gamma_m^{\prime} \le \frac{1}{\beta} \sum_{k=0}^{\infty} 
\frac{(\log \beta)^{km}}{(km)!}.
$$
Taking $\beta=e^{m/e}$, and using Stirling's formula, it follows that 
$\gamma_m^{\prime} \ll e^{-m/e}$.  

\proclaim{Lemma 2.3}  Let  $a_1, \dots , a_n$ and 
$R_1, R_2,\dots , R_n, m\geq 2$ be integers. Then
$$
\#\{ (r_1,\dots ,r_n)\in {\Bbb Z}^n:   \sum_{i=1}^n r_ia_i \equiv 0
\pmod{m},\ \text{with} \ 0\leq r_i\leq R_i-1 \} \geq \frac{R_1\dots R_n}{2^{m-1}}.
$$
\endproclaim

Note that Lemma 2.1 is `best possible' in that if
$R_1=\dots =R_{m-1}=2$ and $n=m-1$ then the only solution has
each $r_i=0$, and thus we get equality above.

On the other hand, we naively expect the proportion typically to be close
to $1/m$ (rather than be as small as $1/2^{m-1}$); and if this is so in
context then we might expect to improve Corollary 2. below.

\demo{Proof} Given the $a_i, R_i$ and $m$ above we define $p(\text{\bf a,R})$,
the proportion of the sums that equal zero $\pmod m$, to be
equal to
$$
\frac 1{R_1\dots R_n}
\#\{ (r_1,\dots ,r_n)\in {\Bbb Z}^n:   \sum_{i=1}^n r_ia_i \equiv 0
\pmod{m},\ \text{with} \ 0\leq r_i\leq R_i-1 \} .
$$
The result that we wish to prove is that $p(\text{\bf a,R})\geq 1/2^{m-1}$.
Let us suppose that we have a counterexample above with 
$s(\text{\bf R}):=\sum_i (R_i-2)$ minimal. 

We will show that we must have each $R_i=2$, else if $R_n\geq 3$ then we 
will construct two new examples \text{\bf b,B} and \text{\bf c,C} with 
$s(\text{\bf R}) > s(\text{\bf B}), s(\text{\bf C})$, and with
$p(\text{\bf b,B})<1/2^{m-1}$ or 
$p(\text{\bf c,C})<1/2^{m-1}$, thus contradicting the minimality
of the purported counterexample $\text{\bf a,R}$. We thus have reduced
proving Lemma 2.3 to the case where every $R_i=2$, which we 
prove in Lemma 2.4 below.

Now we construct \text{\bf b,B} and \text{\bf c,C} as follows:\
Let $b_i=c_i=a_i$ and $B_i=C_i=R_i$ for $1\leq i\leq n-1$.
Let $b_n=b_{n+1}=a_n$ with $B_n=R_n-1$ and $B_{n+1}=2$;
and let $c_n=(R_n-1)a_n$ and $C_n=2$. We see that
$s(\text{\bf B})=s(\text{\bf R})-1$ and 
$s(\text{\bf C})=s(\text{\bf R})-(R_n-2)\leq s(\text{\bf R})-1$.

Now for $0\leq j\leq R_n-1$ we define
$$
\tau_j=\#\{ (r_1,\dots ,r_{n-1})\in {\Bbb Z}^n: \  
ja_n+ \sum_{i=1}^{n-1} r_ia_i \equiv 0 \pmod{m},\ 
\text{with} \ 0\leq r_i\leq R_i-1 \} .
$$
Thus with $R:=R_1\dots R_{n-1}$ we have 
$$
\align
p(\text{\bf a,R}) &=\frac 1{R_nR }
\{ \tau_0 + \tau_1 + \dots + \tau_{R_n-1} \}; \
p(\text{\bf c,C}) = \frac 1{2R } \{ \tau_0 + \tau_{R_n-1} \}; \\
\text{and} \ p(\text{\bf b,B}) &=\frac 1{2(R_n-1)R }
\{ \tau_0 + 2(\tau_1 + \dots + + \tau_{R_n-2}) +\tau_{R_n-1} \} .\\
\endalign
$$
Therefore 
$$
\frac 1{2^{m-1}} > p(\text{\bf a,R}) =\frac 1{R_n}
\{ (R_n-1)p(\text{\bf b,B})+p(\text{\bf c,C})\} 
\geq \min \{ p(\text{\bf b,B}), p(\text{\bf c,C}) \} ,
$$
as required.
\enddemo

\proclaim{Lemma 2.4}  Let  $a_1, \dots , a_n$ and $m\geq 2$ be integers. Then
$$
\#\{A \subseteq \{1, \dots , n\}:  \sum_{i\in A} a_i \equiv 0
\pmod{m}\} \geq 2^{n-(m-1)}.
$$
\endproclaim

\demo{Proof} If $n\leq m-1$ the statement is trivial since we always can take 
$A$ to be the empty set and thus get at least one such sum.  We
will assume henceforth that $n\geq m$.

Let $A_0$ be the subsequence of $a_i$'s which are 
$\equiv 0 \pmod m$,
and then let $A_1=\{1, \dots , n\} \setminus A_0$. We shall define a sequence
of subsets $\{ B_j\}_{j\geq 1}$ of $A_1$, with $B_1\subset B_2 \subset B_3 \subset \dots$ and each $B_j$ having exactly $j$ elements; and we will let
$C_j$ be the set of sums $\pmod m$, of the subsets of $B_j$.

We define $B_1=\{ b_1\}$ where $b_1$ is any element of $A_1$, so that 
$C_1=\{ 0,b_1\}$ has two elements. Given $B_j$ (and thus $C_j$) we attempt
to select $b_{j+1}\in A_1 \setminus B_j$, so that 
$C_{j+1}$ is larger than $C_j$. If this is possible we so construct $B_{j+1}$
(that is, as $B_j\cup \{ b_{j+1}\}$) and move on to attempting
the analogous construction with $j+1$; note that then
$C_{j+1}$ contains at least $j+2$ elements. If this construction is impossible, 
write $j=k$, and note that we must have $b+C_k \subset C_k$ for every 
$b\in A_1 \setminus B_k$.
But since $0\in C_k$ this would imply that $b, 2b, 3b, \dots \in C_k$.
Indeed by repeatedly using the relation $b+C_k \subset C_k$, we see that 
the additive subgroup $S$, generated by the elements of $A_1 \setminus B_k$,
must be a subset of $C_k$. 

In fact there must be such a value of $k$, since if not then 
we would have $m\geq |C_n| \geq n+1$ which gives a contradiction. Note that 
$m\geq |C_k| \geq k+1$.

Now select any subset $R$ of $A_0$, and any subset $T$ of 
$A_1 \setminus B_k$. Note that 
$s:=\sum_{a\in R} a + \sum_{a'\in T} a' \equiv \sum_{a'\in T} a' \pmod m$,
so that $s\in S$, and thus $-s\in S \subset C_k$. Therefore, by the definition
of $C_k$, there exists a subset ${\Bbb U}$ of $B_k$ with 
$\sum_{a''\in {\Bbb U}} a'' \equiv -s \pmod m$. Thus we have 
$\sum_{a\in R\cup T\cup {\Bbb U}} a \equiv 0 \pmod m$, and so
$$
\#\{A \subseteq \{1, \dots , n\}:  \sum_{i\in A} a_i \equiv 0
\pmod{m}\} \geq \sum_{R\subset A_0} \sum_{T\subset A_1 \setminus B_k} 1
= 2^{n-k} \geq  2^{n-(m-1)}.
$$
\enddemo

\proclaim{Corollary 2.5} Let $f$ be a completely multiplicative function where
each $f(p)$ is an $m$th root of unity. Then
$$
\frac{1}{\log x} \sum\Sb n\le x \\ f(n)=1\endSb 
\frac{1}{n} 
\geq \frac 1{2^{m-1}} + o(1) .
$$
\endproclaim

\demo{Proof} If $m=1$ the result is trivial, so assume henceforth
that $m\geq 2$, and such a function $f$ is given.
Given integer $N$ we write 
$N=p_1^{R_1-1} p_2^{R_2-1} \dots p_n^{R_N-1}$ where each $R_i\geq 2$.
Moreover we can write $f(p_j)=e^{2i\pi a_j/m}$ for each $1\leq j\leq n$,
where $a_j$ is some integer.
Thus the number of divisors of $N$ for which $f(d)=1$, is exactly
the number of 1's that appear in the expansion
$$
\sum_{d|N} f(d) =
\prod_{j=1}^n \left( 1+ e^{2i\pi a_j/m} + e^{4i\pi a_j/m} + \dots +
e^{2(R_j-1)i\pi a_j/m} \right) ,
$$
which equals 
$$
\#\{ (r_1,\dots ,r_n)\in {\Bbb Z}^n:   \sum_{i=1}^n r_ia_i \equiv 0
\pmod{m},\ \text{with} \ 0\leq r_i\leq R_i-1 \}.
$$
By Lemma 2.3, this is $\geq R_1\dots R_n/2^{m-1} = d(N)/2^{m-1}$,
where $d(N)$ is the number of divisors of $N$. In other words
$$
\sum\Sb d|N \\ f(d)=1\endSb  1 \geq d(N)/2^{m-1} .
$$

Since $\sum_{N\leq x} \sum_{d|N} f(d) = 
\sum_{d\leq x} f(d) [x/d]$, we deduce that 
$$
\sum\Sb d\leq x \\  f(d)=1\endSb \bigg[ \frac xd \bigg] 
= \sum_{N\leq x}\ \sum\Sb d|N \\ f(d)=1 \endSb 1 
\geq \sum_{N\leq x} \frac {d(N)}{2^{m-1}}
= \frac 1{2^{m-1}} \sum_{d\leq x} \bigg[ \frac xd \bigg] ,
$$
from which we deduce the result. 
\enddemo

In the proof of Corollary 2.5 we made extensive use
of Lemma 2.3. However it may be that `typically' Lemma 2.3 is
not sharp.  We now re-examine the lower bounds for $\sum_{d|N,\ f(d)=1} 1$
given above.
For simplicity, we will assume that $m$ is prime and $N$
is squarefree. Suppose that $J$ is the set of integers for which 
$a_j\ne 0 \pmod m$ (where $a_j$ is as defined above). Thus
$$
\sum_{d|N,\ f(d)=1} 1 =
\frac 1m \sum_{\zeta^m=1} \prod_{j=1}^n (1+\zeta^{a_j}) .
$$
We get a contribution of $2^n$ from the $\zeta=1$ term. Otherwise, if
$\zeta=e^{2i\pi k/m}$ then $|1+\zeta^{a_j}|=2|\cos (\pi a_j k/m)|$. 
Therefore
$$
\sum_{d|N,\ f(d)=1} 1 \geq \frac {2^n}m \left( 1 - \sum_{k=1}^{m-1} 
\prod_{j\in J} |\cos (\pi a_j k/m)| \right) 
 \geq \frac {2^n}m \left( 1 - \sum_{k=1}^{m-1} |\cos^{|J|} (\pi k/m)| \right) ,
$$
by an optimization argument. This is $\gg 2^n/m$ if $|J|\gg m^2$;
thus, if a typical integer $N\leq x$ has $\gg m^2$ prime factors for
which $f(p)\ne 1$ then we might expect to improve considerably the lower 
bound in Corollary 2.5.

\head 3. Basic Results on Integral equations \endhead

\subhead 3a. Existence and uniqueness of solutions and first estimates
 \endsubhead

\noindent We begin with the following simple principle
which we shall use repeatedly.

\proclaim{Lemma 3.1} Let $\alpha$ and $\beta$ be two integrable 
functions from $[0,\infty)$ to ${\Bbb R}$.
Suppose that $\alpha (u)\ge 0$ for all $0 \le u \le 1$, and that 
$\beta_0 \ge \beta(u)\ge 0$ for all $u$.  
If $u\alpha(u)\ge (\beta*\alpha) (u)$ then $\alpha(u)\ge 0$ for all $u$.
In particular, if $u\alpha(u)\ge (1*\alpha)(u)$ then $\alpha(u)\ge 0$ 
for all $u$.
\endproclaim

\demo{Proof}  It suffices to show $\alpha(u) \ge 0$ for those points $u\ge 1$ 
satisfying $\alpha(u) /u^{\beta_0-1} \le \alpha(t)/t^{\beta_0-1}$ 
for all $t\le u$.  For such a $u$,
$$
u\alpha(u) \ge \int_0^u \alpha(t) \beta(u-t) dt \ge 
\int_1^u \alpha(t) \beta(u-t) dt \ge \frac{\alpha(u)}{u^{\beta_0-1}}\int_1^u 
t^{\beta_0-1}\beta(u-t)dt.
$$
If $\alpha(u)<0$ then we must have
$$
u^{\beta_0} \le \int_{1}^{u} t^{\beta_0-1}\beta(u-t)dt \le 
\int_1^u \beta_0 t^{\beta_0-1} dt = u^{\beta_0} -1,
$$
which is a contradiction.
\enddemo

The condition that $\beta$ is bounded may be relaxed.  We need only that 
$\beta$ is bounded on closed intervals.  Thus, for example, the result 
holds for any continuous, non-negative function $\beta$.   

Let $\chi$ be an element of $K({\Bbb U})$.
Our first application of this Lemma is to show the existence and 
uniqueness of solutions to the integral equation
(1.5).  To this end, it is useful to define 
$I_0(u) = I_0(u;\chi)=1$, and for $k\ge 1$,
$$
I_k(u)=I_k(u;\chi) = \int\Sb t_1,\ldots t_k \ge 1\\ t_1+\ldots +t_k \le u\endSb
\frac{1-\chi(t_1)}{t_1} \ldots \frac{1-\chi(t_k)}{t_k} dt_1 \ldots dt_k.
$$
Define for all $k\ge 0$,
$$
\sigma_k(u) = \sum_{j=0}^{k} \frac{(-1)^j}{j!} I_j(u;\chi),
\qquad \text{and}\qquad  \sigma_\infty(u) =\sum_{j=0}^{\infty} 
\frac{(-1)^j}{j!} 
I_j(u;\chi).
$$

\proclaim{Lemma 3.2} For all $j\ge 1$,
$$
u I_j(u) = (1*I_j)(u) + j\left( (1-\chi)*I_{j-1}\right)(u).\tag{3.1}
$$
Further $u\sigma_k(u) = (1*\sigma_k)(u)-((1-\chi)*\sigma_{k-1})(u)$
and $u\sigma_\infty(u)=(\sigma_{\infty}*\chi)(u)$.
\endproclaim
\demo{Proof}  Observe that, for $j\ge 1$, 
$$
\align 
(1*I_j)(u) &= \int_1^u \int\Sb t_1, \ldots, t_j \ge 1
\\ t_1+\ldots+t_j \le t\endSb
\frac{1-\chi(t_1)}{t_1} \ldots \frac{1-\chi(t_j)}{t_j} dt_1\ldots dt_j dt\\
&= \int\Sb t_1, \ldots, t_j \ge 1\\ t_1+\ldots+t_j \le u \endSb
\frac{1-\chi(t_1)}{t_1} \ldots \frac{1-\chi(t_j)}{t_j} (u-t_1-\ldots -t_j)
dt_1 \ldots dt_j \\
&= u I_j(u) - j \int_1^u (1-\chi(t_1)) \int\Sb t_2, \ldots, t_j \ge 1\\ 
t_2+\ldots t_j \le u-t_1\endSb \frac{1-\chi(t_2)}{t_2} \ldots 
\frac{1-\chi(t_j)}{t_j} dt_1 \ldots dt_j\\
&= uI_j(u) - j ((1-\chi)*I_{j-1})(u),\\
\endalign
$$
and (3.1) follows.  Multiply both sides of (3.1) by $(-1)^j/j!$ and 
sum from $j=1$ to $k$.  This gives 
$$
\align
u\sum_{j=1}^{k} \frac{(-1)^j}{j!} I_j(u) &= \sum_{j=1}^{k} 
\frac{(-1)^j}{j!} (1*I_j)(u) + \sum_{j=1}^{k} \frac{(-1)^{j}}{(j-1)!}
((1-\chi)*I_{j-1})(u) \\
&= \sum_{j=1}^{k} \frac{(-1)^j}{j!} (1*I_j)(u) - 
\sum_{j=0}^{k-1} \frac{(-1)^j}{j!} ((1-\chi)*I_j)(u).\\
\endalign
$$
Adding $u$ to both sides we get 
$u\sigma_k = 1*\sigma_k - (1-\chi)*\sigma_{k-1}$.  Summing from $j=1$ to 
$\infty$ (instead of $1$ to $k$) gives $u\sigma_{\infty} = 
\sigma_{\infty}*\chi$. 
\enddemo

\proclaim{Theorem 3.3} For a given $\chi \in K({\Bbb U})$, 
there exists a unique 
solution to the integral equation (1.5).  In fact, $\sigma =
\sigma_{\infty}$ is this 
unique solution, and satisfies $|\sigma(u)|\le 1$ for all $u$.
\endproclaim

\demo{Proof}  By definition $\sigma_{\infty} (u)=1$ for $0\le u\le 1$. 
Since $u\sigma_\infty=\sigma_\infty*\chi$, we see that $\sigma_\infty$ is
a solution to (1.5).  We now show that it is unique.  Let $\sigma$ be another 
solution to (1.5) and put $\alpha(u)= - |\sigma(u)-\sigma_\infty(u)|$.  Note 
that $\alpha(u)=0$ for $0\le u\le 1$ and that
$$
u\alpha(u) = -\biggl|\int_0^u (\sigma(t)-\sigma_\infty(t))\chi(u-t)dt\biggr|
\ge -\int_0^u |\sigma(t)-\sigma_\infty(t)|dt =\int_0^u \alpha(t)dt.
$$
Lemma 3.1 shows that $\alpha(u)\ge 0$ always, whence $\sigma=\sigma_\infty$.  

To show that the unique solution $\sigma$ satisfies $|\sigma(u)|\le 1$ 
for all $u$, we take $\alpha (u)=1-|\sigma(u)|$.   Again $\alpha(u)=0$ for 
$0\le u\le 1$, and 
$$
u\alpha(u) = \int_0^u dt - \biggl|\int_0^u \sigma(t)\chi(u-t)dt\biggr| \ge
\int_0^u (1-|\sigma(t)|)dt =\int_0^u \alpha(t)dt.
$$
Thus $\alpha(u)\ge 0$ for all $u$, by Lemma 3.1, and the proof is complete.
\enddemo

\proclaim{Lemma 3.4}  Let $\chi$ and ${\hat \chi}$ be two elements of 
$K({\Bbb U})$, and let $\sigma$ and $\hsigma$ 
be the corresponding solutions to (1.5).  Then $\sigma(u)$ equals 
$$
\hsigma(u) + \sum_{j=1}^{\infty} \frac{(-1)^j}{j!} \int\Sb 
t_1,\ldots ,t_j \ge 1\\ t_1+\ldots +t_j \le u \endSb
\frac{\hchi(t_1)-\chi(t_1)}{t_1}
\ldots \frac{\hchi(t_j)-\chi(t_j)}{t_j} \hsigma(u-t_1-\ldots-t_j)dt_1\ldots 
dt_j. \tag{3.2}
$$
Consequently, if $|\chi(t)-\hchi(t)|\le \chi_0$ for all $t$ then $|\sigma(u)-\hsigma(u)|\le u^{\chi_0}-1$, for all $u\ge 1$.
\endproclaim

\demo{Proof}  Let ${\hat I}_k$ be the integral corresponding to $\hchi$. Writing $1-\chi$ as $(1-\hchi) + (\hchi - \chi)$ in the definition of $I_k$ we deduce
that 
$$
I_k(u) = \sum_{j=0}^{k}\binom{k}{j} \int\Sb t_1,\ldots,t_j\ge 1\\
t_1+\ldots +t_j\le u\endSb \frac{\hchi(t_1)-\chi(t_1)}{t_1}
\ldots \frac{\hchi(t_j)-\chi(t_j)}{t_j} {\hat I}_k(u-t_1-\ldots t_j)
dt_1\ldots dt_j.
$$
Multiply both sides by $(-1)^k/k!$, and sum from $k=0$ to $\infty$.
This proves (3.2).  

If $|\chi(t)-\hchi(t)|\le \chi_0$ then, from (3.2) and the fact that 
$|\hsigma(t)|\le 1$ always, we obtain for $u\ge 1$
$$
|\sigma(u)-\hsigma(u)| \le \sum_{j=1}^{\infty} \frac{1}{j!} 
\biggl(\int_1^u \frac{\chi_0}{t}dt\biggr)^j = \sum_{j=1}^{\infty} \frac{
(\chi_0\log u)^j}{j!} = u^{\chi_0}-1.
$$
This completes the proof of Lemma 3.4.

\enddemo

\proclaim{Lemma 3.5}  Suppose $\chi \in K({\Bbb U})$ is given
 and let $\sigma$ be the corresponding solution to (1.5).
Then 
$$
A(v):=\frac{1}{v} \int_0^v |\sigma(t)| dt
$$
is a non-increasing function of $v$.  Hence, for all $u\ge v$,
$$
|\sigma(u)| \le A(v) = \frac{1}{v}\int_0^v |\sigma(t)|dt. 
$$
\endproclaim
\demo{Proof}  From (1.5), we have $|\sigma(u)|\le A(u)$ for all $u$.  
Differentiating the definition of $A(v)$, we have
 $A^{\prime}(v) = |\sigma(v)|/v - A(v)/v \le 0$, so $A(v)$ is non-increasing
and therefore $|\sigma (u)| \leq A(u) \leq A(v)$ if $u\geq v$.
\enddemo

\subhead 3b.  Inclusion-Exclusion inequalities \endsubhead

\noindent Our formula for $\sigma_\infty (=\sigma)$ looks 
like an inclusion-exclusion type identity.  For a real-valued
function $\chi$, we now show how to obtain inclusion-exclusion inequalities
for $\sigma$. 

\proclaim{Proposition 3.6}  Suppose $\chi \in K([-1,1])$ is given.  
For all integers $k\ge 0$, and all $u\ge 0$,  
$(-1)^{k+1}(\sigma(u)-\sigma_k(u)) \ge 0$.  Thus 
$\sigma_{2k+1}(u)\le \sigma(u)\le \sigma_{2k}(u)$.
\endproclaim

\demo{Proof} From Lemma 3.2 we know that 
$u\sigma_k= 1*\sigma_k -(1-\chi)*\sigma_{k-1}$, and clearly 
$u\sigma = 1*\sigma- (1-\chi)*\sigma$.  Subtracting these identities we get
$u(\sigma-\sigma_k)= 1*(\sigma-\sigma_k)- (1-\chi)*(\sigma-\sigma_{k-1})$.

Put $\alpha_k(u)=(-1)^{k+1}(\sigma(u)-\sigma_k(u))$ so that the above relation
may be rewritten as 
$$
u\alpha_k = 1*\alpha_k + (1-\chi)*\alpha_{k-1}. \tag{3.3}
$$  
We will show that $\alpha_k(u)\ge 0$ always by induction on $k$.  Since 
$\sigma_0=1$, the case $k=0$ follows from Theorem 3.3. Suppose that 
$\alpha_{k-1}$ has been shown to be non-negative.  Since $(1-\chi)$ 
is always non-negative it follows from (3.3) that 
$u\alpha_k(u) \ge (1*\alpha_k)(u)$.  Clearly $\alpha_k(u)=0$ for $0\le u\le 1$.
Lemma 3.1 now shows that $\alpha_k(u)\ge 0$ always, completing our proof.

\enddemo

We now develop some inclusion-exclusion type inequalities for the
case when $\chi \in K({\Bbb U})$ is complex-valued.  
To state this, we make the following 
definitions: Put $R_0(u)=C_0(u)=1$ and for $k\ge 1$ put
$$
C_k(u)= C_k(u;\chi) =\int\Sb t_1,\ldots, t_k\ge 1\\ t_1+\ldots +t_k\le u\endSb
\frac{|\IM \chi(t_1)|}{t_1}\ldots \frac{|\IM \chi(t_k)|}{t_k}dt_1\ldots dt_k,
$$
and  
$$ 
R_k(u)=R_k(u;\chi)=\int\Sb t_1,\ldots, t_k\ge 1\\ 
t_1+\ldots+t_k\le u\endSb \frac{1-\RE\chi(t_1)}{t_1}\ldots 
\frac{1-\RE\chi(t_k)}{t_k} dt_1\ldots dt_k.
$$

\proclaim{Proposition 3.7} For all $u$, $|\IM \sigma(u)|\le C_1(u)$.  
Let $\hchi =
\RE\chi$ and let $\hsigma$ denote the corresponding solution to (1.5).  
Then for all $u$, $|\RE \sigma(u)- \hsigma(u)| \le C_2(u)/2$.  In particular
$1-R_1(u)-C_2(u)/2\le \RE\sigma(u)\le 1-R_1(u) + (R_2(u)+C_2(u))/2$.
\endproclaim

\demo{Proof}  Observe that 
$$
u|\IM \sigma| = |\IM\sigma *\RE\chi +\RE\sigma*\IM\chi| 
\le |\IM \sigma|*1 + 1*|\IM \chi|.\tag{3.4}
$$
In the same way as we showed $uI_k=1*I_k+k(1-\chi)*I_{k-1}$ (see Lemma 3.2),
it follows that 
$$
uC_k = 1*C_k + k |\IM\chi|*C_{k-1}. \tag{3.5}
$$
Define $\alpha(u) =C_1(u)-|\IM\sigma(u)|$ so that $\alpha(u)=0$ for $u\le 1$.
Taking $k=1$ in (3.5) and subtracting (3.4), 
we get $u\alpha(u) \ge 1*\alpha$.  
Lemma 3.1 shows that $\alpha(u)\ge 0$ always.

Notice that 
$$
u(\RE \sigma(u)-\hsigma(u))= \RE \sigma*\RE\chi - \IM\sigma*\IM\chi - 
\hsigma*\RE\chi
$$
whence, using $|\IM\sigma(t)|\le C_1(t)$,  
$$
u|\RE \sigma(u)-\hsigma(u)| \le |\RE\sigma-\hsigma|*1 + |\IM \chi|*C_1.
\tag{3.6}
$$
Put $\alpha(u)=C_2(u)/2 - |\RE\sigma(u)-\hsigma(u)|$ so that 
$\alpha(u)=0$ for $0\le u\le 1$.  Take $k=2$ in (3.5), divide by 2, and 
subtract (3.6).  This gives $u\alpha(u)\ge 1*\alpha$ so that, by Lemma 3.1,
$\alpha(u)\ge 0$ always.  

By Proposition 3.6 we see that $1-R_1(u)\le \hsigma(u) \le 1-R_1(u)+R_2(u)/2$.
This gives the last assertion of the Proposition.

\enddemo

\head 4.  Proof of the Structure Theorem \endhead

\noindent In this section we discuss the relation between the integral 
equation (1.5) and averages of multiplicative functions.  In particular,
we shall prove the Structure Theorem for the spectrum.

\subhead 4a. Variation of averages of multiplicative functions 
\endsubhead

\noindent In this subsection we establish the following Proposition which 
seeks to show that the average value of a multiplicative function varies
slowly.

\proclaim{Proposition 4.1}  Let $f$ be a multiplicative function with 
$|f(n)|\le 1$ for all $n$.  Let $x$ be large, and suppose $1\le y\le 
{x}$.  Then
$$
\biggl|\frac 1x \sum_{n\le x} f(n) - \frac{1}{x/y} \sum_{n\le x/y} f(n)
\biggr| \ll \frac{\log 2y}{\log x} \exp\biggl(\sum_{p\le x} \frac{|1-f(p)|}{p}
\biggr).
$$
\endproclaim

To prove this Proposition we require a consequence of Theorem 2 of 
Halberstam and Richert [4].  Suppose $h$ is a non-negative 
multiplicative function with $h(p^k)\leq 2\gamma^{k-1}$ for all prime powers $p^k$, for some $\gamma,\ 0<\gamma<2$. 
It follows from Theorem 2 of [4] that 
$$
\sum_{n\leq x} h(n) \leq 
\ \frac{2x}{\log x} \ \sum_{n\leq x} \frac{h(n)}{n} \ 
\left\{ 1 + O\left( \frac{1}{\log x}\right) \right\}. \tag{4.1}
$$
Using partial summation we deduce from (4.1) that for $1\leq y\leq x^{1/2}$,
$$
\align
\sum_{x/y< n\leq x} \frac{h(n)}{n} &\leq \left\{ \frac{1}{\log x} + 
\log \left( \frac{\log x}{\log (x/y)} \right) \right\}
 \ \sum_{n\leq x} \frac{h(n)}{n}  
\left\{ 2 + O\left( \frac{1}{\log x}\right) \right\}\\
&\ll \frac{\log 2y}{\log x} \sum_{n\le x} \frac{h(n)}{n}. 
\tag{4.2}
\\
\endalign
$$
Equipped with (4.1) and (4.2) we proceed to a proof of Proposition 4.1.

\demo{Proof of Proposition 4.1} Since the left side of the
Proposition is trivially $\ll 1$, there is 
nothing to prove if $y>\sqrt{x}$.  Suppose now that $y<\sqrt{x}$. 
Let $g$ be the multiplicative function with 
$g(p^k) =f(p^k)-f(p^{k-1})$ for each prime power.  Then
 $f(n) =\sum_{d|n} g(d)$, and so 
$$
\biggl | \frac 1x \sum_{n\le x} f(n) - \sum_{d\le x} \frac{g(d)}{d}
\biggr| \leq \frac 1x \sum_{d\le x} |g(d)|.
$$
Taking this statement for $x$ and $x/y$,  we get
$$
\biggl| \frac 1x \sum_{n\le x} f(n) - \frac {1}{x/y} \sum_{n\le x/y} f(n) 
\biggr|
\leq \sum_{x/y\le d\le x}\frac{|g(d)|}{d} + \frac 1x \sum_{d\le x} |g(d)|
+ \frac {1}{x/y} \sum_{d\le x/y} |g(d)|.
$$
Since each $|g(p^k)|\leq 2$, it follows from (4.1) and (4.2) that the 
above is 
$$
\align
&\ll \frac{\log 2y}{\log x} \sum_{n\leq x} \frac{|g(n)|}{n} 
\ll \frac{\log 2y}{\log x} \prod_{p\le y} \biggl(1+\frac{|1-f(p)|}{p} +
\frac{2}{p^2} + \frac{2}{p^3} + \ldots \biggr)\\
&\ll \frac{\log 2y}{\log x} \exp\biggl( \sum_{p\leq x} \frac{ |1-f(p)|}{p} 
\biggr) .\\
\endalign
$$
This proves the Proposition. 

\enddemo

\subhead 4b. A useful identity \endsubhead

\proclaim{Lemma 4.2}  Let $f$ be a 
multiplicative function with $|f(p^k)|\le 1$.   Then
$$
\log x \int_0^1 \sum_{n\le x^t} f(n) dt +O\biggl(\frac{x}{\log x} \biggr) 
=\int_0^1 \sum_{n\le x^t} f(n) \sum_{m\le x^{1-t}} \Lambda(m) f(m) dt.
$$
\endproclaim

\demo{Proof} 
Note that 
$$
\align 
\int_0^1 \sum_{n\le x^t} f(n) &\sum_{m\le x^{1-t}} \Lambda(m) f(m) dt 
=\sum_{nm\le x} f(n)f(m) \Lambda(m) \frac{\log (x/nm)}{\log x}\\
&=\sum_{nm\le x} f(mn) \Lambda(m) \frac{\log (x/mn)}{\log x} + 
O\biggl( \sum\Sb mn\le x \\ (m,n) >1 \endSb \Lambda(m) 
\frac{\log (x/mn)}{\log x}\biggr),\\
\endalign
$$
and writing $r=nm$ this is 
$$
\align
&=\sum_{r\le x} f(r) \frac{\log (x/r)}{\log x} \sum_{m|r} \Lambda(m) 
+ O\biggl(\frac{x}{\log x}\biggr)\\
&=\sum_{r\le x} f(r) \log r \frac{\log (x/r)}{\log x}
 +O \biggl(\frac{x}{\log x} \biggr).\\
\endalign
$$

Next observe that 
$$
\align
\sum_{r \le x} f(r) \log r \frac{\log (x/r)}{\log x} &= 
\int_0^1 \sum_{r\le x^t} f(r) \log r dt \\
&=\log x \int_0^1 \sum_{r \le x^t} f(r) dt + O\biggl(\log x \int_0^1 (1-t)x^t 
dt\biggr)\\
&=\log x \int_0^1 \sum_{r\le x^t} f(r) dt + O\biggl(\frac{x}{\log x}\biggr).\\
\endalign
$$
The two identities above establish the Lemma.

\enddemo

As a consequence of Lemma 4.2 we derive a convolution identity for the 
averages of $f$ which will be very useful in our treatment of differential 
delay equations (see the proof of Proposition 1 below).

\proclaim{Proposition 4.3}  Let $f$ be a multiplicative function with 
$|f(p^k)| \le 1$.  Then
$$
\sum_{n\le x} f(n) + O\biggl(\frac{x}{\log x} \exp\biggl(\sum_{p\le x} 
\frac{|1-f(p)|}{p} \biggr)\biggr) = 
\int_0^1 \sum_{n\le x^t} f(n) \sum_{p\le x^{1-t} }  f(p) \log p \,dt.
$$
\endproclaim
\demo{Proof}  Applying Proposition 4.1 we find that 
$$
\sum_{n\le x^t} f(n) = x^{t-1} \sum_{n\le x}f(n) + O\biggl((1-t)\exp\biggl(
\sum_{p\le x} \frac{|1-f(p)|}{p} \biggr)\biggr).
$$
Inserting this in the LHS of Lemma 4.2 we get 
$$
\sum_{n\le x} f(n) + O\biggl(\frac{x}{\log x} \exp\biggl(\sum_{p\le x}
 \frac{|1-f(p)|}{p}\biggr)\biggr) = \int_0^1 \sum_{n\le x^t} f(n) 
\sum_{m \le x^{1-t} } f(m) \Lambda(m)\, dt.
$$
Since 
$$
\sum_{m\le x^{1-t}} f(m) \Lambda(m) = \sum_{p\le x^{1-t}} f(p) \log p 
+ O(x^{(1-t)/2})
$$
and 
$$
\int_0^1 \sum_{n\le x^t} f(n) x^{(1-t)/2} \ll \sqrt{x} \int_0^1 x^{t/2} dt \ll 
\frac{x}{\log x},
$$
we have proved the Proposition.

\enddemo

\subhead 4c. Removing the impact of the small primes\endsubhead

\noindent The main result of this section is the following Proposition which 
separates the contribution of small primes.

\proclaim{Proposition 4.4} Fix $\pi/2\geq \phi>0$.
Suppose that  $f\in {\Cal F}(S)$ where $S\subset {\Bbb U}$ with 
$\text{Ang}(S)\leq \pi/2 -\phi$.
For any $\epsilon \ge \log 2/\log x$, let $g$ be the 
completely multiplicative function with $g(p)=1$ if $p\leq x^\epsilon$,
and $g(p)=f(p)$ otherwise, so that $g\in {\Cal F}(S)$ also. Then
$$
\frac 1x \sum_{n\le x} f(n) 
= \Theta (f,x^\epsilon)\ \frac 1x \sum_{m\le x} g(m) + O_\phi(\epsilon^\eta) 
,\ \ \text{where}\ 
\eta=\eta(\phi)=\frac{\sin \phi}{2\pi} \{ \phi-\sin \phi\} .
$$
\endproclaim

We begin by deriving a weak version 
of Proposition 4.4 as a consequence of Proposition 4.1.  
Using this in conjunction with Lemma 1$'$ we shall prove 
the stronger Proposition 4.4.

\proclaim{Proposition 4.5} For any multiplicative function $f$ 
with $|f(p^k)|\leq 1$ for every prime power $p^k$, 
let $s(f,x):=\sum_{p\leq x} |1-f(p)|/p$.
For any $1> \epsilon \ge \log 2/\log x$, let $g$ be the 
completely multiplicative function with $g(p)=1$ if $p\leq x^\epsilon$,
and $g(p)=f(p)$ otherwise. Then
$$
\frac 1x \sum_{n\le x} f(n) = \Theta (f,x^\epsilon)\ \frac 1x \sum_{m\le x} 
g(m) + O( \epsilon \exp(s(f,x))),
$$
where the implicit constant is absolute.
\endproclaim

\demo{Proof}  Define the multiplicative function $h$ by 
$h(p^k)=f(p^k)-f(p^{k-1})$ if $p\leq x^\epsilon$, and $h(p^k)=0$ otherwise.
Then $f(n)=\sum_{m|n} h(n/m) g(m)$, and so 
$$
\frac 1x \sum_{n\le x} f(n) = \sum_{n\le x} \frac{h(n)}{n} \biggl(\frac nx 
\sum_{m\le x/n} g(m)\biggr). \tag{4.3}
$$
Now, Proposition 4.1 gives 
\footnote{
Strictly speaking the error above must have $\log 2n $ instead of $\log n$; 
but  there is no error in the case $n=1$ and when $n\ge 2$, clearly
$\log 2n \ll \log n$.}
$$
\frac nx \sum_{m\le x/n} g(m) = \frac 1x \sum_{m\le x} g(m) +O\biggl(
\frac{\log n}{\log x} \exp(s(g,x))\biggr).
$$
Using this in (4.3) we obtain
$$
\frac 1x \sum_{n\le x} f(n) = \sum_{n\ge 1} \frac{h(n)}{n}
\frac 1x \sum_{m\le x} g(m) + 
O \biggl( \sum_{n=1}^{\infty} \frac{|h(n)|}{n}
\frac{\log n}{\log x} \exp(s(g,x))\biggr). \tag{4.4}
$$
Since $\sum_{n\ge 1} h(n)/n = \Theta (f,x^\epsilon)$
the main term above corresponds to the main term of the Proposition.
We now show how to handle the error term.
Now $h(n)=0$ if $n$ is divisible by a prime larger than $x^{\epsilon}$ whence
$$
\align
\sum_{n=1}^{\infty} \frac{|h(n)|}{n} \log n &= \sum_{n=1}^{\infty} 
\frac{|h(n)|}{n} \sum_{m|n} \Lambda(m) = 
\sum_{p\le x^{\epsilon} } \log p \sum\Sb n=1 \\ p|n \endSb^{\infty} 
\frac{|h(n)|}{n}   + \sum\Sb p^k \\ k\ge 2 \endSb \log p \sum\Sb n=1 \\ p^k |n
\endSb^{\infty} \frac{|h(n)|}{n}\\
&\ll \sum_{p\le x^{\epsilon}} \frac{\log p}{p} 
\sum_{n=1}^{\infty} \frac{|h(n)|}{n} 
+ \sum\Sb p^k \\ k \ge 2\endSb \frac{\log p}{p^k} \sum_{n=1}^{\infty} 
\frac{|h(n)|}{n} \\
&\ll (\epsilon \log x +1) \exp(s(f,x^{\epsilon})).\\
\endalign
$$
Inserting this in (4.4) we obtain that the error term there is
$$\ll (\epsilon +1/\log x) \exp(s(f,x^{\epsilon}) +s(g,x)) 
\ll \epsilon \exp(s(f,x)).
$$
\enddemo

%


\demo{Proof of Proposition 4.4}  Suppose that $z=e^{2it}$. Then
$1-\text{Re}(z) =2\sin^2t$ and $|1-z|=2|\sin t|$. If we restrict $t$
to the range $\phi\leq t\leq \pi/2$ then we get 
$1-\text{Re}(z) = |1-z|\sin t\geq |1-z|\sin \phi$. Thus we obtain, in the
notation of Proposition 4.5,  
$$
\sum_{p \le x} \frac{1-\text{Re}(f(p))}{p} \geq 
(\sin \phi) \sum_{p \le x} \frac{|1-f(p)|}{p} = s(f,x) \sin \phi .
$$
Now Lemma 1$^{'}$ implies that
$$
\biggl|\sum_{n \le x} f(n)\biggr| \ll_\phi 
x \exp\biggl( - \frac {\sin \phi}\pi ( \phi-\sin \phi ) s(f,x) \biggr),
\tag{4.5}
$$
and similarly,
$$
\biggl|\sum_{n\le x}g(n) \biggr| \ll_\phi x \exp\biggl(-\frac{\sin\phi}{\pi} 
(\phi-\sin\phi)s(g,x)\biggr).
$$
Further,
$$
|\Theta(f,x^{\epsilon})| \ll \exp\biggl(-\sum_{p\le x^{\epsilon}} 
\frac{1-\text{Re }f(p)}{p} \biggr) \ll \exp(-\sin\phi s(f,x^{\epsilon})),
$$
whence
$$
\align
|\Theta (f,x^{\epsilon})| 
\biggl|\sum_{n\le x}g(n) \biggr| &\ll_\phi x \exp\biggl(-\frac{\sin\phi}{\pi} 
(\phi-\sin\phi)(s(f,x^{\epsilon})+s(g,x))\biggr)\\
&\ll_\phi x \exp\biggl( -\frac{\sin\phi}{\pi} 
(\phi-\sin\phi)s(f,x)\biggr).\\
\endalign
$$
Together with (4.5), 
this proves the Proposition in the case $s(f,x) \ge \log 1/\sqrt{\epsilon}$.
The case 
$s(f,x) \le \log 1/\sqrt{\epsilon}$ follows from Proposition 4.5. 
\enddemo

\subhead 4d. Completing the proof of the Structure theorem \endsubhead

\noindent We begin by proving Proposition 1 and its converse.

\demo{Proof of Proposition 1}  
Let $s(u) = [y^{u}]^{-1} \sum_{n\le y^u} f(n)$ so that 
$s(u) =1$ for $u\le 1$.  Proposition 4.3 tells us that 
$$
s(u) =\frac 1u \int_0^u s(u-t) \frac{1}{y^t} \left( \sum_{p\le y^t} f(p)\log p 
\right) dt +O\biggl(\frac{u}{\log y} \biggr).
$$
By the prime number theorem $\vartheta (y^t) = y^t +O(y^t /\log(ey^t))$ 
and so 
$$
s(u)=\frac{1}{u}\int_0^u s(u-t) \chi(t) dt +O\biggl(\frac{u}{\log y}\biggr).
$$
Let $C$ be the implied constant in the above estimate; that is, for 
all $u\ge 1$,
$$
\biggl|s(u) -\frac 1u \int_0^u s(u-t)\chi(t) dt \biggr| \le \frac{Cu}{\log y}.
\tag{4.6}
$$

We will demonstrate that $|\sigma(u)-s(u)| \le 2Cu/\log y$ which proves 
the Proposition.  Put $\alpha(u)= - |\sigma(u)-s(u)| + 2Cu/\log y$.  
Plainly $\alpha(u)\ge 0$ for  $u\le 1$ and note that, using (4.6),
$$
\align
(1*\alpha)(u) &= \frac{Cu^2}{\log y} - \int_0^u |\sigma(u-t)-s(u-t)|dt \\
&\le \frac{Cu^2}{\log y} - \biggl| \int_0^u (\sigma(u-t)-s(u-t))\chi(t) dt
\biggr| \\
&\le \frac{Cu^2}{\log y} -u|\sigma(u)-s(u)| + \frac{Cu^2}{\log y} = 
u\alpha(u).\\
\endalign
$$
By Lemma 3.1, $\alpha(u)\ge 0$ always, proving the Proposition.

\enddemo

\demo{Proof of the converse to Proposition 1}
Let $\chi$ be as in the statement of the 
converse to Proposition 1, and let $\sigma$ 
denote the corresponding solution to (1.5).

Since $\chi$ is measurable and $\chi(t)$ belongs to the convex hull of $S$, 
we can find a  step function $\chi_1$ with the following properties:  
$\chi_1(t) =1$ for $t\le 1$, $\chi_1(t)$ in the convex hull of $S$ 
and $|\chi(t) -\chi_1(t)|\le \epsilon/2$ for almost all $t\in[0,u]$.
\footnote{That is, the inequality is violated only on a set of measure $0$.} 
It is a simple exercise (left to the reader!) that $\chi_1$ exists.  

Next, we choose $y$ large and find $f\in {\Cal F}(S)$ with 
$f(p)=1$ for $p\le y$ and such that if
$$
\psi(t) =  \frac{1}{\vartheta(y^t)} \sum_{p\le y^t} f(p)\log p
$$
then $|\psi(t) -\chi_1(t)| \le \epsilon/2$ for almost all $t \in [0,u]$.  
The existence of $f$ is another straight-forward exercise.  

With this choice, $|\chi(t)- \psi(t)|
\le \epsilon$ for almost all $t \in [0,u]$.  
Let $\tilde{\sigma}$ denote the solution to $u\tilde{\sigma}
(u)= (\tilde{\sigma} *\psi)(u)$ with 
the initial condition $\tilde\sigma(t)=1$ 
for $t\le 1$.  By Proposition 1 we note that for $t\le u$,
$$
\frac{1}{y^t} \sum_{n\le y^t} f(n) = \tilde\sigma(t) 
+ O\Big(\frac{t}{\log y}\Big).
$$
 From Lemma 3.4, we note  that $|\tilde\sigma(t) - \sigma(t)| 
\le t^{\epsilon} -1
\le u^{\epsilon}-1$.  This completes our proof.

\enddemo

We are now in a position to prove the Structure theorem.

\demo{Proof of Theorem 3}  If Ang$(S) =\pi/2$ then $\Gamma(S) = \GT = 
\Lam(S)={\Bbb U}$, and there is nothing to prove.  So we suppose below that 
Ang$(S)<\pi/2$.  

If $z\in \Gamma(S)$ then there 
exist large $x$ and $f\in {\Cal F}(S)$ for which $\frac{1}{x}\sum_{n\le x}
f(n) =z+o(1)$.  Take $y=\exp((\log x)^{\frac23})$ and define $g\in{\Cal F}(S)$
by $g(p)=1$ for $p\le y$, and $g(p)=f(p)$ for $p>y$.  By Proposition 4.4
$\frac{1}{x} \sum_{n\le x} f(n) = \Theta(f,y) \frac 1x \sum_{n\le x} g(n) +o(1)
$.  Take $\chi(t)=1$ for $t\le 1$, and $\chi(t)=\frac{1}{\vartheta(y^t)} 
\sum_{p\le y^t} g(p) \log p$ for $t>1$.  Let $\sigma$ denote the 
corresponding solution to (1.5).  Proposition 1 tells us that $\frac 1x 
\sum_{n\le x} g(n) = \sigma(\frac{\log x}{\log y})+o(1)$.  It follows that 
$$
z=\frac 1x \sum_{n\le x} f(n) +o(1) = \Theta(f,y) \sigma\Big(\frac{\log x}
{\log y}\Big) + o(1).
$$
This shows that $\Gamma(S) \subset \GT \times \Lam(S)$.

Suppose now that $z_{\theta} \in \GT$, and $z_{\sigma} \in \Lam(S)$ are 
given.  Plainly for large $y$ there exists $g\in {\Cal F}(S)$ with $z_\theta
= \Theta(g,y)+o(1)$.  Further suppose $z_\sigma =\sigma(u)$ for 
some $u$, and $\sigma$ a solution to (1.5) for some measurable function 
$\chi$ with $\chi(t)=1$ for $t\le 1$, and $\chi(t)$ in the convex hull of $S$ 
for all $t$.  By Proposition 1 (Converse) we deduce that there exists $h\in
{\Cal F}(S)$ with $h(p)=1$ for $p\le y$ such that 
$z_\sigma =\sigma(u) = \frac{1}{y^u} \sum_{n\le y^u} h(n) + o(1)$.  Define 
$f\in {\Cal F}(S)$ by $f(p) = g(p)$ if $p\le y$, and $f(p)=h(p)$ if $p>y$.  
By Proposition 4.4 it follows that 
$$
\frac{1}{y^u} \sum_{n\le y^u} f(n) = \Theta(g,y) \frac{1}{y^u} 
\sum_{n\le y^u} h(n)  + o(1) = z_\theta z_\sigma +o(1).
$$
Hence $\GT\times \Lam(S) \subset \Gamma(S)$, proving Theorem 3.

\enddemo

\head 5.  Determining the spectrum of $[-1,1]$; Proof of Theorem 1 \endhead

\noindent In this section we shall prove Theorem 1 and Corollary 1.  
In Theorem 3$^{\prime}$ we saw that 
 $\Gamma([-1,1]) = \Lambda([-1,1])$, and we have 
already seen that $\Lambda([-1,1])\supset [\delta_1,1]$.  The following
theorem shows that
$\Lambda([-1,1])\subset [\delta_1,1]$, and more.

\proclaim{Theorem 5.1} Let $\chi \in K([-1,1])$ be given, and 
let $\sigma$ 
denote the corresponding solution to (1.5).  
If $\int_0^u \frac{1-\chi(t)}{t} dt <1$ for all $u$, then $\sigma(u)$ 
is always positive.  On the other hand, if $\int_0^{u_0} \frac{1-\chi(t)}{t} 
dt =1$ for some real number $u_0$, then $\sigma(u) \ge 0$ for 
all $u\le u_0$, and $|\sigma(u)| \le |\delta_1|$ for all $u\ge u_0$.  
Moreover, if $|\sigma(u) - \delta_1| \le \epsilon$ then 
we must have  
$u = (1+1/\sqrt{e}) u_0 +O({\epsilon}^{\frac 14})$ and
$$
\int_0^{u/(1+\sqrt{e})} \frac{1-\chi(t)}{t} dt + \int_{u/(1+\sqrt{e})}^u 
\frac{1+\chi(t)}{t} dt \ll \sqrt{\epsilon}.
$$
\endproclaim

Given Theorem 5.1 we now show how Corollary 1 may be deduced.

\demo{Deduction of Corollary 1}  Given $f\in {\Cal F}([-1,1])$, choose 
$y= \exp( (\log x)^{\frac 23})$, 
and define $g \in {\Cal F}([-1,1])$ by 
$g(p)=1$ for $p\le y$, and $g(p)=f(p)$ for $p>y$.  Define for $t\ge 0$,
$$
\chi(t) = \frac{1}{\vartheta(y^t)} \sum_{p\le y^t} g(p)\log p,
$$
and let $\sigma$ denote the corresponding solution to (1.5).

By Proposition 4.4 (with $S=[-1,1]$, and $\phi=\pi/2$) we have that 
$$
\align
\frac{1}{x} \sum_{n\le x} f(n) &= \Theta(f,y) \frac{1}{x} \sum_{n\le x} g(n) 
+O\Big( \Big(\frac{\log y}{\log x}\Big)^{(\pi/2-1)/(2\pi)}\Big) \\
&= \Theta(f,y) \frac{1}{x} \sum_{n\le x} g(n) +o(1).
\endalign
$$
Appealing now to Proposition 1, this is 
$$
= \Theta(f,y) \Big( \sigma\Big(\frac{\log x}{\log y}\Big) +
O\Big(\frac{\log x}{\log^2 y}\Big) \Big) + o(1) 
= \Theta(f,y) \sigma\Big(\frac{\log x}{\log y}\Big) +o(1).
$$
Since $\Theta(f,y) \in [0,1]$, it follows at once from Theorem 5.1, 
that $\frac 1x \sum_{n\le x} f(n) \ge \delta_1+o(1)$.  Further, if 
equality holds here then we must have $\Theta(f,y)=1+o(1)$, and 
$\sigma(\frac{\log x}{\log y}) = \delta_1 +o(1)$.  The conclusion of the 
corollary now follows upon using our knowledge of 
when equality in Theorem 5.1 can occur.

\enddemo

The remainder of this section will be concerned with the proof of 
Theorem 5.1. Recall from \S 3 the definitions of $I_k(u;\chi)$.  
By Proposition 3.6 with $k=0$ we have $\sigma(u)\ge 1-I_1(u;\chi)$.  
Hence if $I_1(u;\chi)= \int_{0}^{u} \frac{1-\chi(t)}{t} dt<1$ 
for all $u$ then $\sigma(u)>0$ always, which is the 
first case of our Theorem.
So we may suppose that there is a number 
$u_0$ such that $I_1(u_0;\chi)=1$. 
Plainly $\sigma(u) \ge 1-I_1(u;\chi) \ge 1-I_1(u_0;\chi)=0$ if $u\le u_0$.  
Hence it remains to be shown that $|\sigma(u)|\le |\delta_1|$ for all
$u\ge u_0$, and to identify when $\sigma(u)$ is ``close'' to $\delta_1$.

We begin by giving an outline of the underlying ideas of this proof. 
It is helpful first to gain an understanding of the 
extremal function $\rho_{-}(t)$, which we 
discussed briefly in the introduction.  Recall that 
$\rho_{-}(t) =1$ for $t\le 1$, and for $t>1$ is the 
unique continuous solution to the differential-difference equation
$t \rho_{-}'(t) = -2\rho_{-}(t-1)$.  
Alternatively, in terms of integral equations, for $v\ge 1$ we have 
$$
v\rho_{-}(v)  = \int_{v-1}^v \rho_{-}(t) \ dt - 
\int_0^{v-1} \rho_{-} (t) \ dt.
$$ 
By integrating $\rho_{-}'(t)$ appropriately, and 
using the differential-difference relation, we obtain that
$$
\rho_{-}(t) = 1- 2\log t, \qquad \text{for}  \ \ \ 1\le t\le 2,
$$
and that 
$$
\rho_{-}(t) = 1-2\log t + \gamma(t), \qquad \text{for} \ \ \ 2\le t\le 3,
$$
where we put $\gamma(t)=0$ for $t\le 2$, and define for $t\ge 2$ 
$$
\gamma(t) = 4 \int_2^t \frac{\log (v-1)}{v} dv.  \tag{5.1}
$$
Notice that $\rho_{-}(t) \ge 0$ for $t\le \sqrt{e}$, $\rho_{-}(\sqrt{e})=0$,
and that $\rho_{-}(t) \le 0$ for $\sqrt{e}\le t\le 3$. \footnote { 
In fact, $\rho_{-}(t)\le 0$ for all $t\ge \sqrt{e}$ but we do not need this 
fact.}  Hence note that 
$$
\align
\delta_1 &= \rho_{-}(1+\sqrt{e}) = \frac{1}{(1+\sqrt{e})} \Big(\int_{\sqrt{e}}
^{1+\sqrt{e}}  \rho_{-}(t) dt  -  \int_0^{\sqrt{e}}\rho_- (t) dt \Big) \\
&= - \frac{1}{1+\sqrt{e}} \int_0^{1+\sqrt{e}} |\rho_-(t)| dt,\\
\endalign
$$
or alternatively,
$$
|\delta_1|= |\rho_{-}(1+\sqrt{e})| = \frac{1}{1+\sqrt{e}} 
\int_0^{1+\sqrt{e}} |\rho_-(t)| dt. \tag{5.2}
$$
This identity lies at the heart of our proof.

Suppose for simplicity that 
$u>u_0(1+1/\sqrt{e})$;  we seek to show that $|\sigma(u)| \le |\delta_1|$.  
By Lemma 3.5 we note that 
$$
\align
|\sigma(u)| &\le \frac{1}{u_0(1+1/\sqrt{e})} \int_0^{u_0(1+1/\sqrt{e})} 
|\sigma(t)| \ dt \\
&= 
 \frac{1}{u_0(1+1/\sqrt{e})} \biggl(
\int_0^{u_0} \sigma(t) dt + \int_{u_0}^{u_0(1+1/\sqrt{e})} |\sigma(t)| dt 
\biggr).\\
\endalign
$$
Our idea is essentially to 
compare $|\sigma(t)|$ with $|\rho_{-}(t\sqrt{e}/u_0)|$.  We shall show 
that $|\sigma(t)|$ is smaller on average than $|\rho_-(t\sqrt{e}/u_0)|$.  
 From this and the above inequality it would follow that 
$|\sigma(u)| \le \frac{1}{(1+\sqrt{e})} \int_0^{1+\sqrt{e}} |\rho_{-}(t)|dt$,
and from (5.2) the result follows.

In order to carry this out, we introduce the parameters
$$
\lambda = I_1(u_0(1-1/\sqrt{e}),\chi), \qquad \text{ and} \qquad 
\tau=I_1(u_0/\sqrt{e},\chi), \tag{5.3}
$$
which satisfy $0\le \lambda \le \tau \le 1$.  In \S 5b we present an argument
which maximizes $\int_0^{u_0} \sigma(t) \ dt$ under the constraint (5.3).  
We show there that 
$$
\frac{1}{u_0} \int_0^{u_0} \sigma(t) dt \le 2-\frac{2}{\sqrt{e}} - 
E_1(\lambda,\tau) = \frac{1}{\sqrt{e}}
\int_0^{\sqrt{e}} \rho_-(t)dt - E_1(\lambda,\tau), 
$$
where $E_1(\lambda,\tau)$ is an explicit 
non-negative function of $\tau$ and
$\lambda$ (see Corollaries 5.6 and 5.7 below).  

For $u_0 \le t \le u_0(1+1/\sqrt{e})$, we use the 
inclusion-exclusion inequalities of Proposition 3.6
to obtain estimates of the form 
$$
|\sigma(t)| \le |\rho_-(\sqrt{e}t/u_0)| + E_2(\lambda,\tau,t/u_0),
$$
for some non-negative function $E_2(\lambda,\tau,t/u_0)$.  
The key is to obtain very precise bounds for 
$E_2(\lambda,\tau,t/u_0)$ such that
$$
\frac{1}{u_0(1+1/\sqrt{e})} \int_{u_0}^{u_0(1+1/\sqrt{e})} 
E_2(\lambda,\tau,t/u_0)dt \le \frac{E_1(\lambda,\tau)}{1+1/\sqrt{e}}.
$$
In fact, we shall see that equality above holds only when 
$\lambda = \tau =0$.  
Combining this with our bound on $\int_0^{u_0} \sigma(t) \ dt$, we 
shall have 
shown that $|\sigma(t)|$ is smaller than $|\rho_{-}(t\sqrt{e}/u_0)|$ 
on average; as desired.

\subhead 5a. Preliminaries \endsubhead

\noindent Throughout $\lambda$ and $\tau$ are as in (5.3).  
We shall find it useful to 
consider the function $\hchi(t)=\chi(t)$ if $t\le u_0$
and $\hchi(t)=1$ for $t >u_0$.  Let $\hsigma$ denote the 
corresponding solution to (1.5).  Below,  ${\hat I}_k(u)$ 
will denote $I_k(u;\hchi)$.  Note that ${\hat I}_1(u) \le 1$ 
for all $u$, and 
so by Proposition 3.6, it follows that $\hsigma(u) \ge 1- {\hat I}_1(u) \ge 0$ 
always.  

\proclaim{Lemma 5.2} In the range $u_0 \le u \le 2u_0$ we have
$$
\max\Big(-2\log \frac{u}{u_0}, -2\log \frac{u}{u_0} +\frac{{\hat I}_2(u)}{2}
- \frac{{\hat I}_3(u)}{6} \Big)\le \sigma(u) \le \frac{{\hat I}_2(u)}{2} .
$$
\endproclaim

\demo{Proof} By Lemma 3.4 we see that in this range
$$
\sigma(u) = \hsigma(u) - \int_{u_0}^{u} \frac{1-\chi(t)}{t} \hsigma(u-t) dt.
$$
Since $0\le 1-\chi(t)\le 2$, and $0\le \hsigma(u-t)\le 1$, it follows 
that 
$$
-2\log \frac{u}{u_0} + \hsigma(u) \le \sigma(u) \le \hsigma(u).
$$
Moreover $\hchi$ has been designed so that ${\hat I}_1(u)=1$ for all 
$u\geq u_0$. Therefore, by Proposition 3.6, 
we see that $\hsigma(u)\le 1-{\hat I}_1(u) +{\hat I}_2(u)/2
={\hat I_2}(u)/2$, and also that $\hsigma(u) \ge 
\max(0, 1- {\hat I}_1(u) + {\hat I}_2(u)/2 -{\hat I}_3(u)/6) 
= \max(0,{\hat I}_2(u)/2 - {\hat I}_3(u)/6)$.  The Lemma follows.

\enddemo

In order to use Lemma 5.2 successfully, we require estimates 
for ${\hat I_2}(u)$, and ${\hat I_3}(u)$.  We 
develop these in the next two Lemmas.

\proclaim{Lemma 5.3}  In the range $u_0\le u\le u_0(1+1/\sqrt{e})$ we have
$$
\align
{\hat I}_2(u) &\le \min\Big(1,
2\gamma\Big(\frac{u\sqrt{e}}{u_0}\Big) 
+ 2\lambda-\tau^2 +2(\tau-\lambda) {\hat I}_1\left(u-u_0(1-1/\sqrt{e})\right)
\Big) 
\\
&\le 
\min\Big(1, 2\gamma\Big(\frac{u\sqrt{e}}{u_0}\Big) + 2\tau-\tau^2\Big),
\\
\endalign
$$
and 
$$
{\hat I}_2(u) \ge 2\gamma\Big(\frac{u\sqrt{e}}{u_0}\Big).
$$
\endproclaim

\demo{Proof}  Write $\hchi_1(t) = (1-\hchi(t))/t$.  We define $\psi_0(t) =
0$ for $t\le u_0/\sqrt{e}$ and $\psi_0(t) = \hchi_1(t)$ for $t>u_0/\sqrt{e}$.
Define $\psi_1(t)= 0$ if $t\le u_0/\sqrt{e}$ or if $t>u_0$ and 
$\psi_1(t) = 2/t$ if $u_0/\sqrt{e}< t\le u_0$.  Notice that 
$\psi_0(t)\leq \psi_1(t)$ for all $t\leq u_0$, and so
$(1*\psi_0*\psi_0)(u)\leq (1*\psi_1*\psi_1)(u)=2\gamma(u\sqrt{e}/u_0)$
for $u\leq u_0(1+1/\sqrt{e})$.  

Since $\hchi_1 \ge 0$ and ${\hat I}_1=1*\hchi_1 \le 1$, 
we get ${\hat I}_2 = 1*\hchi_1*\hchi_1 \le 1*\hchi_1 \le 1$.  Hence
$$
\align
{\hat I}_2 &= 1*\hchi_1*\hchi_1 = 1*\psi_0*\psi_0 
+ 1*(\hchi_1-\psi_0)*(\hchi_1+\psi_0)\\
&\le 2\gamma (u\sqrt{e}/u_0) + 1*(\hchi_1-\psi_0)*(\hchi_1+\psi_0).\\
\endalign
$$
Now, for $u$ in this range,
$$
\align
1*(\hchi_1-\psi_0)*&(\hchi_1+\psi_0) = \int_{1}^{u_0/\sqrt{e}} \hchi_1(t_1)
\Big(\int_{1}^{u-t_1} (\hchi_1(t_2)+\psi_0(t_2))dt_2\Big)dt_1\\
&\le \Big(\int_1^{u_0(1-1/\sqrt{e})}\hchi_1(t_1)dt_1\Big)\Big(\int_1^u 
(\hchi_1(t_2)+\psi_0(t_2)) dt_2 \Big)\\
&\ \ + \Big(\int_{u_0(1-1/\sqrt{e})}^{u_0/\sqrt{e}}
\hchi_1(t_1) dt_1\Big)\Big(
 \int_1^{u-u_0(1-1/\sqrt{e})} (\hchi_1(t_2)+\psi_0(t_2))dt_2\Big)\\
&= \lambda(1+1-\tau)+(\tau-\lambda)(2{\hat I}_1(u-u_0(1-1/\sqrt{e})) -\tau).
\\
\endalign
$$
This shows the middle upper bound of the Lemma, from which the last upper 
bound of the lemma follows as $\lambda\le \tau$ and ${\hat I}_1(u-u_0
(1-1/\sqrt{e})) \le 1$.

Observe that 
$$
\align
{\hat I}_2 &= 1*\hchi_1*\hchi_1 = 1*\psi_1*\psi_1+ 1*(\hchi_1-\psi_1)
*(\hchi_1+\psi_1)\\
&=2\gamma(u\sqrt{e}/u_0) + 1*(\hchi_1-\psi_1)*(\hchi_1+\psi_1).\\
\endalign
$$
We now show that $1*(\hchi_1-\psi_1)\ge 0$ which would show the lower bound.
If $t>u_0$ then $(1*(\hchi_1-\psi_1))(t) 
= {\hat I}_1(u_0)-\int_{u_0/\sqrt{e}}^{u_0} 2dv/v =0$.  If 
$t\le u_0/\sqrt{e}$ then $(1*(\hchi_1-\psi_1))(t)=(1*\hchi_1)(t)\ge 0$.
Lastly if $u_0/\sqrt{e}\le t\le u_0$, then
$$
1*(\hchi_1-\psi_1) = \int_0^{t} (\hchi_1(v)-\psi_1(v))dv =\int_t^{u_0}
\biggl(\frac{2}{v}-\frac{1-\hchi(v)}{v}\biggr)dv\ge 0.
$$

\enddemo

\proclaim{Lemma 5.4} If $u\le u_0(1+1/\sqrt{e})$ then ${\hat I}_3(u)\le 
3 \lambda {\hat I_2}(u) + 3\tau^2$.
\endproclaim
\demo{Proof}  By definition 
$$
{\hat I}_3(u) = \int\Sb t_1,t_2,t_3\ge 1\\ t_1+t_2+t_3\le u\endSb 
\frac{1-\hchi(t_1)}{t_1}\frac{1-\hchi(t_2)}{t_2} \frac{1-\hchi(t_3)}{t_3}
dt_1dt_2dt_3.
$$
Since $u\le u_0(1+1/\sqrt{e})$ it follows that either one of $t_1$, $t_2$, 
$t_3$ is $\le u_0(1-1/\sqrt{e})$ or at least two of $t_1$, $t_2$, $t_3$ must be
$\le u_0/\sqrt{e}$.  The first case contributes $\le 3\lambda {\hat I}_2(u)$ 
and the second contributes $\le 3\tau^2$.  
\enddemo

\subhead 5b. Bounding $\int_0^{u_0} |\sigma(u)| du$ \endsubhead

\noindent Note that if $u\le u_0$ then $\sigma(u)=\hsigma(u) \ge 
1-{\hat I}_1(u) \ge 0$.  Since $\sigma(u)\le 1-I_1(u)+I_2(u)/2$ 
by Proposition 3.6,  
we obtain
$$
\frac{1}{u_0}\int_0^{u_0}|\sigma(u)|du = \frac{1}{u_0} \int_0^{u_0} 
\sigma(u)du \le 1 - \frac{(1*I_1)(u_0)}{u_0}+\frac{(1*I_2)(u_0)}{2u_0}.
$$
Observe that $1*I_2 = 1*(1*\chi_1*\chi_1) = 1*\chi_1*1*\chi_1 =I_1*I_1$.
Hence
$$
\frac{1}{u_0} \int_0^{u_0} |\sigma(u)|du \le \frac{1}{2} +\frac{1}{2u_0} 
\left((1-I_1)*(1-I_1)\right)(u_0). \tag{5.4}
$$
If we had lower bounds for $I_1(t)$ then we could use those in (5.4) to 
get an upper bound on $\int_0^{u_0} |\sigma(u)|du$.  

\proclaim {Lemma 5.5} For $0 \le u \le u_0$, $I_1(u) \ge \psi_1(u)$
where
$$
\psi_1(u) := 
\cases
\max\left(0,\lambda + 2 \log((u/u_0)/(1-1/\sqrt{e}))\right) 
&\text{if $0 \le u \le (1-1/\sqrt{e})u_0$}\cr
\max(\lambda, \tau+2\log (\sqrt{e}u/u_0)) &\text{if $(1-1/\sqrt{e})u_0 
\le u \le u_0/\sqrt{e}$}\cr
\max(\tau, 1+ 2\log (u/u_0)) &\text{if $u_0/\sqrt{e} \le u \le u_0$.}
\endcases
$$
Note that if  $\tau \ge 2\log 2 -1$ then
$$
\psi_1(u) \ge \psi_2(u) := 
\cases
\max(0,2\log (2u/u_0)) &\text{if $0 \le u \le e^{\tau/2}u_0/2$}\cr
\max (\tau,1+2\log (u/u_0)) &\text{if $e^{\tau/2}u_0/2 \le u \le u_0$.}
\endcases
$$
\endproclaim

\demo{Proof} We denote each of the above ranges in the definition 
of $\psi_1$ by $[u_1,u_2]$.  Since $I_1$ is  non-decreasing, we know that 
that $I_1(u)\geq I_1(u_1)$ (which gives the lower bounds
$0$, $\lambda$, $\tau$ respectively).  Further $I_1(u) \ge I_1(u_2) - 
\int_u^{u_2} 2dt/t = I_1(u_2) -2\log (u_2/u)$, which gives the other lower
bound for that range.

The bounds $\psi_1(u)\geq \psi_2(u)$ follow from the definitions.
\enddemo

We could plug in the lower bound $\psi_1$ in (5.4) to obtain an upper 
bound for $\int_0^{u_0} |\sigma(u)|du$.  However the resulting 
expression is complicated and we prefer to obtain simpler, but still 
sufficiently strong, bounds.  We first use  Lemma 5.5 to deal
with the simpler case when $\tau \ge 2\log 2 -1$.  

\proclaim{Corollary 5.6}  Suppose $\tau \ge 2\log 2 -1$, then
$$
\frac{1}{u_0}\int_0^{u_0} |\sigma(u) | du  
\le 1+ e^{\tau/2} \Big(1-\frac{2}{\sqrt{e}}\Big).
$$
\endproclaim

\demo{Proof}  Since $I_1(t) \ge \psi_2(t)$, and as $\psi_2(t)=0$ for 
$t\le u_0/2$ we see that 
$$
\frac{1}{u_0} \int_0^{u_0} |\sigma(u)|du \le \frac{1}{2} + \frac{1}{2u_0} 
((1-\psi_2)*(1-\psi_2))(u_0) = 1 - \frac{(1*\psi_2)(u_0)}{2u_0}.
$$
The corollary follows upon calculating $(1*\psi_2)(u_0)$.
\enddemo

\proclaim{Corollary 5.7}  We have
$$
\frac{1}{u_0}\int_0^{u_0} |\sigma(u)|du 
\le 2-\frac{2}{\sqrt{e}} - \frac{\tau^2}{2\sqrt{e}} - 
\cases
0 &\text{if $\tau \le  2\log 2 -1$} \cr
\lambda/12  &\text{if $\tau \le 3/10$.} \cr
\endcases 
$$
\endproclaim

\demo{Proof}  Throughout the calculations in this proof we make
use of the hypothesis that $0\leq \lambda\leq \tau \leq  2\log 2 -1$.
 From (5.4), we know that the desired integral is $\le 1/2 
+((1-\psi_1)*(1-\psi_1)(u_0))/(2u_0)$.

With some calculation we verify that
$$
1- \frac{(1*\psi_1)(u_0)}{u_0} 
= 2\left( 2-e^{-\lambda/2} - e^{-1/2} \left(e^{\tau/2}
+e^{(\lambda-\tau)/2} - e^{-\lambda/2}\right)\right). \tag{5.5}
$$
Since $\tau\le 2\log 2 -1$ we have $\psi_1(u)= \lambda$ for $u_0(1-1/\sqrt{e})
\le u\le u_0/2$.  Moreover, if $u_0-u >u_0/\sqrt{e}$ and $\psi_1(u_0-u) 
\neq \tau$ then, by definition,
we must have 
$u \le u_0(1-e^{(\tau-1)/2})\le u_0(1-1/\sqrt{e}) e^{-\lambda/2}$,
so that $\psi_1(u)=0$. Thus
$$
\frac{\psi_1*\psi_1(u_0)}{2u_0} 
=\frac{\tau}{u_0} \int_0^{u_0(1-1/\sqrt{e})} \psi_1(u) du 
+ \frac{\lambda}{u_0} \int_{u_0(1-1/\sqrt{e})}^{u_0/2} \psi_1(u_0-u) du.
$$
With some calculation one can verify that this equals
$$
\tau \left( 1-\frac{1}{\sqrt{e}}\right) (2e^{-\lambda/2} -2 +\lambda) + \lambda^2 \left( \frac{1}{\sqrt{e}}- \frac 12\right)+ 
\frac{\lambda}{\sqrt{e}} (2e^{(\lambda
-\tau)/2} -2 + (\tau-\lambda)). \tag{5.6}
$$
Thus we know $(5.4) \le (5.5)+(5.6)$.  We now obtain some simple upper
bounds for the expressions in (5.5) and (5.6).

By observing that 
$1-e^{-\lambda/2} \le \lambda/2 -\lambda^2/8 +\lambda^3/48$,
and that
$$
e^{\tau/2} +e^{(\lambda-\tau)/2} -e^{-\lambda/2} \ge 1+\lambda +
\frac{\tau^2}{8} + \frac{\tau^3}{48} + \frac{(\lambda -\tau)^2}{8} + 
\frac{(\lambda -\tau)^3}{48} - \frac{\lambda^2}{8} \ge 1+ \lambda + 
\frac{\tau^2}{4}  - \frac{\tau \lambda}{4} ,
$$
we deduce the upper bound
$$
(5.5) \le 2-\frac{2}{\sqrt{e}} -\frac{\tau^2}{2\sqrt{e}} - 
\lambda\biggl( \frac{2}{\sqrt{e}} -1 -\frac{\tau}{2\sqrt{e}} \biggr) - 
\lambda^2\biggl( \frac 14 -\frac{\lambda}{24}\biggr). \tag{5.7}
$$

Since $2e^{-\xi/2} -2 +\xi \le \xi^2/4$ for all $\xi \ge 0$, we get
$$
\align
(5.6)
&\le \frac{\tau  \lambda^2}{4}(1-1/\sqrt{e})
+\lambda^2(1/\sqrt{e}-1/2) + \frac{\lambda}{4\sqrt{e}} (\tau-\lambda)^2 \\
&\le \lambda \frac{\tau^2}{4\sqrt{e}} + \lambda^2 \biggl(\frac{1}{\sqrt{e}}
 -\frac 12 +\frac{\lambda}{4\sqrt{e}}\biggr). \\
\endalign  
$$
Combining this upper bound with (5.7), we get 
$$
(5.4)\le 2\biggl(1-\frac{1}{\sqrt{e}}\biggr) -\frac{\tau^2}{2\sqrt{e}}
 -\lambda \biggl( \frac{2}{\sqrt{e} } -1 - 
\frac{\tau}{2\sqrt{e}} - \frac{\tau^2}{4\sqrt{e}}\biggr) 
 - \lambda^2
\biggl(\frac 34 - \frac{1}{\sqrt{e}} - \frac{\lambda}{24} \biggl( 
1+\frac{6}{\sqrt{e}}\biggr) \biggr). 
$$

We deduce our result by noting that
$$
\frac 34 - \frac{1}{\sqrt{e}} - \frac{\lambda }{24} \biggl( 1+\frac{6}
{\sqrt{e}}\biggr) \ge \frac 34 - \frac 1{\sqrt{e}} - \frac{(2\log 2 -1)}
{24} \biggl(1 + \frac{6}{\sqrt{e}} \biggr) > 0,
$$
and
$$
\frac{2}{\sqrt{e}} -1 - \frac{\tau}{2\sqrt{e}} -\frac{\tau^2}{4\sqrt{e}} \ge
\cases
0 &\text{if $\tau \le  2\log 2 -1$} \cr
1/12  &\text{if $\tau \le 3/10$.} \cr
\endcases
$$
\enddemo

\subhead 5c. Proof of Theorem 5.1 for large $\tau (\ge 29/100)$ \endsubhead

\noindent Define $\alpha =\exp(|\delta_1|/2-10^{-6})$.  By Lemmas 5.2 and 
5.3 we know that for $(1+1/\sqrt{e})u_0\ge u\ge u_0$, 
$|\sigma(u)|\le 
\max(1/2,2\log(u/u_0))$.  It follows that Theorem 1 holds in the range
$0 \le u \le \alpha u_0$.  If $u > \alpha u_0$ then by Lemma 3.5
$$
|\sigma(u)|\le \frac{1}{\alpha u_0}\int_0^{\alpha u_0} |\sigma(t)|dt   
\tag{5.8}
$$
so it suffices to show that this integral is $<|\delta_1| -10^{-6}$
to complete the proof of Theorem 5.1 in this range of $\tau$. 

We estimate this integral by bounding it 
in various ranges using several results
from previous sections:\  First, we bound $\int_0^{u_0} |\sigma(t)| dt$
by Corollary 5.6 when $2\log 2 -1 \le \tau \le  1$, and by the first part
of Corollary 5.7 when $29/100 \le \tau \le  2\log 2 -1$.

Second, since $|\sigma(u)|\le \max(1/2,2\log (u/u_0))$ for 
$2u_0/\sqrt{e} \le u \le u_0\alpha$,  we get the bound
$$
\frac{1}{u_0}\int_{2u_0/\sqrt{e}}^{\alpha u_0} 
|\sigma(t)|dt \le \int_{2/\sqrt{e}}^{e^{1/4}} 
\frac{dt}{2} + \int_{e^{1/4}}^{\alpha} 2\log t dt. 
$$

Finally, we have 
$|\sigma(u)|\le \max(\tau -\tau^2/2,2\log (u/u_0))$ 
for $u_0 \le u \le 2u_0/\sqrt{e}$,  by Lemmas 5.2 and 5.3 
and since $\gamma(u)=0$
in this range. Thus if $\exp(\tau/2-\tau^2/4)\ge 2/\sqrt{e}$ 
(which happens when $\tau \ge 0.5231 \ldots$) then 
$$
\frac{1}{u_0}\int_{u_0}^{2u_0/\sqrt{e}} |\sigma(t)| dt 
\le (2/\sqrt{e}-1)(\tau-\tau^2/2) 
$$
whereas if $\exp(\tau/2-\tau^2/4) \le 2/\sqrt{e}$ 
(which happens when $\tau \le 0.5231 \ldots$)  then
$$
\int_{u_0}^{2u_0/\sqrt{e}} |\sigma(t)| dt \le (\exp(\tau/2-\tau^2/4) -1) 
(\tau-\tau^2/2) +\int_{\exp(\tau/2-\tau^2/4)}^{2/\sqrt{e}} \log t dt.
$$
Combining the above upper bounds on the integrals in the appropriate ranges,
we deduce, after several straightforward calculations, that 
the integral on the right side of (5.8) is indeed 
$< |\delta_1| -10^{-6}$
and so Theorem 5.1 
follows (for $\tau \ge 29/100$).  

\smallskip

{\sl Henceforth we suppose that 
$\tau \le 29/100$. }

\subhead 5d.  The range $u_0 \le u \le (2/\sqrt{e})u_0$ \endsubhead 

\noindent   We suppose in this section that $u_0\le u\le 
(2/\sqrt{e})(u_0)$.  Note that $\gamma(u\sqrt{e}/u_0)=0$ in this range.  
Observe that 
$$
\align
{\hat I}_1(u-u_0(1-1/\sqrt{e})) &= \tau + 
\int_{u_0/\sqrt{e}}^{u-u_0(1-1/\sqrt{e})} \hchi_1(t)dt \le \tau +2 
\log (1+\sqrt{e} (u/u_0-1))\\
&\le \tau + 2\sqrt{e} \log (u/u_0),\\
\endalign
$$
where the last inequality follows because $\log (1+\sqrt{e}(e^x-1))
\le\sqrt{e}x$ for all $x\ge 0$ (in fact, (right side)$-$(left side) is an increasing function
of $x$).  Inserting this in the middle bound of Lemma 5.3, we see 
by Lemma 5.2 that 
$$
\sigma(u)\le \frac{{\hat I}_2(u)}{2} 
\le \lambda(1-\tau)+\frac{\tau^2}2 
+2\sqrt{e} (\tau-\lambda) \log \frac{u}{u_0}.  
$$

Now define 
$$
\nu = \frac{\lambda(1-\tau)+\tau^2/2}{2(1-\sqrt{e}(\tau-\lambda))}.
$$
 From our assumption that $\tau \le 29/100$, it is 
easy to show that $\nu \le \tau/2-\tau^2/4 < 0.15$.   
Since $\sigma(u) \ge -2\log (u/u_0)$ by Lemma 5.2, we 
see by our upper bound above for $\sigma(u)$ that 
that 
$$
\align
|\sigma(u)| &\le \max(2\log (u/u_0), \lambda(1-\tau) +\tau^2/2
+2\sqrt{e}(\tau-\lambda)\log(u/u_0))\\
&= 
\cases 
2\log (u/u_0) &\text{if  } u\ge u_0 e^{\nu}\\
2\log (u/u_0) + 2(1-\sqrt{e}(\tau-\lambda)) (\nu-\log (u/u_0)) &\text{if  } 
u\le u_0e^{\nu}.\\
\endcases\\
\endalign
$$

Using the above upper bounds, we deduce that
$$
\align
\frac{1}{u_0}\int_{u_0}^{2u_0/\sqrt{e}} |\sigma(u)|du &\le \int_1^{2/\sqrt{e}}
2\log t \ dt + 2(e^{\nu}-1-\nu) (1-\sqrt{e}(\tau-\lambda)) \\
&\le \int_1^{2/\sqrt{e}} 2\log t \ dt + \frac{4}{15}
\frac{(\lambda(1-\tau)+\tau^2/2)^2}{1-\sqrt{e}(\tau-\lambda)},\\
\endalign
$$
by the definition of $\nu$ and since $e^x-1-x\le 8x^2/15$ if $0\le x\le 0.15$.

We simplify this a little by observing that, since $\nu \le \tau/2-\tau^2/4$,
$$
\align
\frac{(\lambda(1-\tau)+\tau^2/2)^2}{1-\sqrt{e}(\tau-\lambda)} &= 
\lambda(1-\tau) \frac{\lambda(1-\tau)}{1-\sqrt{e}(\tau-\lambda)} + 
\frac{\tau^2}{2} \frac{\lambda(1-\tau)}{1-\sqrt{e}(\tau-\lambda)} 
+ \tau^2 \nu\\
&\le \lambda\tau (1-\tau)^2 + \frac{\tau^2}{2} \tau(1-\tau) + \tau^2 
\left(\frac{\tau}{2}-\frac{\tau^2}{4}\right).\\
\endalign
$$
Inserting this in the previous estimate, and using $\tau\le 29/100$ we 
obtain
$$
\frac{1}{u_0}\int_{u_0}^{2u_0/\sqrt{e}} |\sigma(u)|\ du 
\le \int_{1}^{2/\sqrt{e}}
2\log t \ dt + \frac{\tau^2}{16} + \frac{\lambda}{25}
= \int_1^{2/\sqrt{e}} |\rho_{-}(t\sqrt{e})| dt +
\frac{\tau^2}{16}+ \frac{\lambda}{25}.  \tag{5.9}
$$

\subhead 5e. The range $2u_0/\sqrt{e} \le u\le (1+1/\sqrt{e})u_0$ \endsubhead

\noindent From Lemmas 5.2, 5.3 and 5.4 we 
see that in this range 
$$
\sigma(u) \le \frac{{\hat I}_2(u)}{2} \le \gamma\Big(\frac{u\sqrt{e}}{u_0}
\Big) 
+\tau -\frac{\tau^2}{2},
$$
and 
$$
\sigma(u) \ge -2\log \frac{u}{u_0} + \frac{{\hat I}_2(u)}{2}
-\frac{{\hat I}_3(u)}{6} \ge
-2\log \frac{u}{u_0} + (1-\lambda)\gamma\Big(\frac{u\sqrt{e}}{u_0}\Big) 
-\frac{\tau^2}{2}.
$$
Hence 
$$
|\sigma(u)| \le \max\Big( 2\log \frac{u}{u_0} - (1-\lambda)\gamma\Big(
\frac{u\sqrt{e}}{u_0}\Big) 
+\frac{\tau^2}2,\,\,\,
\gamma\Big(\frac{u\sqrt{e}}{u_0}\Big)+\tau-\frac{\tau^2}2\Big). 
\tag{5.10}
$$

\proclaim{Lemma 5.8} The function $2\log (t) -(1-\lambda)\gamma(t\sqrt{e})+\tau^2/2$ is 
increasing in the range $2/\sqrt{e} \le t < (1+1/\sqrt{e})$, 
and it is $>\tau-\tau^2/2+\gamma(t\sqrt{e})$.
\endproclaim

\demo{Proof}  For $2/\sqrt{e}\le t < (1+1/\sqrt{e})$ we have
$$
\frac{d}{dt} \left( 2\log (t) -(1-\lambda)\gamma(t\sqrt{e}) +\tau^2/2\right) =
\frac1t \left(2-4(1-\lambda) \log (t\sqrt{e}-1)\right) > 0,
$$
which gives the first statement.  Now,
$$
\frac{d}{dt} \left( \log (t) -\gamma(t\sqrt{e})\right) = 
\frac{1}{t}\left(1 - 4
\log (t\sqrt{e}-1)\right),
$$
which is positive in $(2/\sqrt{e}, e^{-1/2}+e^{-1/4})$ and negative in 
$(e^{-1/2}+e^{-1/4},(1+1/\sqrt{e}))$.  So the minimum of
$\log (t) -\gamma(t\sqrt{e})$  is attained at one of the end points 
$2/\sqrt{e}$ or $(1+1/\sqrt{e})$. The values taken by 
$\log (t) -\gamma(t\sqrt{e})$ at these two points are
$0.19\ldots$ and $0.1829\ldots$, respectively, which are both larger than
$1/8 \geq (\tau-\tau^2)/2$. Thus $\log t-\gamma(t\sqrt{e}) 
\ge (\tau-\tau^2)/2$ throughout our range. Doubling this and adding 
$\tau^2/2+\gamma(t\sqrt{e})$
to both sides implies the second assertion of the lemma, since 
$\lambda \gamma(t\sqrt{e})\geq 0$
\enddemo

By (5.10) and Lemma 5.8 we see that for any $2u_0/\sqrt{e} \le u\le
(1+1/\sqrt{e})u_0$ we have
$$
\frac{1}{u_0} \int_{2u_0/\sqrt{e}}^{u} |\sigma(t)|\ dt 
\le \int_{2/\sqrt{e}}^{u/u_0} 
\Big(2\log t- (1-\lambda)\gamma(t\sqrt{e})+\frac {\tau^2}2\Big) dt
$$
and that
$$
\align
\Big(1+\frac{1}{\sqrt{e}}-\frac{u}{u_0} \Big) |\sigma(u)| 
&\le 
\Big(1+\frac{1}{\sqrt{e}}-\frac{u}{u_0} \Big) \Big(2\log \frac{u}{u_0}
-(1-\lambda)\gamma(u\sqrt{e}) +\frac{\tau^2}2\Big)\\
&\le \int_{u/u_0}^{1+1/\sqrt{e}} \Big(2\log t- (1-\lambda) 
\gamma(t\sqrt{e})+\frac{\tau^2}2\Big)dt.\\
\endalign
$$
Adding these two inequalities, and noting that 
$\int_{2/\sqrt{e}}^{1+1/\sqrt{e}} \gamma(t\sqrt{e})dt = 0.0416\ldots <1/24$
we arrive at
$$
\align
\frac{1}{u_0} \int_{2u_0/\sqrt{e}}^{u} |\sigma(t)|dt &+
\biggl(1+\frac{1}{\sqrt{e}}-\frac{u}{u_0} \biggr) |\sigma(u)|
\\
&\le 
\int_{2/\sqrt{e}}^{1+1/\sqrt{e}} 
(2\log t -\gamma(t\sqrt{e}))dt 
+\frac{\lambda}{24} + \frac{\tau^2}{2}\biggl(1-\frac{1}{\sqrt{e}}\biggr)\\
&= \int_{2/\sqrt{e}}^{1+1/\sqrt{e}} 
|\rho_{-}(t\sqrt{e})| dt 
+\frac{\lambda}{24} + \frac{\tau^2}{2}\biggl(1-\frac{1}{\sqrt{e}}\biggr) .
\tag{5.11}\\
\endalign
$$

\subhead 5f. Completion of the proof of Theorem 5.1 \endsubhead

\noindent Recall that $\alpha=\exp(|\delta_1|/2-10^{-6})$ and that 
by Lemmas 5.2 and 5.3 we have 
$|\sigma(u)|\le \max(1/2,2\log (u/u_0))$ when $u_0\le u\le 
(1+\sqrt{e})u_0$. Moreover $\tau\leq 29/100$.
 Thus Theorem 5.1 holds in the range $u\le \alpha u_0$; and so,
below, we suppose that $u>\alpha u_0 > 2u_0/\sqrt{e}$.  
Put $v=\min(u_0(1+1/\sqrt{e}),u)$ so that by Lemma 3.5
$$
\align
\biggl(1+\frac{1}{\sqrt{e}}\biggr)|\sigma(u)| &= \frac{v}{u_0} |\sigma(u)| +
\biggl(1+\frac{1}{\sqrt{e}}-\frac{v}{u_0}\biggr)|\sigma(v)|\\ 
&\le \frac{1}{u_0} \int_0^v |\sigma(t)|dt 
+ \biggl(1+\frac{1}{\sqrt{e}}-\frac{v}{u_0}\biggr)|\sigma(v)|.\\
\endalign
$$  
Using the second part of Corollary 5.7 together with (5.9) and (5.11), 
and recalling that $2-2/\sqrt{e} = \int_0^1 |\rho_{-}(t\sqrt{e})|dt$, we 
see that this is 
$$
\le \int_0^{1+1/\sqrt{e}} |\rho_{-}(t\sqrt{e})| dt 
- \tau^2 \biggl(\frac{1}{\sqrt{e}}-\frac{9}{16}\biggr) -\frac{\lambda}{600}
\le \biggl(1+\frac{1}{\sqrt{e}}\biggr)|\delta_1| -\frac{\tau^2}{25}
-\frac{\lambda}{600},
$$
because of the identity (5.2).
It follows from this that $|\sigma(u)|\le |\delta_1|$ for all $u\ge u_0$.
Further, $|\sigma(u) - \delta_1| \le \epsilon$ implies that  
$\lambda\le \tau \ll \sqrt{\epsilon}$.  Since $I_1(u_0;\chi)=1$, we must have 
$$
\int_0^{u_0/\sqrt{e}} \frac{1-\chi(t)}{t} dt + \int_{u_0/\sqrt{e}}^{u_0} 
\frac{1+\chi(t)}{t} dt =2\tau \ll \sqrt{\epsilon}.  \tag{5.12}
$$

We now try to pinpoint further the case when 
$|\sigma(u)-\delta_1|\le \epsilon$.  Put $\chi_-(t) = 1$ for $t\le 
u_0/\sqrt{e}$, and $\chi_-(t)=-1$ for $t>u_0/\sqrt{e}$.  Note that 
the corresponding solution to (1.5) is $\rho_-(t\sqrt{e}/u_0)$.  Using 
Lemma 3.4 it follows that 
$$
\sigma(u) = \rho_- \Big(\frac{u\sqrt{e}}{u_0}\Big) + \sum_{k=1}^{\infty} 
\frac{(-1)^k}{k!} D_k,
$$
say, where 
$$
D_k = \int\Sb t_1,\ldots,t_k \ge 1\\ 
t_1+\ldots+t_k \le u\endSb \frac{\chi_-(t_1)-\chi(t_1)}{t_1} 
\ldots \frac{\chi_-(t_k)-\chi(t_k)}{t_k} \rho_- \Big( \frac{(u-t_1-\ldots-t_k)
\sqrt{e}}{u_0}\Big)dt_1\ldots dt_k.
$$

Suppose first that $u_0 < u\le 3u_0/\sqrt{e}$.  Notice that when $k\ge 2$,  
at least one of the $t_i$ must be $\le u_0$.  It follows that for $k\ge 2$,
$$
|D_k| \le k \Big(\int_0^{u_0} \frac{|\chi_-(t)-\chi(t)|}{t} dt\Big) 
\Big(\int_0^{3u_0/\sqrt{e}} \frac{|\chi_-(t)-\chi(t)|}{t}dt\Big)^{k-1}.
$$
By (5.12), we see that the first factor above is $\le C\sqrt{\epsilon}$ 
for some constant $C$; and clearly the second factor is 
$\le (C\sqrt{\epsilon} + 2\int_{u_0}^{3u_0/\sqrt{e}} \frac{dt}{t})^{k-1} 
\le 2^{k-1}$, if $\epsilon$ is small enough.  So, for $k\ge 2$, 
$|D_k|\le 2^{k-1}C k \sqrt{\epsilon}$ whence it follows by (5.12) that 
$$
\align
\sigma(u) &= \rho_-\Big(\frac{u\sqrt{e}}{\epsilon}\Big) - \int_0^{u} 
\frac{\chi_-(t)-\chi(t)}{t}  \rho_-\Big(\frac{(u-t)\sqrt{e}}{u_0}\Big) dt 
+ O(\sqrt{\epsilon})\\
&= \rho_-\Big(\frac{u\sqrt{e}}{u_0}\Big) + \int_{u_0}^{u} \frac{1+\chi(t)}{t} 
\rho_-\Big(\frac{(u-t)\sqrt{e}}{u_0}\Big) dt + O(\sqrt{\epsilon}).\\
\endalign
$$
Notice that $\rho_{-}$ is a non-increasing function in the range 
$[0,3-\sqrt{e}]$.  Thus the 
$\rho_{-}((u-t)\sqrt{e}/u_0)$ term in the right side of the equation above lies 
between $\rho_{-}(0)=1$ and $\rho_{-}(3-\sqrt{e}) = 0.3978\ldots$.  
Hence when $u_0 \le u\le 3u_0/\sqrt{e}$, we conclude that 
$$
\frac{1}{3} \int_{u_0}^{u} \frac{1+\chi(t)}{t} dt +O(\sqrt{\epsilon}) 
\le \sigma(u) -\rho_{-}
\Big(\frac{u\sqrt{e}}{u_0}\Big)
\le \int_{u_0}^{u} \frac{1+\chi(t)}{t} dt +O(\sqrt{\epsilon}). \tag{5.13}
$$
In the range $t\in [0,3]$ we know that $\rho_-$ has its minumum of 
$\delta_1$ at $t= (1+\sqrt{e})$, and further it is easy to check 
that $(\rho_-(t)-\delta_1) \gg (t-(1+\sqrt{e}))^2$.  From this 
and (5.13) it follows that if $|\sigma(u)-\delta_1| \le \epsilon$ 
for $u$ in the range $[u_0,3u_0/\sqrt{e}]$ then we must have 
$u/u_0 =(1+1/\sqrt{e}) + O(\epsilon^{\frac 14})$, and 
that $\int_{u_0}^u \frac{1+\chi(t)}{t} dt = O(\sqrt{\epsilon})$.  This 
proves Theorem 5.1 in this range of $u$.

Now suppose that $u>3u_0/\sqrt{e}$.  By Lemma 3.5, we note 
that 
$$
|\sigma(u)| \le \frac{1}{3u_0/\sqrt{e}} \int_0^{3u_0/\sqrt{e}} |\sigma(t)| dt.
$$
Now by (5.13) and a simple computation, we see that for $t\le 3u_0/\sqrt{e}$, 
$$ 
\rho_{-}\Big(\frac{t\sqrt{e}}{u_0}\Big) +O(\sqrt{\epsilon}) \le \sigma(t) 
\le \rho_-\Big(\frac{t\sqrt{e}}{u_0}\Big) + 2\log \frac{u}{u_0} +O(\sqrt{\epsilon})
\le -\rho_{-}\Big(\frac{t\sqrt{e}}{u_0}\Big)+O(\sqrt{\epsilon})
$$
whence $|\sigma(t)| \le |\rho_{-}(t\sqrt{e}/u_0)| +O(\sqrt{\epsilon})$.  
Inserting this into our bound for $|\sigma(u)|$, we get 
$$
|\sigma(u)| \le 
\frac{\sqrt{e}}{3u_0} \int_0^{3u_0/\sqrt{e}} |\rho_{-}(t\sqrt{e}/u_0)|dt 
+O(\sqrt{\epsilon})
\le 0.61 +O(\sqrt{\epsilon}).
$$
Thus if $\epsilon$ is small enough, then 
$|\sigma(u)-\delta_1|\le \epsilon$ is impossible for $u>3u_0/\sqrt{e}$.  
This completes our proof of Theorem 5.1.

\head 6. The Euler Product Spectrum \endhead

\subhead 6a.  Proof of Theorem 4 \endsubhead

\noindent  Suppose $z\in \GT$ so that there exists $f\in {\Cal F}(S)$ 
with $\Theta(f,\infty) = z$.  Suppose $\alpha \in S$ and define $g=g_y \in 
{\Cal F}(S)$ by $g(p)=f(p)$ for $p\le y$, $g(p)= \alpha$ for $y< p \le y^{e^k}$
and $g(p) = 1$ for $p>y^{e^k}$; here $k\ge 0$ is some real number.  Then 
by the prime number theorem $\Theta(g,\infty) =\Theta(f,y) e^{-k(1-\alpha)} 
+ O(1/\log y)$.  Letting $y\to \infty$ we have shown that $e^{-k(1-\alpha)} z
\in \GT$ for all $\alpha \in S$ and all $k\ge 0$.  

Now suppose $\alpha =\sum_{j=1}^{l} k_j \alpha_j$ belongs to the 
convex hull of $S$; where $\alpha_j \in S$ and $k_j\ge 0 $ with 
$\sum_{j=1}^l k_j =1$.  If $z\in \GT$ we see,
 by using the result of the preceding 
paragraph $l$ times, that for all $k\ge 0$,  
$$
e^{-k(k_1(1-\alpha_1)+\ldots +k_l(1-\alpha_l))} z = e^{-k(1-\alpha)} z \in 
\GT. 
$$
This shows that $\GT \supset {\Cal E}(S) \times \GT$.  Since $1\in \GT$ 
and $1\in {\Cal E}(S)$ we have $\GT = {\Cal E}(S)\times \GT
\supset {\Cal E}(S)$.

To demonstrate that $\GT \subset {\Cal E}(S) \times [0,1]$, we 
require the following technical lemma.

\proclaim{Lemma 6.1} Let $p>1$ and  
$\alpha=a+ib\in {\Bbb U}$ with $\alpha\ne \pm 1$.  Let
$$
s:= -\frac 1b \arg(p-\alpha) = 
\frac 1{|b|} \arctan \biggl( \frac{|b|}{p-a} \biggr),\
\text{and put }\
r:=\left(\frac{p-1}{p-\alpha}\right)
e^{(1-\alpha)s}.
$$ 
Then $r$ is a real number in the range
$0\le r \le 1$.  
\endproclaim

The upper bound for $r$ is tight and is attained in the situation 
$\alpha \to 1$.  We can also show that $r \ge e^{-2/(p+1)} ((p-1)/(p+1))$
which is attained when $\alpha \to -1$; but this is not necessary for our 
applications.

\demo{Proof} Since $-bs = \arg(p-\alpha)$ we see that $r$ is a non-negative 
real number.  Since $\arctan (t) \le t$ for all $t\in [0,\infty)$ ($\arctan$ 
is to lie between $0$ and $\pi/2$ here) we get
$$
|e^{(1-\alpha)s}| = e^{(1-a)s} = 
\exp\biggl(\frac{1-a}{|b|}\arctan\biggl(\frac{|b|}{p-a}\biggr)\biggr) 
\le \exp\biggl( \frac{1-a}{p-a}\biggr).
$$
Next, as $e^{-t} \ge 1-t$, we get $|(p-1)/(p-\alpha)| \le (p-1)/(p-a) \le 
\exp(-(1-a)/(p-a))$.  Hence
$$
r= |r| = \frac{p-1}{|p-\alpha|} |e^{(1-\alpha)s}| \le
 \exp\biggl(-\frac{1-a}{p-a} + \frac{1-a}{p-a}\biggr) =1,
$$
as desired.

\enddemo

For $f\in {\Cal F}(S)$, let $k_p = 0$ if $f(p)$ is real, 
and $k_p = - \arg(p-f(p))/\IM f(p) \ge 0$ otherwise.  By 
Lemma 6.1, we 
may conclude that
$$
\frac{p-1}{p-f(p)} = r_p e^{-k_p(1-f(p))}
$$
for a real number $0\le r_p \le 1$.  Taking the product over all primes we 
get
$$
\Theta(f,\infty) = \biggl(\prod_p r_p\biggr) \exp\biggl(-\sum_p k_p (1-f(p))
\biggr).
$$
This shows that $\GT \subset {\Cal E}(S) \times [0,1]$.

For every $1\neq \alpha \in {\Bbb U}$ we note that $(1-\RE \alpha)/|\IM \alpha|
= \cot (|\arg(1-\alpha)|)$.  Hence for all $1\neq \alpha$ in the 
convex hull of $S$ we have $(1-\RE \alpha)/|\IM \alpha| \le 
\cot(\text{Ang }(S))$.  This shows that if $z\in {\Cal E}(S)$ then 
$|z| \le \exp(-|\arg z| \cot (\text{Ang }(S)))$. Since 
$\GT \subset {\Cal  E}(S) \times [0,1]$ the same upper bound holds for all 
$z\in \GT$.  

Now suppose that $1\neq \beta $ is a real number in the convex hull of $S$.  
Then 
$$
{\Cal E}(S) \supset {\Cal E}(S) \times \{e^{-k(1-\beta)}: \ \ k\ge 0\} 
= {\Cal E}(S)\times [0,1].
$$
Hence ${\Cal E}(S) = {\Cal E}(S) \times [0,1]$ in this case, which completes
the proof of Theorem 4.

\subhead 6b.  Proof of Corollary 3 \endsubhead
\noindent By Theorems 3 and 4, 
$$
\Gamma(S) = \GT \times \Lambda(S) = {\Cal E}(S) \times \GT \times \Lambda(S)
= {\Cal E}(S) \times \Gamma(S).
$$
If $S=\{1\}$ then $\Gamma(S)=\{1\}$ and Ang$(S)=0$ so that (ii) is 
immediate in this case.  Suppose then that $1\neq \alpha \in S$.  
If $z\in \Gamma(S)$ then we know that $e^{-k(1-\alpha)}z \in \Gamma(S)
$ for all $k\ge 0$.  Letting $k$ vary from $0$ to $\infty$ we get a 
spiral connecting $z$ to $0$.  This shows that $\Gamma(S)$ is connected.  
If the convex hull of $S$ contains a real point other that $1$ then we know 
from Theorem 4 that $\GT = [0,1] \times \GT$ is starlike.  Hence
$\Gamma(S) = \GT \times \Lambda(S) = [0,1] \times \GT \times \Lambda(S)
= [0,1] \times \Gamma(S)$ is starlike as well.  This completes 
the proof of part (i).  

Suppose $\chi(t) =1$ for $t\le 1$ and $\chi(t) =\alpha$ for $t>1$, and
let $\sigma$ denote the corresponding solution to (1.5).  Then by 
Theorem 3.3, we see that 
$$
\sigma (u) = 1- \int_1^u \frac{1-\chi(t)}{t} dt = 1 -(1-\alpha) \log u,\ \
\text{for  } 1\le u\le 2.
$$
If $\alpha$ is in the convex hull of $S$, this shows that $1-(1-\alpha)\log u
\in \Lambda(S)$ for $1\le u \le 2$.
Suppose now that $\pi/2 > \text{Ang }(S) >0$ and that $\pm 1 \neq \zeta \in S$
with Ang$(\zeta) =\text{Ang }(S) = \theta$, say.  
Suppose $1< u \le 2$ and let $z=1- (1-\zeta)\log u$  
so that $z\in \Lambda(S)$.  If $|\arg z| =\nu$ 
then a simple geometric consideration shows that $|z| = \sin \theta/ 
\sin(\theta +\nu)$.  On the other hand, if $z\in \GT$ then by Theorem 4
$$
|z| \le \exp(-\nu \cot \theta) \le \frac{1}{1+\nu \cot\theta} 
< \frac{1}{\cos \nu + \sin \nu \cot \theta} 
= \frac{\sin \theta}{\sin(\theta+\nu)},
$$
which is a contradiction. This proves part (ii).

If $1$, $e^{i\alpha}$ and $e^{i\beta}$ are distinct elements of $S$ then
for all $k$, $l \ge 0$ we know that $e^{-k(1-e^{i\alpha}) -l (1-e^{i\beta})} 
\in {\Cal E}(S)$.  Now let us fix
the real part of $k(1-e^{i\alpha}) + l(1-e^{i\beta})$; that is let us 
fix $2k \sin^2(\alpha/2) +2l \sin^2(\beta/2) = r$, say.  Then as $k$ varies 
from $0$ to $r/(2\sin^2(\alpha/2))$ we see that the imaginary part of
$k(1-e^{i\alpha})+l(1-e^{i\beta})$, which is $-k\sin \alpha -l \sin \beta$,
varies  continuously from $-r\cot(\beta/2)$ to $-r\cot(\alpha/2)$.  If the 
variation in the imaginary  part is larger than $2\pi$ in magnitude 
then ${\Cal E}(S)$ clearly contains the circle with center $0$ and radius 
$e^{-r}$; this happens provided $r \ge 2\pi/|\cot(\alpha/2) -\cot(\beta/2)|$.
Hence we have proved (iii).

\subhead 6c. Proof of Theorem 3$'$ \endsubhead

\noindent Suppose $z_\sigma \in \Lam(S)$.  We shall show that 
for any $k\ge 0$ and any $\alpha \in S$, $e^{-k(1-\alpha)} z_\sigma 
\in \Lambda(S)$.  Using this repeatedly (as in \S 6a) it would 
follow that $z_\sigma e^{-k(1-\alpha)} 
\in \Lam(S)$ for any $k\ge 0$ and $\alpha$ in the convex hull of $S$.  
This means that $\Lam(S) \supset \Lam(S) \times {\Cal E}(S)$, 
and since $1\in {\Cal E}(S)$ it would follow that $\Lam(S) = \Lam(S) 
\times {\Cal E}(S)$.  

Since $z_\sigma \in \Lam(S)$, we know that there is a measurable function 
$\chi \in K(S)$, and $u\ge 1$ such that $z_\sigma= \sigma(u)$, 
where $\sigma$ is the corresponding solution to (1.5).  By Proposition 1 
(Converse), for large $x$ we may find $f\in {\Cal F}(S)$ with $f(p)=1$ 
for $p\le x^{\frac 1u}$ such that 
$$ 
\frac{1}{x} \sum_{n\le x} f(n) = \sigma(u) +o(1) = z_{\sigma} +o(1).
$$
Now put $y=\exp((\log x)^{\frac 23})$, and define $g\in {\Cal F}(S)$ by 
$g(p)=1$ for $p\le y$ or $p>y^{e^k}$  and $g(p)=\alpha$ for $y< p\le y^{e^k}$. 
Hence $\Theta(g,x^{\frac 1u}) = e^{-k(1-\alpha)} +o(1)$.  
Consider $h\in {\Cal F}(S)$ defined by 
$h(p)= g(p)$ if $p\le x^{\frac 1u}$, and $h(p)= f(p)$ for $p>x^{\frac 1u}$.  
By Proposition 1, we see easily that $\frac 1x \sum_{n\le x} h(n) +o(1)$ 
belongs to $\Lam(S)$.  On the other hand, appealing to Proposition 4.5 we
obtain
$$
\align
\frac{1}{x} \sum_{n\le x} h(n) &= \Theta(g,x^{\frac{1}{u}}) \frac{1}{x} 
\sum_{n\le x} f(n) + O\Big( \frac{\log y}{\log x} \exp\Big( 
\sum_{y< p\le y^{e^k}} \frac{2}{p} + \sum_{x^{\frac 1u} < p\le x} \frac{2}{p}
\Big) \Big) \\
&= e^{-k(1-\alpha)} z_\sigma + o(1).\\
\endalign
$$
Letting $x\to \infty$, we conclude that $e^{-k(1-\alpha)} z_\sigma \in 
\Lam(S)$, as desired.

Plainly $\Lambda(S) \subset \Gamma(S)$.  Further, by the 
result just established and Theorems 3 and 4,  
$$
\Lam(S) \times [0,1]= \Lam(S) \times {\Cal E}(S) \times [0,1] 
\supset \Lam(S) \times \GT = \Gamma(S).
$$
If the convex hull of $S$ contains a 
real point other than $1$ then by Theorem 4, ${\Cal E}(S) ={\Cal E}(S) \times 
[0,1] = \GT$ and hence $\Lambda(S) = \Lambda(S) \times [0,1] =\Gamma(S)$.

\head 7. Angles and Projections of the Spectrum \endhead

\subhead 7a. Proof that Ang$(\Gamma(S))\ll \text{Ang}(S)$ \endsubhead

\noindent Let Ang$(S)=\theta$ and we seek to show that Ang$(\Lambda(S))=
\text{Ang}(\Gamma(S)) \ll \theta$.  Suppose $\chi \in K(S)$ is given,
and let $\sigma$ denote the corresponding solution to (1.5).  We 
need to show that Ang$(\sigma(u)) \ll \theta$ for all $u$; or, in other words,
$|\IM \sigma(u)|\ll \theta(1-\RE\sigma(u))$.

We may suppose that $\theta$ is sufficiently 
small, else the result is trivial.  We let 
$\hchi =\RE \chi$ and let $\hsigma$ denote the corresponding solution to (1.5).
Recall that, 
in the notation of \S 3b, 
$R_1(u)=R_1(u;\chi)= \int_1^u \frac{1-\RE \chi(t)}{t}dt 
=\int_1^u \frac{1-\hchi(t)}{t}dt$.  By applying Lemma 1$'$, taking $D$ there 
to be the region $\{ z\in {\Bbb U}: \text{Ang}(z) \le \pi/4\}$,  
we see that for all $v$,
$$
|\sigma(v)|, |\hsigma(v)| \le c_1 \exp(-c_2 R_1(v)) \tag{7.1}
$$
where $c_1$ and $c_2$ are absolute positive constants.

By simple trigonometry, for any $z\in {\Bbb U}$, Ang$(z)\le \arcsin |z|$.  Hence
we may assume that $|\sigma(u)| \gg \sin \theta \gg \theta$, whence $R_1(u)
\ll \log (1/\theta)$, by (7.1).  

By Lemma 3.4
$$
\sigma(u)=\hsigma(u)- \int_1^u \frac{\hchi(t)-\chi(t)}{t} \hsigma(u-t) dt
+O\biggl(\sum_{j=2}^{\infty} \frac{1}{j!} \biggl(\int_1^u \frac{|\hchi(t)-
\chi(t)|}{t} dt \biggr)^j\biggr).
$$
Since $\chi$ is in the convex hull of $S$ and Ang$(S)=\theta$, $|\IM\chi(t)|
=|\hchi(t)-\chi(t)| \le \tan\theta (1-\RE \chi(t))$.  Hence $\int_{1}^{u}
\frac{|\hchi(t)-\chi(t)|}{t} dt \le \tan\theta R_1(u)$.  Using this above, 
and as $R_1(u)\ll \log (1/\theta)$, we get 
$$
\sigma(u)=\hsigma(u) + i\int_1^u \frac{\IM\chi(t)}{t} \hsigma(u-t) dt 
+ O(\theta^2R_1(u)^2).\tag{7.2}
$$

If $R_1(u)\ge 1$ then by Theorem 5.1, $|\hsigma(u)|\le |\delta_1|$.  If 
$R_1(u)\le 1$ then by the inclusion-exclusion 
inequalities of Proposition 3.6, 
$\hsigma(u)\le \hsigma_2(u) \le 1 - R_1(u) +R_1(u)^2/2
\le 1-R_1(u)/2$.  Hence, using (7.2) and as $R_1(u)\ll \log (1/\theta)$,
$$
1-\RE \sigma(u) = 1-\hsigma(u)+O(\theta^2 R_1(u)^2) 
\geq  R_1(u)+O(\theta^2 R_1(u)^2) \gg \min(R_1(u),1),
$$
since $\theta$ is sufficiently small.

Taking imaginary parts in (7.2), and recalling $|\IM \chi(t)|
\ll \theta(1-\RE \chi(t))$ and (7.1), we see that
$$
\align
|\IM \sigma(u)|&\le \int_1^{u} \frac{|\IM \chi(t)}{t} |\hsigma(u-t)| dt 
+O(\theta^2 R_1(u)^2)\\
&\ll \theta\int_{u/2}^{u} \frac{1-\RE \chi(t)}{t} dt +\theta \int_{1}^{u/2} 
\frac{1-\RE\chi(t)}{t} e^{-c_2R_1(u-t)} dt + \theta^2 R_1(u)\\
&\ll \theta (R_1(u)-R_1(u/2)) + \theta R_1(u/2)\exp(-c_2 R_1(u/2)) 
+\theta^2 R_1(u).
\\
\endalign
$$
Since $R_1(u)-R_1(u/2)\le 2\log 2\ll 1$, and $R_1(u)\ll \log (1/\theta)$,
the above shows that $|\IM \sigma(u)| \ll \theta
\min(R_1(u),1)$.  Combining this with our lower bound for $1-\RE \sigma(u)$, 
gives Ang$(\sigma(u))\ll \tan($Ang$(\sigma(u)))= |\IM \sigma(u)|/(1-\RE \sigma(u)) \ll \theta$, completing the proof.

\subhead 7b.  The maximal projection of $S=\{\pm 1, \pm i \}$ \endsubhead

\noindent In this section we shall prove Theorem 7(i).  
The result for $S=\{1,-1\}$ follows from Theorem 5.1, and so we 
may restrict ourselves to the case $S=\{\pm 1,\pm i\}$.  By Theorem 3$'$ 
we see that $\Gamma(S)=\Lambda(S)$ so we shall work here with $\Lambda(S)$.
Let $\chi \in K(\{\pm 1, \pm i\})$ be given, and let $\sigma$ be the 
corresponding solution to (1.5). 
We shall show that for all $u$, $\RE \sigma(u)\ge 
-(1+|\delta_1|)/2$ and $|\IM \sigma(u)|\le (1+|\delta_1|)/2$, so that 
the maximal projection of $\{ \pm 1, \pm i \}$ is $(1+|\delta_1|)/2$ 
as conjectured; that is, Theorem 7(i).

\proclaim{Lemma 7.1} Let $\chi^{\prime}$ be any real-valued measurable function
satisfying 
$$
|\RE \chi(t)| \le \frac{1+\chi^{\prime}(t)}{2} \qquad \text{and} 
\qquad |\IM \chi(t)|\le \frac{1-\chi^{\prime}(t)}{2}.
$$
for all $t$.  Let $\sigma^{\prime}$ be the corresponding solution 
to (1.5).  Then, for all $u$,
$$
|\RE \sigma(u)| \le \frac{1+\sigma^{\prime}(u)}{2} \qquad \text{and} 
\qquad |\IM\sigma(u)|\le \frac{1-\sigma^{\prime}(u)}{2}.
$$
\endproclaim

\demo{Proof} Let $\beta(u):= (1+\sigma^{\prime}(u))/2-|\RE\sigma(u)|$ and 
$\gamma(u):=(1-\sigma^{\prime}(u))/2-|\IM \sigma(u)|$. Since
$$
\align
u|\RE\sigma(u)| &\leq |\RE \chi| * |\RE\sigma| + |\IM \chi| * |\IM\sigma|\\
&\leq 
\frac{1+\chi^{\prime}}{2} * |\RE\sigma| + \frac{1-\chi^{\prime}}{2}  * |\IM\sigma| ,\\
\endalign
$$
we deduce that $u\beta(u) \geq (1+\chi^{\prime})/2 * \beta + (1-\chi^{\prime})/2 *\gamma$. Similarly, by bounding $|\IM\sigma(u)|$ we get 
$u\gamma(u) \geq (1-\chi^{\prime})/2 * \beta + (1+\chi^{\prime})/2 *\gamma$.
Taking $\alpha(u) = \min \{ \beta(u), \gamma(u)\}$ we have $\alpha(u)=0$
for $0\leq u\leq 1$, and we deduce from the above that 
$u\alpha(u) \geq 1*\alpha$. Therefore $\alpha(u)\geq 0$ for all $u$,
by Lemma 3.1.
\enddemo

\demo{Proof of Theorem 7(i)}  
We wish to show that $|\IM \sigma(u)|$ and $-\RE \sigma(u)$ are both
$\leq  (1-\delta_1)/2$. Note that $\chi'$ exists, as in Lemma 7.1 since
the convex hull of $S$ is described by the conditions
$|\RE \chi(t)|+|\IM \chi(t)|\leq 1$.
By Theorem 5.1, we know that $\sigma^{\prime}(u)\ge \delta_1$ always.  
Hence by Lemma 7.1, $|\IM\sigma(u)| \le (1-\delta_1)/2 = (1+|\delta_1|)/2$. 
Further, if $I_1(u;\chi^{\prime}) \ge 1$ then $|\sigma^{\prime}(u)|\le 
|\delta_1|$ by Theorem 5.1, so that $|\RE \sigma(u)|\le (1+|\delta_1|)/2$.

We now handle the case when $I_1(u,\chi^{\prime})\le 1$.  Put $\hchi=\RE \chi$
and let $\hsigma$ be the corresponding real-valued solution to (1.5). 
By Proposition 3.7 and Theorem 5.1, 
$$
\RE\sigma(u)\ge \hsigma(u)-\frac{C_2(u)}{2}
\ge \delta_1 -\frac{C_2(u)}{2} 
\ge \delta_1 -\frac{C_1(u)^2}{2}.
$$
Now 
$$
C_1(u) =\int_1^u \frac{|\IM\chi(t)|}t dt 
\le \int_1^u \frac{(1-\chi^{\prime}(u))}{2t} dt \le \frac 12,
$$ 
and so $\RE \sigma(u)\ge \delta_1 -1/8 > - (1+|\delta_1|)/2$, which 
completes our proof.

\enddemo

\subhead 7c.  Towards the proofs of Theorems 5, 6(ii), and 7(ii) \endsubhead

\noindent  In the following subsections, we suppose that $S$ is a 
given subset of ${\Bbb U}$ with Ang$(S)=\theta<\pi/2$.  
Suppose that $\chi \in K(S)$ is given, and that 
$\sigma(u)$ is the corresponding solution 
to (1.5).  Define 
$$
P(u) = \int_{0}^{u} \min(2,(1-\RE \chi(t))\sec^2\theta )\frac{dt}{t}.
$$
Let $u_0$ be such that $P(u_0)+P(u_0/2)=1$; if no such point exists, 
set $u_0=\infty$.

\proclaim{Lemma 7.2} With these notations $P(u)\cos^2\theta \le 
R_1(u) \le P(u)$, where $R_i, C_i$ are as in section 3b.  Further 
$$
R_1(u)^2 + C_1(u)^2 \le R_1(u)P(u), \tag{7.3}
$$
$$
R_2(u)+C_2(u) \le \min(R_1(u)P(u), 2R_1(u)\sqrt{P(u/2)P(u)}), \tag{7.4}
$$
$$
|\IM\sigma(u)|\le \sqrt{R_1(u)(P(u)-R_1(u))}, \tag{7.5}
$$
and
$$
R_1(u)(1-P(u)/2) \le 1-\RE\sigma(u)\le R_1(u)(1+P(u)/2).
\tag{7.6}
$$
\endproclaim

\demo{Proof}  It is clear from the definitions that $P(u)\cos^2\theta \le
R_1(u) \le P(u)$. Since $\chi(t)$ lies in the convex 
hull of $S$, and Ang$(S)=\theta$, 
we have 
$$
|\IM \chi(t)| \le \min(\sqrt{1-(\RE\chi(t))^2},(1-\RE\chi(t))\tan\theta ).
$$
Using Cauchy's inequality we obtain
$$
C_1(u)^2 \le R_1(u) \int_1^u \min(1+\RE \chi(t), (1-\RE\chi(t))\tan^2\theta )
\frac{dt}{t}.
$$
Adding $R_1(u)^2$ to the above, we obtain (7.3).  By Proposition 3.7,
$|\IM\sigma(u)|\le C_1(u)$, and so we deduce (7.5).

Plainly $R_2(u)\le R_1(u)^2$, and $C_2(u)\le C_1(u)^2$.  So the 
first bound in (7.4) follows from (7.3).  Further, from the definition of 
$R_2$, we have $R_2(u)\le 2R_1(u/2) R_1(u)- R_1(u/2)^2 \le 2R_1(u/2)R_1(u)$,
and similarly $C_2(u)\le 2C_1(u/2)C_1(u)$.  By Cauchy's inequality, and (7.3),
$$
\align
R_2(u) + C_2(u) &\le 2 (R_1(u/2)^2+ C_1(u/2)^2)^{1/2} (R_1(u)^2+C_1(u)^2)^{1/2}
\\
&\le 2\sqrt{R_1(u/2)P(u/2)R_1(u)P(u)},\\
\endalign
$$
and the second bound of (7.4) follows as $R_1(u/2)\le R_1(u)$.

By Proposition 3.7 we know that 
$$
R_1(u) -\frac{R_2(u)+C_2(u)}{2} \le 1-\RE\sigma(u) \le R_1(u) + 
\frac{R_2(u)+C_2(u)}{2},
$$
and using the first bound of (7.4), we obtain (7.6).
\enddemo

We next prove a technical Lemma which will be useful in the proof of 
Lemma 7.4.

\proclaim{Lemma 7.3}  If $a,b\geq c>0$ and $0\le x,y\leq 1$ then
$$ 2ax+2by - (\sqrt{a}x+\sqrt{b}y)^2 \geq c (x+y)(2-x-y).$$
\endproclaim

\demo{Proof} Without loss of generality assume $a\geq b$. We shall
prove that result for $c=b$, and then the more general statement follows.
First note that $(\sqrt{a}+\sqrt{b})(2-x)\geq 2\sqrt{b}
\geq 2\sqrt{b}y$. Multiplying this through by $(\sqrt{a}-\sqrt{b})x$
and adding $b(2y-y^2+2x-x^2)$ to both sides, we get
$a(2x-x^2)+b(2y-y^2)\geq 2xy\sqrt{ab}+b(2(x+y)-(x+y)^2)$
after some re-arranging. This directly implies the result.
\enddemo

\proclaim{Lemma 7.4} Suppose that $u\ge u_0$.  Then 
$$
|\sigma(u)|^2 \le 1 - \frac{\cos^2\theta}{u_0} \int_0^{u_0/2} 
(P(t)+ P(u_0-t))(2-P(t)-P(u_0-t))dt.
$$
\endproclaim
\demo{Proof}  If $t\le u_0$ then $R_1(t)\le P(t)\le 1$ by Lemma 7.2, and so by Proposition 3.7,
$$
\align
|\RE \sigma(t)|&\le \max\biggl(\!1-R_1(t)+\frac{R_2(t)+C_2(t)}{2},
-1+R_1(t)+\frac{C_2(t)}{2}\!\biggr) \\
&= 1- R_1(t) + \frac{R_2(t)+C_2(t)}{2}.\\
\endalign
$$
Using (7.5) we deduce
$$
\align
|\sigma(t)|^2 &\le 1 - 2R_1(t) + R_1(t) P(t) +R_2(t)+ C_2(t)\\
& \qquad + 
(R_2(t)+C_2(t))\biggl(\frac{R_2(t)+C_2(t)}{4}- R_1(t)\biggr).\\
\endalign
$$
By (7.4), $R_2(t)+C_2(t)\le R_1(t) P(t)\le R_1(t)$
and so for $t\le u_0$, we have shown
$$
|\sigma(t)|^2 \le 1- 2R_1(t) + R_1(t)P(t) + R_2(t) +C_2(t).\tag{7.7}
$$
By Lemma 3.5, Cauchy's inequality, and the above bound 
we obtain for $u\ge u_0$
$$
\align
|\sigma(u)|^2 &\le \biggl(\frac{1}{u_0} \int_0^{u_0} |\sigma(t)|dt \biggr)^2
\le \frac{1}{u_0} \int_0^{u_0} |\sigma(t)|^2 dt \\
&\le \frac{1}{u_0} \int_0^{u_0} (1-2R_1(t) + R_1(t)P(t) +R_2(t)+C_2(t))dt.
\tag{7.8}\\
\endalign
$$

Denote $\chi_1(t) = (1-\RE \chi(t))/t$, so that $R_1=1*\chi_1$ and $R_2
=1*\chi_1*\chi_1$.  It follows that $1*R_2= 1*\chi_1*1*\chi_1=R_1*R_1$.  
In like manner, $1*C_2=C_1*C_1$.  Using this, Cauchy's inequality, and (7.3),
we obtain
$$
\align
\int_0^{u_0} (R_2(t)+&C_2(t))dt = (R_1*R_1)(u_0)+(C_1*C_1)(u_0) \\
&\le (\sqrt{R_1^2+C_1^2}*\sqrt{R_1^2+C_1^2})(u_0)
 \le (\sqrt{R_1P}*\sqrt{R_1P})(u_0) \\
&= 2\int_0^{u_0/2} \sqrt{R_1(t)P(t)R_1(u_0-t)P(u_0-t)}dt.\\
\endalign
$$
Using this in (7.8) we deduce that $|\sigma(u)|^2 \le 1-J$ where 
$$
J= \frac{1}{u_0} \int_0^{u_0/2} \biggl(
2R_1(t)+2R_1(u_0-t) 
- \left(\sqrt{R_1(t)P(t)}+\sqrt{R_1(u_0-t)P(u_0-t)}\right)^2
\biggr)dt.
$$

For $0\le t\le u_0/2$, take $a=R_1(t)/P(t)$, $b=R_1(u_0-t)/P(u_0-t)$, 
so that $a$ and $b$ are $\ge \cos^2\theta$ by Lemma 7.2.  
Take $x=P(t)$ and $y=P(u_0-t)$,
so that both $x$ and $y$ are $\le 1$.  Applying Lemma 7.3, the
integrand in the definition of $J$ is $\ge \cos^2\theta (P(t)+P(u_0-t))
(2-P(t)-P(u_0-t))$; which proves the Lemma.

\enddemo

Using Lemma 7.4 we can get an explicit bound on $|\sigma(u)|$ when $u\ge u_0$.

\proclaim{Proposition 7.5}  If $u\ge u_0$ then $|\sigma(u)|\le 1- (56/411)
\cos^2 \theta$.
\endproclaim

\demo{Proof}  Put $\alpha=P(u_0/2)$ so that $0\le \alpha \le 1/2$, and
$P(u_0)=1-\alpha$.  For $0\le t\le u_0/2$, note that 
$$
P(t) \ge P(u_0/2) -\int_t^{u_0/2} 2\frac{dv}{v} = \alpha - 2\log (u_0/(2t)),
$$
and also $P(t)\ge 0$.  Similarly
$$
P(u_0-t) \ge P(u_0) -\int_{u_0-t}^{u_0} 2\frac{dv}{v} 
= 1-\alpha - 2\log(u_0/(u_0-t)),
$$
and also $P(u_0-t) \ge P(u_0/2) =\alpha$.  Thus if we put 
$$
m(t)= 
\cases 
1-\alpha +2\log (1-t/u_0) &\text{for } t/u_0\le 1- e^{\alpha}/\sqrt{e}\\
\alpha &\text{for } 1-e^{\alpha}/\sqrt{e} \le t/u_0 \le 1/(2e^{\alpha/2})\\
2\alpha +2\log (2t/u_0) &\text{for }  1/(2e^{\alpha/2} ) \le t/u_0 \le 1/2,\\
\endcases
$$
then $P(t)+P(u_0-t)\geq m(t) \geq 0$ for each $0 \le t\le u_0/2$. 

Note that $P(t)+P(u_0-t) \le P(u_0/2)+P(u_0) =1$, and that the function 
$y(2-y)$  is increasing in the range $0\le y \le 1$.  Hence 
$$
\align
\frac{1}{u_0} \int_0^{u_0/2} (P(t)+P(u_0-t)) &(2-P(t)-P(u_0-t)) dt 
\ge \frac{1}{u_0} \int_0^{u_0/2} m(t) (2-m(t)) dt \\
&=(12-4\alpha)\frac{e^{\alpha}}{\sqrt{e}} - 13 -2\alpha^2 + (6-2\alpha)
e^{-\alpha/2},\\
\endalign
$$
after some calculations.
This function of $\alpha$ attains a unique minimum 
in the range $(0,1/2)$, at $\alpha_0=0.08055\ldots$, at which point
its value is $\geq  0.272516916\ldots \ge 112/411$.  Inserting this 
into Lemma 7.4, and taking square roots of both sides 
we obtain the result.
\enddemo

For convenience, in the next three subsections we put $\lambda =\lambda_\theta
= (28/411) \cos^2\theta$.

\subhead 7d.  Proof of Theorem 5 \endsubhead

\noindent For all $u$, we seek to show that the distance of $\sigma(u)$ 
from $\lambda$ is $\le 1-\lambda$.  Suppose $u\ge u_0$.  By the triangle
inequality the distance of $\sigma(u)$ from $\lambda$ is $\le \lambda$ plus
the distance from $\sigma(u)$ to $0$.  By Proposition 7.5, the latter 
distance is $\le 1-2\lambda$, so that our claim holds in this case.

Suppose $u\le u_0$.  Observe that $\min\{ t,2\sqrt{t(1-t)}\} \leq 2-t-(2+t)(28/411)$
for all $0\leq t\leq 1$. Taking $t=P(u)$, multiplying through by $R_1(u)$ and observing that $P(u/2)\leq 1-t$, we obtain 
$R_2(u)+C_2(u)\le R_1(u)(2-P(u))-2\lambda (1-\RE\sigma(u))$,
from (7.4) and the second inequality in (7.6). By (7.7) we deduce that
$2\lambda(1-\RE\sigma(u))\le 1-|\sigma(u)|^2$ and so, re-arranging,
$(\RE \sigma(u)-\lambda)^2+(\IM \sigma(u))^2 \le (1-\lambda)^2$.

It follows that $\Lambda(S)$ is contained in the 
circle centered at $\lambda$ with radius $1-\lambda$, and 
Theorem 5 follows since $\Gamma(S) \subset [0,1]\times \Lambda(S)$.

\subhead 7e. Proof of Theorem 6(ii) \endsubhead

\noindent We shall 
show that Ang$(\sigma(u)) \le \frac{\pi}2 - \frac{\sin \delta}{2}$.  
Suppose first that $u\ge u_0$.  Note that
 Ang$(\sigma(u)) \le \arcsin(|\sigma(u)|) \le \arcsin (1-2\lambda)$, 
by Proposition 7.5.  Now $\arcsin (1-2\lambda)\le \pi/2 -\sqrt{4\lambda}$,
and our claim follows in this case since $\cos\theta =\sin\delta$, and
$\sqrt{112/411}>1/2$.  

Thus we may suppose $u<u_0$.  By definition, Ang$(\sigma(u))= \arctan(
|\IM \sigma(u)|/(1-\RE \sigma(u)))$.  By (7.5) and (7.6),
$$
\frac{|\IM \sigma(u)|}{(1-\RE\sigma(u))}
\le \frac{\sqrt{R_1(u)(P(u)-R_1(u))}}{R_1(u)(1-P(u)/2)} 
\le 2\sqrt{P(u)/R_1(u)-1} \le 2\tan\theta,
$$
since $P(u)\le 1$ as $u<u_0$, and $P(u)/R_1(u)\le \sec^2 \theta$ by Lemma 7.2.  

For $0\le x< 1$ we have $(1+x)/(1-x) \ge 4x/(1-x^2)$.  Taking 
$x=\tan(\theta/2)$ we deduce that $2\tan\theta \le \tan(\pi/4+\theta/2)$.  
Thus Ang$(\sigma(u))\le \arctan(2\tan \theta) \le \pi/4+\theta/2 = \pi/2
-\delta/2 \le \frac{\pi}{2} - \frac{\sin \delta}{2}$, as desired.  
We have shown that Ang$(\Lambda(S)) \le \frac{\pi}{2} -\frac{\sin \delta}{2}$,
and Theorem 6(ii) follows.

\subhead 7f. Proof of Theorem 7(ii) \endsubhead

\noindent We show that the projection of $\sigma(u)$ on $S$ is 
$\le 1- 2\lambda$.  From this it follows that the maximal projection of 
$\Lambda(S)$ (and hence of $\Gamma(S)$ by Theorem 3$'$) is $\le 1-2\lambda$, proving 
Theorem 7(ii).  If $u\ge u_0$ then the projection of $\sigma(u)$ on $S$ is 
$\le |\sigma(u)|\le 1-2\lambda$, by Proposition 7.5, and our claim follows.
Thus we may suppose that $u<u_0$.  
Since Ang$(S)=\theta =\pi/2-\delta$, we need to show that $\RE 
(e^{-i\gamma} \sigma(u)) \le 1-2\lambda$ for $2\delta \le |\gamma| \le 
\pi$ (taking the projection along $\zeta=e^{i\gamma}$).

Recall that by (7.7), $|\sigma(u)|^2 
\le 1-2R_1(u) + R_1(u)P(u) + R_2(u)+C_2(u)$.  Using (7.4) together with 
the bound $P(u/2)\le 1-P(u)$, we deduce that, 
since $R_1(u)\ge P(u)\cos^2\theta$ by Lemma 7.2, 
$$
\align
|\sigma(u)|^2 &\le 1 - 2R_1(u)(1-P(u)/2 -\min\{ P(u)/2,\sqrt{P(u)(1-P(u))}\}
)\\
&\le 1- 2\cos^2\theta P(u)(1-P(u)/2-\min\{P(u)/2,\sqrt{P(u)(1-P(u))}\}) \\
&\le 1-4\lambda\\
\endalign
$$
in the range $1/6\le P(u)\le 1$, as may be verified using Maple.  Thus 
in this range of $P(u)$, the projection of $\sigma(u)$ on $S$ is 
$\le |\sigma(u)|\le 1-2\lambda$, as desired.

Now suppose $P(u)\le 1/6$.  From the above argument we know that 
$|\sigma(u)|^2 \le 1- 2R_1(u)(1-P(u))\le 1- 5R_1(u)/3$, so that $|\sigma(u)|
\le 1-5R_1(u)/6 \le 1-2\lambda$ if 
$R_1(u)> 12\lambda/5 = (112/685)\cos^2\theta$.

So we are left with the case $P(u)\le 1/6$, and $R_1(u)\le \cos^2\theta/6 
\le 1/6$.  By (7.5),
$|\IM \sigma(u)| \le \sqrt{R_1(u)P(u)} \le (1/6)\cos \theta$,
and by (7.6), $\RE \sigma(u) \ge 1-R_1(u)(1+P(u)/2)  
\ge 1-(1/6)(13/12) =59/72$.  Hence 
$\tan (|\arg \sigma(u)|) = |\IM \sigma(u)|/\RE \sigma(u) < \cos \theta \le \cot \theta =\tan\delta$.  
Thus $|\arg \sigma(u)|\le \delta$, and so if $2\delta \le |\gamma|\le \pi$,
$|\arg (e^{-i\gamma} \sigma(u))|>\delta$.  So the 
projection of $\sigma(u)$ on $e^{i\gamma}$ is $\le \cos \delta =\sin \theta
\le 1-(1/2)\cos^2 \theta$.  This completes the proof of Theorem 7(ii).

\head 8. Generalized notions of the spectrum:  The Logarithmic spectrum
\endhead

\noindent We may generalize the notion of spectrum by considering the 
values
$$ 
\biggl(\sum_{n\le N} \kappa (n)\biggr)^{-1} \sum_{n\le N} f(n) \kappa (n)
$$
for $f\in {\Cal F}(S)$ as $N\to \infty$, where $\kappa (n)$ is a given
positive valued function (we considered the case $\kappa =1$ above).
In this setting one quickly becomes curious about the 
weights $\kappa (n)=1/n^\sigma$ for a given real number $\sigma \ge 0$. 
If $\sigma>1$ then the
sum converges absolutely and so we obtain the set
of Euler products $\zeta(\sigma)^{-1} \prod_p (1-f(p)/p^{\sigma})^{-1}$.
If $\sigma <1$ then the new spectrum is exactly the same as 
$\Gamma(S)$, since if $f\in {\Cal F}(S)$ is completely multiplicative 
then, for any given $\sigma<1$, we have 
$$
\biggl(\sum_{n\le x} \frac 1{n^{\sigma}}\biggr)^{-1} \sum_{n\le x} \frac{f(n)}
{n^{\sigma}} = \frac 1x \sum_{n\le x} f(n) + o(1) .\tag{8.1}
$$
To see this note that if $\sum_{p\leq x} (1-\RE f(p))/p \to \infty$
then both sides of the equation are $o(1)$ by Lemma 1' and partial summation.
Thus we may assume that 
$\sum_{p\leq x} |1-f(p)|/p \asymp_S \sum_{p\leq x} (1-\RE f(p))/p \ll 1$.
Let $g(p^k)=f(p^k)-f(p^{k-1})$ for each prime power. By (4.1) we have 
$\sum_{n\le t} |g(n)| \ll (t/\log t)\exp( \sum_{p\leq t} |g(p)|/p ) \ll t/\log t$ when $t\leq x$; and so
$\sum_{d\le x} |g(d)|/d^{\sigma} \ll x^{1-\sigma}/(1-\sigma)\log x$
by partial summation. Therefore, since
$\sum_{n\le t} n^{-\sigma} = t^{1-\sigma}/(1-\sigma) + O(1)$,
we obtain
$$\align
\sum_{n\le x} \frac{f(n)}{n^{\sigma}} &= \sum_{d\le x} \frac{g(d)}{d^{\sigma}}
\sum_{n\le x/d} \frac{1}{n^{\sigma}} = \sum_{d\le x} \frac{g(d)}{d^{\sigma}}
\biggl(\frac{1}{1-\sigma}\biggl(\frac{x}{d}\biggr)^{1-\sigma} +O(1) \biggr) \\
&=\frac{x^{1-\sigma}}{1-\sigma}\sum_{d\le x} \frac{g(d)}{d}+O\biggl(
\sum_{d\le x}\frac{|g(d)|}{d^{\sigma}}\biggr) \\
&=\biggl( \sum_{n\leq x} \frac 1{n^\sigma} \biggr) \biggl\{
\sum_{d\le x} \frac{g(d)}{d} + O\biggl( \frac{1}{\log x}\biggr)
\biggr\} . \\
\endalign
$$
Comparing the formula at $\sigma$ with the formula at $\sigma=0$ gives
(8.1).

This leaves us with the case $\sigma=1$, that is $\kappa(n)=1/n$, which 
gives rise to the logarithmic spectrum $\Gamma_0(S)$ mentioned 
in the introduction. We now proceed to a study of this spectrum, 
beginning with some general results on logarithmic means.  
Elsewhere we will apply these methods to 
obtain upper bounds on $L(1,\chi)$.\footnote{In the 
spirit of P.J. Stephens [12] who showed that $|L(1,\chi_d)| \le 
\frac{1}{4}(2-\frac{2}{\sqrt{e}}+o(1)) \log |d|$ where $\chi_d$ is a
quadratic character with conductor $|d|$.  We establish similar 
results for higher order characters.}  


One may also consider other other choices of $\kappa(n)$; for 
example, $\kappa(n)=d_k(n)$, the $k$th divisor function.
It would be interesting to determine this spectrum when $S=\{\pm 1\}$.

\subhead 8a.  Generalities on logarithmic means \endsubhead

\proclaim{Proposition 8.1}  Let $f$ be a multiplicative function with $|f(n)|
\le 1$ for all $n$, and put $g(n) =\sum_{d|n} f(d)$.  Then 
$$
\frac{1}{\log x} \Big| \sum_{n\le x} \frac{f(n)}{n} \Big| \le 2e^{2\gamma} 
\prod_{p\le x} \Big(1-\frac 1p \Big)^2 \Big(1 + \frac{|g(p)|}{p} 
+ \frac{|g(p^2)|}{p^2} + \ldots \Big) +O \Big(\frac{1}{\log x}\Big).
$$
\endproclaim 

\demo{Proof}  Since
$$
\sum_{n\le x} g(n)  = \sum_{n\le x} \sum_{d|n} f(d) 
= \sum_{d\le x} f(d) \Big(\frac{x}{d} + O(1)\Big) = 
x \sum_{d\le x} \frac{f(d)}{d} + O(x),
$$
we see that 
$$
\frac{1}{\log x} \Big| \sum_{n\le x} \frac{f(n)}{n} \Big| \le \frac{1}{x\log x}
\sum_{n\le x} |g(n)| + O\Big(\frac{1}{\log x}\Big).
$$
Note that $|g(n)|$ is a non-negative multiplicative function with $|g(n)| 
\le d(n)$ for all $n$.  Hence by Theorem 2 of Halberstam and Richert [4] (see 
(4.1) above) we obtain
$$
\align
&\frac{1}{x\log x} \sum_{n\le x} |g(n)| \le \frac{2}{x\log^2 x} 
\sum_{n\le x} \frac{|g(n)|}{n} + O\Big(\frac{1}{\log x}\Big)
\\
\le &\frac{2}{x\log^2 x} \prod_{p\le x} \Big(1+ \frac{|g(p)|}{p} + 
\frac{|g(p^2)|}{p^2} +\ldots \Big) + O\Big(\frac{1}{\log x}\Big).
\\
\endalign
$$
The result follows from Mertens' theorem. 
\enddemo

Since $2-|1+z| \le (1-\RE z)/2$ whenever $|z|\leq 1$, 
the right side of the equation in Proposition 8.1 is 
$$
\ll \exp\Big( - \sum_{p\le x} \frac{2-|g(p)|}{p} \Big) \le \exp \Big( -\frac{1}{2} \sum_{p\le x} \frac{1-\RE f(p)}{p}\Big).
$$
More precisely one obtains
$$
\frac{1}{\log x} \Big| \sum_{n\le x} \frac{f(n)}{n} \Big| \le 
\frac{26e^{2\gamma}}{\pi^2} 
\exp \Big( -\frac{1}{2} \sum_{p\le x} \frac{1-\RE f(p)}{p}\Big) ,
$$
a weak, but relatively easy and effective, analogue of Lemma 1$'$
for logarithmic means.  Moreover this has the advantage that $f$ need not be
restricted to a subset of ${\Bbb {\Bbb U}}$ since the case $f(n)=n^{it}$ does not impede
us here (since $\sum_{n\leq x} n^{i\alpha -1} \ll 1$).

Next we derive analogues of Propositions 4.1, 4.4, and 4.5.  As the above 
example indicates, the situation here is much simpler.  
For example, the analogue of Proposition 4.1 is the trivial estimate
$$
\frac{1}{\log x} \sum_{n\le x} \frac{f(n)}{n} -\frac{1}{\log (x/y)} 
\sum_{n\le x/y}\frac{f(n)}{n} \ll \frac{\log 2y}{\log x}, 
$$
which is valid for all functions $f$ with $|f(n)|\le 1$, and all 
$1\le y\le \sqrt{x}$.  
Using this estimate (in place of Proposition 4.1) and 
arguing exactly as in the proof of Proposition 4.5
we arrive at the following Proposition (see also Lemma 5 of Hildebrand [11]).

\proclaim{Proposition 8.2}  Let $f$ be any multiplicative function 
with $|f(n)|\le 1$.  Let $g$ be the ccompletely multiplicative function 
defined by $g(p)=1$ for $p \le y$ and $g(p)=f(p)$ for $p>y$.  Then
$$
\frac{1}{\log x} \sum_{n\le x} \frac{f(n)}{n} = \Theta(f,y)
\frac{1}{\log x} \sum_{n\le x} \frac{g(n)}{n} + O\biggl(\frac{\log y}{\log x}
\exp(s(f,y))\biggr),
$$
where $s(f,y)=\sum_{p\le y}|1-f(p)|/p$.  The remainder 
term above is $\ll (\log y)^3/\log x$. 
\endproclaim


We omit the proof of Proposition 8.2 since it is almost identical to that of Proposition 4.5.

Observe that 
$$
\frac{1}{\log y^u} \sum_{n\le y^u} \frac{f(n)}{n} 
= \frac{1}{u} \int_0^u \frac{1}{[y^t]} 
\sum_{n\le y^t} f(n) dt + O\biggl(\frac{1}{u\log y}\biggr),
$$
which implies that $\Gamma_0(S)$ lies inside the convex hull of
$\Gamma(S)$.
 From this equation, we deduce the following analogues of Proposition 1 and its converse.

\proclaim{Proposition 8.3}  Let $f$ and $\chi$ be as in Proposition 1.  
Then 
$$
\frac{1}{\log y^u} \sum_{n\le y^u} \frac{f(n)}{n} = \frac{1}{u} \int_0^u 
\sigma(t) dt + O\biggl(\frac{u}{\log y}\biggr).
$$
\endproclaim

\proclaim{Proposition 8.3 (Converse)} Let $f$ and $\chi$ be as 
in the converse of Proposition 1.  Then for all $1/\log y\le t \le u$
$$
\frac{1}{t} \int_0^t \sigma(v) dv = \frac{1}{\log y^t} \sum_{n\le y^t} 
\frac{f(n)}{n} + O (u^{\epsilon}-1) + O\biggl(\frac{u}{\log y}\biggr).
$$
\endproclaim

Let $S$ be a closed subset of ${\Bbb U}$ with $1\in S$.  
We define $\Lam_0(S)$ 
to be the set of values $\frac 1u\int_0^u \sigma(t) dt$ obtained as follows:
Let $\chi$ be any element of $K(S)$, and 
let $\sigma$ denote the corresponding solution to (1.5).  Then $\Lam_0(S)$
is the set of all values $\frac 1u \int_0^u \sigma(t) dt= \frac 1u 
(1*\sigma)(u)$ for 
all $u>0$, and all such choices of $\chi$.  Note that $\Lam_0(S)$ is
in the convex hull of $\Lam(S)$.

Combining Proposition 8.2, with Proposition 8.3 and its Converse, we 
obtain the following Structure Theorem for the logarithmic spectrum. 

\proclaim{Theorem 8.4} Let $S$ be a closed subset of ${\Bbb U}$ with $1\in S$. 
Then $\Gamma_0(S) = \GT \times \Lam_0(S)$.  Further $\Lam_0(S) 
=\Lam_0(S) \times {\Cal E}(S)$, and so  
$$
\Lam_0(S) \subset \Gamma_0(S) \subset \Lam_0(S) \times [0,1].
$$  
\endproclaim

Theorem 8.4 is proved exactly in the same way as Theorems 3 and 3$'$; 
so we omit its proof.  We end this subsection by making the 
following useful observation:\footnote{More generally, 
$u(\sigma_1*\sigma_2)(u) = ((\chi_1+\chi_2)*(\sigma_1*\sigma_2))(u)$.}
$$
\align
u(1*\sigma)(u) &= \int_0^u (u-t)\sigma(t) dt + \int_0^u t\sigma(t) dt 
= (1*(1*\sigma))(u) + (1*(t\sigma(t)))(u)
\\
&= (1*1*\sigma)(u) + (1*\chi*\sigma)(u)
= ((1*\sigma)*(1+\chi))(u). \tag{8.2} \\
\endalign
$$

\subhead 8b. Bounding $\Gamma_0(S)$: Proof of Theorem 8 \endsubhead

\noindent If $S=\{1\}$ then $\Gamma_0(S)={\Cal R}=\{1\}$, and 
there is nothing to prove.  Suppose that $S$ contains an element $\alpha\neq 
1$.  Then $(\frac{1+\alpha}{2})^n \in {\Cal R}$ for all $n\ge 1$.   As 
$n\to \infty$ this sequence of points converges to $0$,  
and since ${\Cal R}$ is closed,  we deduce that $0\in {\Cal R}$.  By convexity
it follows that ${\Cal R}={\Cal R}\times [0,1]$.  Hence, by Theorem 8.4,
 we need only show that 
$\Lam_0(S) \subset {\Cal R}$  in order to establish Theorem 8.

We define, for any complex number $z$, its ${\Cal R}$-norm
$\| z\|_{\Cal R} := \min_{r\in {\Cal R}} |z-r|$; that is $\|z\|_{\Cal R}$ 
is the shortest distance from $z$ to ${\Cal R}$.  We first make a 
couple of general observations about this norm:

Let $X$ be a measurable subset of 
the real line, and suppose $f$ is a non-negative measurable function with 
$\int_X f(x)dx =1$.  Then for any measurable function $g$,
\footnote {An analogous convexity result holds for 
sums: If $a_i\ge 0$ with $\sum a_i=1$, then 
$\|\sum a_i z_i\|_{\Cal R} \le \sum a_i \|z_i\|_{\Cal R}$.} 
$$
\Big\| \int_X f(x) g(x)\ dx \Big\|_{\Cal R} \le \int_X f(x)  \|g(x)\|_{\Cal R} 
\ dx.
\tag{8.3}
$$
To see (8.3), suppose $r(x)$ is a point in ${\Cal R}$ closest to $g(x)$.  
Then $\int_X f(x) r(x)dx $ is a convex combination 
of the points $r(x)$, and so is an element of ${\Cal R}$.  Therefore
$$
\align
\Big\|\int_X f(x)g(x)dx\Big\|_{\Cal R} &\le \Big|
\int_X f(x)g(x) dx - \int_X f(x)r(x)dx\Big|
\le \int_X f(x) |g(x)-r(x)|dx \\
&= \int_X f(x) \|g(x)\|_{\Cal R} dx,
\\
\endalign
$$
which proves (8.3).

Let $s$ be any point in the convex hull of $S$ and let 
$r$ be a point in ${\Cal R}$ closest to given $z$. By the definition 
of ${\Cal R}$,  we know that $r\frac{1+s}{2}$ is also a point in ${\Cal R}$, 
and so
$$\left\|z \frac{1+s}{2}\right\|_{\Cal R} \le \left|z\frac{1+s}{2}-r\frac{1+s}{2}\right| 
= \left|\frac{1+s}{2}\right| \|z\|_{\Cal R}\le \|z\|_{\Cal R}. \tag{8.4}
$$

Suppose $\chi \in K(S)$ is given and let $\sigma$ 
be the corresponding solution 
to (1.5).  We shall show that $\frac 1u 
(1*\sigma)(u) \in {\Cal R}$ for all $u$.  
This proves that $\Lam_0(S)$ (and so $\Gamma_0(S)$) is contained in 
${\Cal R}$. 
Define $\alpha(u)= - u\| \frac{1}{u} (1*\sigma)(u)\|_{\Cal R}$.  
Plainly $\alpha(u)=0$ for $u\le 1$, and we shall show below that it 
is always non-negative so that $\alpha(u)=0$ for all $u$, which 
proves that $\frac{1}{u}(1*\sigma)(u)\in {\Cal R}$. 

By (8.2) we see that 
$$
\Big\|\frac{1}{u}(1*\sigma)(u)\Big\|_{\Cal R} 
=\Big\| \frac{1}{u^2} \int_0^u 2v \Big(\frac{1}{v} (1*\sigma)(v)\Big) 
\frac{1+\chi(u-v)}{2} dv \Big\|_{\Cal R}.
$$
Applying (8.3) with $X=[0,u]$, and $f(x)=2v/u^2$, we deduce that 
the above is 
$$ 
\le \frac{1}{u^2} \int_0^u 2v \Big\| 
\Big(\frac{1}{v}(1*\sigma)(v)\Big) 
\frac{1+\chi(u-v)}{2}\Big\|_{\Cal R} dv, 
$$
which by (8.4) is 
$$
\le \frac{1}{u^2} \int_0^u 2v \Big\| \frac 1v (1*\sigma)(v)\Big\|_{\Cal R} dv.
$$
It follows that $u\alpha(u)\ge (2*\alpha)(u)$, and so by Lemma 3.1, 
$\alpha(u)$ is always non-negative,
as desired.  This completes the proof of Theorem 8.

\subhead 8c.  Proof of Corollary 4 \endsubhead

\noindent If $S=[-1,1]$ then $\frac{1+s}{2} \in [0,1]$ for all $s\in S$, 
and so it follows that ${\Cal R} = [0,1]$ here.  Hence $\Gamma_0([-1,1]) 
\subset [0,1]$.  Since $\Gamma_0([-1,1]) \supset {\Cal E}([-1,1]) 
= [0,1]$, it follows that $\Gamma_0([-1,1]) = [0,1]$, proving part (i).

Part (ii) is proved in the same way as Corollary 3(ii):  Take 
$\chi(t)=1$ for $t\le 1$ and $\chi(t)=\alpha$ for $t>1$.  Then 
by Theorem 3.3, $\sigma(t)=1-(1-\alpha)\log t$ for $1\le t\le 2$.  
Hence, for $1\le u\le 2$, 
$$
\frac{1}{u}\int_0^{u} \sigma(t) dt  =1 - (1-\alpha) \frac{1}{u}
\int_1^u \log t dt 
= 1- (1-\alpha)\Big(\log u -1 +\frac 1u\Big) \tag{8.5}
$$
belongs to $\Lam_0(S)$, as desired.  If $0<\text{Ang}(S) <\pi/2$, then 
take $1<u \le 2$, and $\alpha \in S$ such that Ang$(\alpha)=
\text{Ang}(S)$. The argument given in \S 6b (proof of Corollary 3(ii)) 
shows that such elements of $\Lam_0(S)$ are not in $\GT$.

Notice that taking $\alpha \in S$ such that Ang$(\alpha)=$ Ang$(S)$ 
in the construction (8.5), we obtain that 
Ang$(S) \le $ Ang$(\Gamma_0(S))$.\footnote 
{Alternately, 
this follows because Ang$(\Gamma_0(S)) \ge $ Ang$({\Cal E}(S)) =$ Ang$(S)$.}  
We now show that Ang$({\Cal R}) \le \text{Ang}(S)$, so that 
by Theorem 8, we have $\text{Ang}(\Gamma_0(S)) =\text{Ang}(S)$.  
Suppose Ang$(S)=\frac{\pi}{2}-\delta$, so that $S$ is contained in the 
convex hull of $\{1\} \cup \{e^{i\theta}: \ \ 2\delta \le |\theta| \le \pi\}$.
Each product $\prod_{j=1}^{m} \frac{1+s_j}{2}$, where $s_j$ is in the 
convex hull of $S$, is easily expressed as a convex combination 
of elements of the form $\prod_{j=1}^{n} \frac{1+e^{i\theta_j}}{2}$ 
where $2\delta \le |\theta_j|\le \pi$.  Hence ${\Cal R}$ is 
contained in the convex hull of $1$ and points of the 
form $\prod_{j=1}^{n} \frac{1+e^{i\theta_j}}{2} 
= \prod_{j=1}^{n} \cos (\theta_j/2) e^{i\theta_j/2}$ where 
$2\delta \le |\theta_j|\le \pi$.  Such a product has magnitude $\le 
(\cos \delta)^n \le \cos \delta$ if $n\ge 1$.  Thus ${\Cal R}$ 
is in the convex hull of $\{1\} \cup \{|z|\le \cos \delta\}$.  
If $|z|\le \cos\delta$ then Ang$(z) \le \arcsin(|z|) \le \frac{\pi}{2} 
-\delta$, and so it follows that Ang$(\Cal R) \le \frac{\pi}{2}-\delta$.  
This proves (iii).

To prove (iv), we first observe that $f(x):= (\cos x)^{\frac{1}{x}}$ 
is decreasing in $(0,\frac{\pi}{2}]$.  Differentiating $f$ logarithmically, 
we need to show that $-(\log \cos x)/x^2 -\tan x/x \le 0$, or
equivalently, that $g(x):=x\tan x + \log \cos x \ge 0$. Now $g'(x) =x\sec^2 x$ 
is positive in $(0,\frac{\pi}{2}]$, and so $g(x) \ge g(0) =0$, as desired.  
It follows that if $\delta \le \theta   \le \frac{\pi}{2}$ then 
$\cos \theta \le (\cos \delta)^{\frac{\theta}{\delta}}$.

 From the proof of (iii), we 
know that ${\Cal R}$ is contained in the convex hull of 
$1$ and products of the form $\prod_{j=1}^{n} \cos (\theta_j/2) 
e^{i\theta_j/2}$ where each $\theta_j \in [2\delta,{\pi}]$.  
If such a product has argument $\nu$, then we must have $\sum_{j=1}^{n}
\theta_j \ge 2\nu$.  By the previous paragraph, the 
magnitude of such a product is $\le \prod_{j=1}^n (\cos 
\delta)^{\frac{\theta_j}{2\delta}} \le (\cos \delta)^{\frac{\nu}{\delta}}$.  
Thus ${\Cal R}$ is contained in the set $\{z: \ \ |z| \le (\cos \delta)^
{\frac{|\arg z|}{\delta}}\}$, which proves (iv).

\head 9. Quadratic residues and nonresidues revisited:  Proof of Theorem 9
 \endhead

\noindent  Throughout this section $D$ denotes a 
fundamental discriminant.
 
\proclaim{Proposition 9.1} Let $B$ be fixed, and $X$ be large, 
and suppose $1\le z\le \frac 14(\log X)$. 
Let $f(n)$ be a completely multiplicative function 
satisfying $f(p)=\pm 1$ for $p\le z$, and $f(p)=0$ for $p>z$.  Put 
$P = 4\prod_{p\le z} p$ and let $a \pmod {P}$ be an arithmetic progression 
(with $a\equiv 1$, or $5 \pmod 8$) such that $\fracwithdelims() {a}{p} 
=f(p)$ for each $p\le z$.  With ${\Cal N}(X;a,P)$ denoting the 
number of fundamental discriminants $0< D\le X$ with $D\equiv a\pmod P$, 
we have 
$$
\frac{1}{{\Cal N}(X;a,P)}  \sum\Sb 0 <D \le X\\
                                D\equiv a\pmod P\endSb
\sum_{n\le (\log X)^B} \fracwithdelims() {D}{n} = 
\sum_{n\le (\log X)^{B} } f(n) + O\Big(\frac{(\log X)^B}{z}\Big).
$$
\endproclaim

\demo{Proof}  We write $n=rs$ where each prime dividing $r$ is $\le z$, and
each prime dividing $s$ is $>z$.  Thus $\fracwithdelims() {D}{n}
=f(r) \fracwithdelims(){D}{s}$, and so 
$$
\sum\Sb 0< D\le X \\ D\equiv a\pmod P\endSb 
\sum_{n\le (\log X)^B} \fracwithdelims() {D}{n} 
= \sum\Sb r\le (\log X)^B\\ 
        p|r \implies p \le z \endSb f(r) 
\sum\Sb s\le (\log X)^B/r\\
        p|s \implies p>z\endSb \sum\Sb 0< D\le X\\ D\equiv a\pmod P\endSb
\fracwithdelims() {D}{s}.\tag{9.1}
$$

The $s=1$ term in (9.1) contributes 
$$
\sum\Sb r\le (\log X)^B \\ p|r\implies p\le z\endSb f(r) {\Cal N}(X;a,P)
= {\Cal N}(X;a,P) \sum_{n\le (\log X)^B} f(n). \tag{9.2}
$$
The terms $s>1$ with $s=\square$ ($\square$ denotes the square of an integer) 
contribute
$$
\align
&\ll {\Cal N}(X;a,P) \sum\Sb r\le (\log X)^B \\ p|r\implies p\le z\endSb 
\sum\Sb 1<s \le (\log X)^B/r\\ s\in {\Bbb Z}^2 \\ p|s\implies p>z \endSb 1 \\
&\ll {\Cal N}(X;a,P) \sum\Sb r\le (\log X)^B/z^2\endSb 
\sqrt{\frac{(\log X)^B}{r} } 
\ll {\Cal N}(X;a,P) \frac{(\log X)^B}{z}.\tag{9.3} \\
\endalign
$$
Finally we consider the contribution of the terms $s\neq \square$ to (9.1). 
For such an $s$, $\fracwithdelims(){\cdot}{s}$ is a non-principal 
character of conductor $\le s$, and so we may expect substantial 
cancellation in the sum over $D$ in (9.1).  Indeed,
we have using $\mu(m)^2 =\sum_{l^2 |m} \mu(l)$
$$
\sum\Sb 0< D\le X\\ D\equiv a\pmod P\endSb \fracwithdelims() {D}{s} 
= \sum\Sb m \le X \\ m \equiv a\pmod P\endSb \mu(m)^2 \fracwithdelims(){m}{s}
=\sum\Sb l \le \sqrt{X} \endSb \mu(l) 
\sum\Sb m\le X\\ m\equiv a\pmod P\\ l^2 |m \endSb \fracwithdelims() {m}{s}.
$$
By the P{\' o}lya-Vinogradov inequality the 
inner sum over $m$ above is $\ll \sqrt{s} \log s$.  Hence the 
sum over $D$ above is $\ll \sqrt{Xs} \log s$.  This demonstrates that 
the $s\neq \square$ terms in (9.1) contribute an amount 
$$
\ll \sum\Sb r\le (\log X)^B \endSb \sum\Sb s\le (\log X)^B/r\endSb 
\sqrt{Xs}\log s \ll X^{\frac 12 +\epsilon}. 
$$                              
Combining this with the estimates (9.2), and (9.3), we see by (9.1) that
$$
\frac{1}{{\Cal N}(X;a,P)}  \sum\Sb 0 <D \le X\\
                                D\equiv a\pmod P\endSb
\sum_{n\le (\log X)^B} \fracwithdelims() {D}{n} = 
\sum_{n\le (\log X)^{B} } f(n) + O\Big(\frac{(\log X)^B}{z} + \frac{X^{\frac12 
+\epsilon}}{{\Cal N}(X;a,P)}\Big).
$$
Since 
$$
{\Cal N}(X;a,P) \sim \frac{X}{P} \frac{6}{\pi^2} \prod_{p|z} 
\Big(1-\frac{1}{p^2}\Big)^{-1},
$$
the second error term above is $\ll P/X^{\frac 12 -\epsilon}$. 
Since $z\le \frac{1}{4}\log X$, $P\ll X^{\frac 14 +\epsilon}$ by the 
prime number theorem, and so the second error term above is $\ll X^{-\frac 14
+\epsilon}$, which is subsumed by the error term of the Proposition.  

\enddemo

Armed with Proposition 9.1, we 
now show that $\beta(B)\le \gamma(B)$ for all $B$.  
Let $X$ be large, and choose $z=\frac 14 \log X$.  By Proposition 9.1, we know 
that there is a fundamental 
discriminant $D$ with $X/\log X\ll D \le X$, such that
$$
\frac{1}{(\log D)^B} \sum_{n\le (\log D)^B} \fracwithdelims() {D}{n} 
\le \frac{1}{(\log D)^B} \sum_{n\le (\log D)^B} f(n) +o(1), \tag{9.4}
$$
where $f$ is any completely multiplicative function as in Proposition 9.1.  
Suppose we are given $\chi \in {\Cal C}(u)$.  Put $y=z^{\frac 1u}$,
and choose $f \in {\Cal F}(\{0,\pm 1\})$ as in the converse of Proposition 1.  
Thus choose $f$ so that $f(p)=1$ for $p\le y$, $f(p)=0$ for $p>y^u=z$, 
and such that for almost all $0\le t\le u$,
$$
\Big| \chi(t) - \frac{1}{\vartheta(y^t)} \sum_{p\le y^t} f(p)\log p\Big|
\le \epsilon.
$$
 From Proposition 1 (Converse) it follows that the right hand side of 
(9.4) is $\sigma(Bu) + O(u^{B\epsilon}-1) + O(u/\log y) +o(1)$.  
Letting $\epsilon \to 0$, and $X\to \infty$ (so that $y\to \infty$), it 
follows that $\beta(B) \le \sigma(Bu)$.  Now varying $u$, and $\chi \in {\Cal 
C}(u)$, we deduce that $\beta(B) \le \gamma(B)$.

To complete the proof of Theorem 9, it remains now to show that 
$-\rho(B)\le \gamma(B) <0$.  We first show that 
$$
|\sigma(Bu)| \le \rho(B) \qquad \text{for all } B \text{ and all } 
\chi \in {\Cal C}(u). \tag{9.5}
$$
To prove (9.5), suppose $\chi \in {\Cal C}(u)$ is given, and 
put $a(B)= \rho(B) - |\sigma(Bu)|$.  Since $\rho(B)=1$ for $B\le 1$, 
it follows that $a(B)\ge 0$ for $B\le 1$.  Define $b(t)=1$ for 
$t\le 1$, and $b(t)=0$ for $t>1$.  From the definition of 
the Dickman function $B\rho(B)=(b*\rho)(B)$, and so we have 
$$
\align
B \ a(B)& = B\rho(B) - |B\sigma(Bu)| = (b*\rho)(B)
- \frac{1}{u} \Big|\int_{(B-1)u}^{Bu} \sigma(t) \chi(Bu-t)dt\Big|
\\
&\ge (b* \rho)(B) - \int_{B-1}^{B} |\sigma(ut)| dt
= (b*a)(B).
\\
\endalign
$$
 From Lemma 3.1 it follows that $a(B)\ge 0$ always, which establishes (9.5).

 From (9.5) we see that $-\rho(B) \le \gamma(B)$, and 
it remains now to show that $\gamma(B) <0$.  We prove this by considering the 
following example:  Put $\chi_{-}(t) =1$ for $t\le 1$, $\chi_{-}(t)=-1$ 
for $1\le t\le 2$, and $\chi_{-}(t)=0$ for $t>0$, so that $\chi_{-} \in 
{\Cal C}(u)$ for all $u\ge 2$.  Hence, if $\sigma_{-}$ denotes the 
solution to $w\sigma_{-}(w) = \sigma_{-}*\chi_{-}$ then 
$\gamma(B) \le \min_{u\ge 2} \sigma_{-}(Bu) = \min_{w\ge 2B} \sigma_{-}(w)$.
We will now show that $\sigma_{-}(w)$ changes sign infinitely often; hence
there are arbitrarily large $w$ with $\sigma_{-}(w)<0$ which 
shows that $\gamma(B)<0$ for all $B$.

Suppose $\sigma_{-}(w)$ maintains sign from some point on: precisely,
suppose $|\sigma_{-}(w_0)|>0$ and that $\sigma_{-}(w)$ has 
the same sign as $\sigma_{-}(w_0)$
for all $w\ge w_0$.  Define $F(w) =\int_{w-1}^{w} \sigma_{-}(t) dt$.
Note that 
$$
w\sigma_{-}(w) = (\sigma_{-}*\rho_{-})(w) = 
\int_{w-1}^w \sigma_{-}(t) dt -\int_{w-2}^{w-1} \sigma_{-}(t)dt = 
F(w)-F(w-1). \tag{9.6}
$$
Since $F(w)$ has the same sign as $\sigma(w_0)$ for all $w\ge w_0+1$, 
we deduce from (9.6) that $|F(w+1)| = |F(w)+(w+1)\sigma_-(w+1)| 
=|F(w)| + |(w+1)\sigma_-(w+1)| \ge |F(w)|$ for all $w\ge w_0+1$.  
Hence
$$
\liminf_{n\to \infty} |F(w_0+n)| \ge |F(w_0+1)| = \int_{w_0}^{w_0+1} 
|\sigma_-(t)| dt >0.
$$
However, from (9.5) we see that 
$$
|F(w)| = \int_{w-1}^{w} |\sigma_{-}(t)| dt \le \int_{w-1}^{w} \rho(t/2) dt 
\le \rho((w-1)/2) 
$$
and so $|F(w)| \to 0$ as $w\to \infty$.  This contradiction proves that 
$\sigma_{-}$ must change sign infinitely often, and 
completes our proof of Theorem 9.

\head Acknowledgements 
\endhead
We'd like to thank Richard Hall and
Hugh Montgomery for exciting our interest in this problem, 
Mark Watkins for the
question that inspired Theorem 9 and Malcolm Adams,
Kevin Clancey, Steven Finch, Jim Haglund, Seva Lev and Carl Pomerance 
for several useful remarks.  We are also grateful to the 
referee for some useful suggestions on style.

\enddocument